\documentclass[11pt]{amsart}
\usepackage{amssymb} 
\usepackage{amsmath} 
\usepackage{graphicx} 
\usepackage[english]{babel} 
\usepackage{epsfig}
\usepackage{enumerate}
\usepackage{color}
\usepackage{amsmath}

\usepackage[all]{xy}
\usepackage{mathrsfs}
\usepackage{graphics}
\usepackage{amsthm}
\usepackage{graphicx}

\usepackage[usenames,dvipsnames]{pstricks}
 \usepackage{pst-grad} 
\usepackage{pst-plot} 

\swapnumbers

\newenvironment{proo}{\begin{trivlist} \item{\emph{Proof.}}}
 {\hfill $\square$ \end{trivlist}}
\swapnumbers

\theoremstyle{definition}
\newtheorem{definition}[subsection]{Definition}
\newtheorem{example}[subsection]{Example}
\newtheorem{notation}[subsection]{Notation}
\newtheorem{remark}[subsection]{Remark}
\theoremstyle{plain}
\newtheorem{theorem}[subsection]{Theorem}
\newtheorem{proposition}[subsection]{Proposition}
\newtheorem{lemma}[subsection]{Lemma}
\newtheorem{corollary}[subsection]{Corollary}

\def\ZZ{{\mathbb{Z}}}
\def\K{{\mathbb K}}
\def\ot{\otimes}

\def\noi{\noindent}
\def\Dy{\mbox {\it Dyck}}
\def\D{\mathcal{D}}

\def\ov{\overline}

\def\up{\mathcal{UP}}
\def\dw{\mathcal{DW}}
\def\t{\times}
\def\ot{\otimes}
\def\D{\mathcal D}
\def\Dot2{\ov{\D_m^+\ot \D_m^+}}

\def\de{\Delta}
\def\lam{{\underline {\lambda}}}
\def\Dot2{\ov{\D_m^+\ot \D_m^+}}

\begin{document}

\author[D. L\'opez N., L.-F. Pr\'eville-Ratelle, M. Ronco]{Daniel L\'opez N., Louis-Fran\c cois Pr\'eville-Ratelle, Mar\'\i a Ronco}
\address{DL: Departamento de Matem\'aticas , Universidad de Concepci\'on\\
Concepci\'on, Chile}
\email{dlopezn@udec.cl}
\address{LFPR: Instituto de Matem\'atica y F\'\i sica\\
Campus Norte, Camino Lircay s/n,
\\Talca, Chile}
\email{lfprevilleratelle@inst-mat.utalca.cl}
\address{MOR: IMAFI, Universidad de Talca\\ Campus Norte, Avda. Lircay s/n\\ Talca, Chile}
\email{maria.ronco@inst-mat.utalca.cl}

\title{Algebraic structures defined on $m$-Dyck paths}


\thanks{\thanks{Our joint work was partially supported by the projects Fondecyt Postdoctorado 3140298 and Fondecyt Regular 1130939}
}


\begin{abstract} We introduce natural binary set-theoretical products on the set of all $m$-Dyck paths, which led us to define a non-symmetric algebraic operad $\Dy^m$. Our construction is closely related to the $m$-Tamari lattice, so the products defining $\Dy^m$ are given by intervals in this lattice. For $m=1$, we recover the notion of dendriform algebra introduced by J.-L. Loday in \cite{Lod}, and there exists a natural operad morphism from the operad ${\mbox {\it Ass}}$ of associative algebras into the operad  $\Dy^m$, consequently $\Dy ^m$ is a Hopf operad. We give a description of the coproduct in terms of $m$-Dyck paths in the last section. As an additional result, for any composition of $m+1\geq 2$ in $r+1$ parts, we get a functor from the category of $\Dy ^m$ algebras into the category of $\Dy ^r$ algebras. \end{abstract}

\maketitle

\section*{Introduction} \label{section:introduction}  

For $m\geq 1$, the $m$-Dyck paths are a particular family of lattice paths counted by Fuss-Catalan numbers, which are connected with the (bivariate) diagonal coinvariant spaces of the symmetric group. These representations are also called the Garsia-Haiman spaces, and they can be defined for an arbitrary number of sets of variables. Our work is motivated by the combinatorics of these spaces and by the Loday-Ronco Hopf algebra on binary trees.\\ 
\medskip

The Garsia-Haiman spaces have influenced the work of many combinatorialists in the past 20 years (see for instance \cite{Hai1994}, \cite{Hai2001}, \cite{Hai2002},  \cite{GarHag}), and they are still a very active area of research today (see \cite{BMChPR2013}, \cite{PRvie}, \cite{StuThoWil2015}) with many open problems. Note that the previous two lists of references are far from exhaustive. In particular we refer to the books of Bergeron (see \cite{BergBook}) and Haglund (\cite{HagBook}) for more explanations and references. Motivated by the combinatorics of the Garsia-Haiman spaces (see \cite{Hai1994}, \cite{Hai2001}, \cite{Hai2002}) and by an enumerative formula of Chapoton counting intervals in the Tamari lattice (see \cite{Cha2006}), F. Bergeron introduced the $m$-Tamari lattice, where the case $m=1$ is the usual Tamari lattice.  F. Bergeron and the second author (see \cite{BerPre}) showed that the trivariate diagonal coinvariant spaces are related to the intervals and the labelled intervals of the $m$-Tamari lattice. These labelled intervals are some generalizations of parking functions, where the latter is another family of combinatorial objects related with the (bivariate) Garsia-Haiman spaces. The $m$-Tamari lattice is the starting poing of our work.\\

In \cite{Lod}, J.-L. Loday introduced the notion of dendriform algebra and proved that the algebraic operad of dendriform algebras is naturally described on the vector space $\K[{\mathcal Y}_{\infty}]$ spanned by planar rooted binary trees. 
Dendriform algebras are associative algebras whose product splits as the sum of two binary operations. In many associative algebras already known in literature, as the algebras defined by shuffles (see \cite{EilMac} or \cite{MalReu}) and the Rota-Baxter algebras (see \cite{Agu}), the associative product comes from a dendriform structure. In \cite{LodRon}, J.-L. Loday and the third author, proved that any free dendriform algebra has a natural structure of bialgebra, which is described in terms of {\it admissible cuttings} of trees.  

The main goal of our work is to introduce a non-symmetric Hopf operad $\Dy ^m$ such that the space of $n$-ary operations of the theory is precisely the vector space $\K[\Dy _n^m]$, spanned by all the $m$-Dyck paths of size $n$, for any $m\geq 1$. When $m = 1$, we recover the operad of dendriform algebras. 

Given an $m$-Dyck path of size $n$, there is a unique way to color its down steps with elements of the set $\{ 1,\dots ,n\}$ in such a way that F. Bergeron\rq s covering relation consists in increasing the level of a down step without changing its color. This condition characterizes the order and is the key ingredient of our construction. The operad $\Dy ^m$ is spanned by $m+1$ binary operations $*_0,\dots , *_m$, which are given by intervals of F. Bergeron\rq s $m$-Tamari lattice.  For readers interested in algebraic operads, let us point out that the operads $\Dy^m$ are easily seen to be Koszul.

We also introduce the notion of $\Dy ^m$-bialgebra and described the coproduct on the vector space $\K[\Dy ^m]$, spanned by the set of $m$-Dyck paths, in terms of {\it admissible cuttings} of the Dyck path, which seem to be a particular case of the cuttings of rooted trees introduced by R. Grossman and R. Larson in \cite{GroLar}. 

For $m=1$, we know that the subspace of primitive elements of a dendriform bialgebra has a natural structure of brace algebra. For $m > 1$, the space of primitive elements of a $\Dy^m$ algebra is a brace algebra equipped with some additional structure. 
In a forthcoming work we describe the operads associated to the primitive elements of  $\Dy^m$ bialgebras.

Before giving a more precise description of the  contents of the manuscript, let us point out that in \cite{NovThi}, J.-C. Novelli and J.-Y. Thibon introduced the notion of $m$-permutations and defined the Sylvester congruence in this new context. These construction led them to define $m$-trees as the classes 
of $m$-permutations modulo the generalized Sylvester congruence. In a second work, see \cite{Nov}, J.-C. Novelli introduced 
the notion of $m$-dendriform algebra and showed that the vector space spanned by $m$-trees provide a natural description of this operad. 
Even if the dimension of the operad of $m$-dendriform algebras in degree $n$ is the number of $m$-Dyck path of size $n$ and both of them are generated by $m+1$ products, J.-C. Novelli\rq s operad  is different from $\Dy ^m$. In particular, our $\Dy ^m$ operad is defined by only two types of relations. A nice bijection between Dyck paths and $m$-trees still needs to be defined in order to compare both structures.
\medskip

{\it Contents} 

In the first section we recall some basic definitions and constructions of Dyck paths, needed in the sequel. 

In Section $2$ we introduce basic operations $\t _j$ on the set of $m$-Dyck paths, and the notion of coloring of a Dyck path. The basic constructions of this section are used in Section $3$ to define binary products $*_0, \dots ,*_m$ on the space $\K[\Dy ^m]$, spanned by the set of $m$-Dyck paths, and to prove the relations between them. 

In Section $4$, we show that the $\Dy ^m$ algebra structure on the space spanned by Dyck paths is related to the $m$-Tamari lattice by the formulas: 
$$P * _i Q = \sum _{P/_i Q \leq Z\leq P\backslash _iQ} Z,$$
for any pair of Dyck paths $P$ and $Q$ and any integer $0\leq i\leq m$. 

We introduce the formal definition of $\Dy ^m$ algebra in Section $5$, and prove that the space $\K[\Dy ^m]$, equipped with the products $*_i$ introduced in the previous section, is the free $\Dy ^m$ algebra spanned by one generator. As the operad of $\Dy ^m$ algebras is regular, the whole operad is described by the free object spanned by one generator, so the combinatorial properties of $m$-Dyck paths define completely the operad. We show that, given two non-negative integers $h <m$, there is a natural way to define for any composition ${\underline r}$ of $m$ in $h+1$ parts, an operad homomorphism ${\mathbb F}_{\underline r}$ from $\Dy ^h$ into $\Dy ^m$, which is compatible with the refinement of compositions.  In particular, any $\Dy ^m$ algebra has an underlying associative structure, which describes the Hopf operad structure of $\Dy ^m$. To end Section $5$ we prove that the image of a free $\Dy^m$ algebra under the functor ${\mathbb F}_{\underline r}$ is a free $\Dy^h$ algebra, for any composition ${\underline r}$ of $(m+1)$ in $(h+1)$ parts.

The last section is devoted to define the coproduct on Dyck paths in terms of admissible cuttings.
\bigskip

\subsection*{Acknowledgements}  The third author wants to thank M. Livernet for helpful discussions on lattice path operads which provided a first motivation for this work, as well as V. Dotsenko for his enthousiastic ideas, and F. Bassino for her kind interest. The economic support of the Universit\'e Paris 13 during a visit to the LAGA Laboratory in February 2015 was fundamental for the sketch of a reseach plan whose first result is the present manuscript. D. L\'opez N. and M. Ronco want to thank specially Prof. Antonio Laface for his interest and support. The second author would like to thank Luc Lapointe for many fruitful discussions.
\section*{Preliminaries}
\medskip

All the vector spaces considered in the present work are over $\K$, where $\K$ is a field. For any set $X$, we denote by $\K[X]$ the vector space spanned by $X$. For any $\K$-vector space $V$, we denote by $V^+:=\K\oplus V$ the augmented vector space.  The set of non-negative integers is denoted $\ZZ_+$.

\section{$m${-}Dyck paths}
\medskip

In the present section we introduce basic notions of the combinatorial and algebraic structures we shall need in the rest of the work. For more detailed constructions and the proofs of the results we refer to \cite{BerPre}, \cite{BouFusPre} and \cite{BMChPR2013}.

\begin{definition} \label{Dyckpath} For $m, n\geq 1$, an {\it $m$-Dyck path of size $n$} is a path on the real plan 
${\mathbb R}^2$, starting at $(0,0)$ and ending at $(2nm,0)$, consisting of up steps $(m,m)$ and down steps $(1,-1)$, 
which never goes below the $x$-axis. Note that the initial and terminal points of each step lean on $\ZZ_+^2$.\end{definition}

We denote by $\Dy _n^m$ the set of all $m$-Dyck paths of size $n$.

The number of elements of the set $\Dy _n^m$ is $d_{m,n}:= \frac 1{mn+1}\binom {(m+1)n}n$.

\begin{example}\label{example1} For $m = 2$, we get that \begin{enumerate}
\item the unique element of $\Dy_1^2$ is 

\medskip

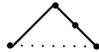
\scalebox{0.7} 
{
\begin{pspicture}(0,-0.47)(1.78,0.47)
\psdots[dotsize=0.12](0.06,-0.37)
\psdots[dotsize=0.12](0.92,0.39)
\psdots[dotsize=0.12](1.3,0.01)
\psdots[dotsize=0.12](1.7,-0.39)
\psline[linewidth=0.04cm](0.1,-0.39)(0.9,0.41)
\psline[linewidth=0.04cm](0.9,0.41)(1.7,-0.39)
\psline[linewidth=0.04cm,linestyle=dotted,dotsep=0.16cm](0.06,-0.37)(1.7,-0.39)
\end{pspicture} 
}
\medskip

\item the elements of $\Dy_2^2$ are

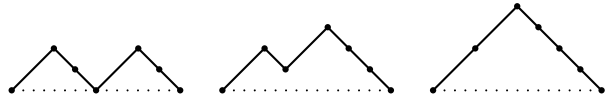
\scalebox{1.4} 
{
\begin{pspicture}(0,-0.44)(5.67,0.44)
\definecolor{color34}{rgb}{0.00392156862745098,0.00392156862745098,0.00392156862745098}
\psline[linewidth=0.02cm,linecolor=color34](0.03,-0.4)(0.43,0.0)
\psline[linewidth=0.02cm,linecolor=color34](0.43,0.0)(0.83,-0.4)
\psline[linewidth=0.02cm,linecolor=color34](0.83,-0.4)(1.23,0.0)
\psline[linewidth=0.02cm,linecolor=color34](1.23,0.0)(1.63,-0.4)
\psline[linewidth=0.02cm,linecolor=color34](2.03,-0.4)(2.43,0.0)
\psline[linewidth=0.02cm,linecolor=color34](2.43,0.0)(2.63,-0.2)
\psline[linewidth=0.02cm,linecolor=color34](2.63,-0.2)(3.03,0.2)
\psline[linewidth=0.02cm,linecolor=color34](3.03,0.2)(3.63,-0.4)
\psline[linewidth=0.02cm,linecolor=color34](4.03,-0.4)(4.83,0.4)
\psline[linewidth=0.02cm,linecolor=color34](4.83,0.4)(5.63,-0.4)
\psdots[dotsize=0.06,linecolor=color34](0.03,-0.4)
\psdots[dotsize=0.06,linecolor=color34](0.43,0.0)
\psdots[dotsize=0.06,linecolor=color34](0.63,-0.2)
\psdots[dotsize=0.06,linecolor=color34](0.83,-0.4)
\psdots[dotsize=0.06,linecolor=color34](1.23,0.0)
\psdots[dotsize=0.06,linecolor=color34](1.43,-0.2)
\psdots[dotsize=0.06,linecolor=color34](1.63,-0.4)

\psdots[dotsize=0.06,linecolor=color34](2.03,-0.4)
\psdots[dotsize=0.06,linecolor=color34](2.43,0.0)
\psdots[dotsize=0.06,linecolor=color34](2.63,-0.2)
\psdots[dotsize=0.06,linecolor=color34](3.03,0.2)
\psdots[dotsize=0.06,linecolor=color34](3.23,0.0)
\psdots[dotsize=0.06,linecolor=color34](3.43,-0.2)
\psdots[dotsize=0.06,linecolor=color34](3.63,-0.4)
\psline[linewidth=0.02cm,linestyle=dotted,dotsep=0.08cm](2.03,-0.4)(3.63,-0.4)
\psdots[dotsize=0.06,linecolor=color34](4.03,-0.4)
\psdots[dotsize=0.06,linecolor=color34](4.43,0.0)
\psdots[dotsize=0.06,linecolor=color34](4.83,0.4)
\psdots[dotsize=0.06,linecolor=color34](5.03,0.2)
\psdots[dotsize=0.06,linecolor=color34](5.23,0.0)
\psdots[dotsize=0.06,linecolor=color34](5.43,-0.2)
\psdots[dotsize=0.06,linecolor=color34](5.63,-0.4)
\psline[linewidth=0.02cm,linestyle=dotted,dotsep=0.08cm](4.03,-0.4)(5.63,-0.4)
\psline[linewidth=0.02cm,linestyle=dotted,dotsep=0.08cm](0.03,-0.4)(1.63,-0.4)
\end{pspicture} 
}

\end{enumerate}

\end{example}
\medskip

In order to define constructions on Dyck paths, we use a notation similar as the one employed by M. Bousquet-M\'elou, E. Fusy and the second author in \cite{BouFusPre}.
\medskip

\begin{notation} \label{updown} Let $P$ be an $m$-Dyck path. We denote by $\up(P)$ the set of up steps of $P$ and by $\dw(P)$ the set of down steps of $P$.\end{notation}

\begin{definition} \label{excursion} Let $u\in \up(P)$ be an up step of an $m$-Dyck path $P$, 
the {\it rank} of $u$ is $k$ if $u$ is the $k^{\rm th}$ up step of $P$, counting from left to right.

The shortest (translated) Dyck path which starts with $u$  is called the {\it excursion of $u$} in $P$, 
and is denoted $P_u$. The down step $w_u\in \dw(P)$ {\it matches} $u$ if it is the final step of the excursion of $u$ in $P$. 

Finally, a down step $d$ is at {\it level} $k$ if the last up step $u$ preceding $d$ has rank $k$. \end{definition}
\medskip

\begin{example} The up step $u$ in the path $P$ has rank $3$ and the down step $d$ matches it.

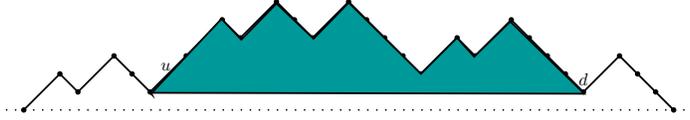
\scalebox{0.6} 
{
\begin{pspicture}(0,-1.28)(15.42,1.28)
\definecolor{color213b}{rgb}{0.0,0.6,0.6}
\psdots[dotsize=0.12](0.4,-1.2)
\psdots[dotsize=0.12](1.2,-0.4)
\psdots[dotsize=0.12](1.6,-0.8)
\psdots[dotsize=0.12](2.4,0.0)
\psdots[dotsize=0.12](2.8,-0.4)
\psdots[dotsize=0.12](3.2,-0.8)
\psdots[dotsize=0.12](4.0,0.0)
\psdots[dotsize=0.12](4.8,0.8)
\psdots[dotsize=0.12](5.2,0.4)
\psdots[dotsize=0.12](6.0,1.2)
\psdots[dotsize=0.12](6.4,0.8)
\psdots[dotsize=0.12](6.8,0.4)
\psdots[dotsize=0.12](7.6,1.2)
\psdots[dotsize=0.12](8.0,0.8)
\psdots[dotsize=0.12](8.4,0.4)
\psdots[dotsize=0.12](8.8,0.0)
\psdots[dotsize=0.12](9.2,-0.4)
\psdots[dotsize=0.12](10.0,0.4)
\psdots[dotsize=0.12](11.2,0.8)
\psdots[dotsize=0.12](10.4,0.0)
\psdots[dotsize=0.12](11.6,0.4)
\psdots[dotsize=0.12](12.0,0.0)
\psdots[dotsize=0.12](12.4,-0.4)
\psdots[dotsize=0.12](12.8,-0.8)
\psdots[dotsize=0.12](13.6,0.0)
\psdots[dotsize=0.12](14.0,-0.4)
\psdots[dotsize=0.12](14.4,-0.8)
\psdots[dotsize=0.12](14.8,-1.2)
\psline[linewidth=0.04cm,linestyle=dotted,dotsep=0.16cm](0.0,-1.2)(15.4,-1.2)
\psline[linewidth=0.04cm](0.4,-1.2)(1.2,-0.4)
\psline[linewidth=0.04cm](1.2,-0.4)(1.6,-0.8)
\psline[linewidth=0.04cm](1.6,-0.8)(2.4,0.0)
\psline[linewidth=0.04cm](2.4,0.0)(3.2,-0.8)
\psline[linewidth=0.04cm](3.2,-0.8)(4.8,0.8)
\psline[linewidth=0.04cm](4.8,0.8)(5.2,0.4)
\psline[linewidth=0.04cm](5.2,0.4)(6.0,1.2)
\psline[linewidth=0.04cm](6.0,1.2)(6.8,0.4)
\psline[linewidth=0.04cm](6.8,0.4)(7.6,1.2)
\psline[linewidth=0.04cm](7.6,1.2)(9.2,-0.4)
\psline[linewidth=0.04cm](9.2,-0.4)(10.0,0.4)
\psline[linewidth=0.04cm](10.0,0.4)(10.4,0.0)
\psline[linewidth=0.04cm](10.4,0.0)(11.2,0.8)
\psline[linewidth=0.04cm](11.2,0.8)(12.8,-0.8)
\psline[linewidth=0.04cm](12.8,-0.8)(13.6,0.0)
\psline[linewidth=0.04cm](13.6,0.0)(14.8,-1.2)
\usefont{T1}{ptm}{m}{n}
\rput(3.551455,-0.235){$u$}
\pspolygon[linewidth=0.04,fillstyle=solid,fillcolor=color213b](3.22,-0.84)(4.8,0.84)(4.78,0.82)(5.22,0.36)(6.0,1.16)(6.82,0.38)(7.6,1.18)(9.2,-0.4)(10.0,0.4)(10.38,-0.02)(11.18,0.76)(12.82,-0.84)(12.82,-0.84)(3.22,-0.82)(3.24,-0.86)
\usefont{T1}{ptm}{m}{n}
\rput(12.747485,-0.515){ $d$}
\end{pspicture} 
}

\end{example}
\medskip

Any up step $u$ in a Dyck path $P$ is determined by its rank, from now on we identify them, 
and denote the set of up steps of a Dyck path of size $n$ as $\up(P) =\{ 1, \dots ,n\}$.
\bigskip

\bigskip

\section{Operations on Dyck paths}
\medskip

We want to describe basic operations on Dyck paths that we need in the sequel.

\begin{notation}  \label{not:Ln(P)} For a path $P\in \Dy_n^m$ and an integer $1\leq k\leq n$, we denote by $\dw _k(P)$ the set of down steps of level $k$ of $P$ and by  $L_k(P)$ the number of elements of $\dw _k(P)$. When no confusion is possible, we shall denote the last term of the sequence $L_n(P)$ simply by $L(P)$.\end{notation}

Note that $0\leq \sum _{i=1}^j L_i(P)\leq mj$, for $1\leq j \leq n$. A Dyck path $P$ is uniquely determined by the sequence $(L_1(P),\dots , L_n(P))$.

\begin{definition} \label{suspension-concatenation}
Let $P$ and $Q$ be two $m$-Dyck paths of sizes $n_1$ and $n_2$, respectively. For $0\leq i\leq L(P)$, 
define the $i^{th}$-{\it concatenation} of $P$ and $Q$, denoted $P\t_i Q$, as the Dyck path of size $n = n_1+n_2$ 
obtained in the following way:
\begin{enumerate}
\item if $d_{1}^P, \dots ,d_{L(P)}^P$ denotes the ordered sequence of down steps of level $n_1$ of $P$, cut $P$ at the final vertex of $d_{L(P){-}i}^P$,
\item glue the initial point of $Q$ (translated) at the the final point of $d_{L(P){-}i}^P\in \dw(P)$,
\item glue the down steps $d_{L(P){-}i+1}^P, \dots ,d_{L(P)}^P$ at the end point of $Q$. \end{enumerate}\end{definition}
\medskip

\begin{example} Let $P$ and $Q$ be the $2$-Dyck paths

\medskip
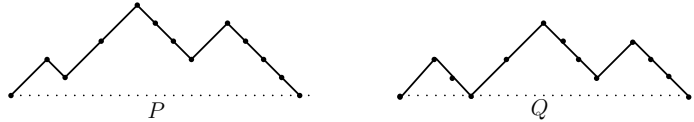
\scalebox{0.6} 
{
\begin{pspicture}(0,-1.1)(15.281015,1.1)
\psdots[dotsize=0.12](0.18101563,-0.98)
\psdots[dotsize=0.12](0.9810156,-0.18)
\psdots[dotsize=0.12](1.3810157,-0.58)
\psdots[dotsize=0.12](2.1810157,0.22)
\psdots[dotsize=0.12](2.9810157,1.02)
\psdots[dotsize=0.12](3.3810155,0.62)
\psdots[dotsize=0.12](3.7810156,0.22)
\psdots[dotsize=0.12](4.1810155,-0.18)
\psdots[dotsize=0.12](4.9810157,0.62)
\psdots[dotsize=0.12](5.381016,0.22)
\psdots[dotsize=0.12](5.7810154,-0.18)
\psdots[dotsize=0.12](6.1810155,-0.58)
\psdots[dotsize=0.12](6.5810156,-0.98)
\psline[linewidth=0.04cm,linestyle=dotted,dotsep=0.16cm](0.18101563,-0.98)(6.7810154,-0.98)
\psline[linewidth=0.04cm](0.18101563,-0.98)(0.9810156,-0.18)
\psline[linewidth=0.04cm](0.9810156,-0.18)(1.3810157,-0.58)
\psline[linewidth=0.04cm](1.3810157,-0.58)(2.9810157,1.02)
\psline[linewidth=0.04cm](2.9810157,1.02)(4.1810155,-0.18)
\psline[linewidth=0.04cm](4.1810155,-0.18)(4.9810157,0.62)
\psline[linewidth=0.04cm](4.9810157,0.62)(6.5810156,-0.98)
\usefont{T1}{ptm}{m}{n}
\rput(3.4,-1.3){\Large $P$}
\usefont{T1}{ptm}{m}{n}
\rput(11.9,-1.3){\Large $Q$}
\psdots[dotsize=0.12](8.801016,-1.0)
\psdots[dotsize=0.12](9.561016,-0.18)
\psdots[dotsize=0.12](10.381016,-1.0)
\psdots[dotsize=0.12](11.1610155,-0.18)
\psdots[dotsize=0.12](11.981015,0.62)
\psdots[dotsize=0.12](12.421016,0.22)
\psdots[dotsize=0.12](9.961016,-0.6)
\psdots[dotsize=0.12](12.761016,-0.18)
\psdots[dotsize=0.12](13.1610155,-0.6)
\psdots[dotsize=0.12](13.961016,0.2)
\psdots[dotsize=0.12](14.361015,-0.18)
\psdots[dotsize=0.12](14.761016,-0.56)
\psdots[dotsize=0.12](15.201015,-1.02)
\psdots[dotsize=0.12](8.801016,-1.02)
\psline[linewidth=0.04cm](8.821015,-1.02)(9.601016,-0.14)
\psline[linewidth=0.04cm](9.601016,-0.16)(10.421016,-1.04)
\psline[linewidth=0.04cm](10.381016,-0.98)(11.981015,0.62)
\psline[linewidth=0.04cm](11.981015,0.62)(13.181016,-0.58)
\psline[linewidth=0.04cm](13.181016,-0.58)(13.981015,0.22)
\psline[linewidth=0.04cm](13.981015,0.22)(15.181016,-0.98)
\psline[linewidth=0.04cm,linestyle=dotted,dotsep=0.16cm](8.781015,-0.98)(15.181016,-0.98)
\end{pspicture} 
}

\bigskip

the element $P\t_2 Q$ is the following $2$-Dyck path:

\medskip

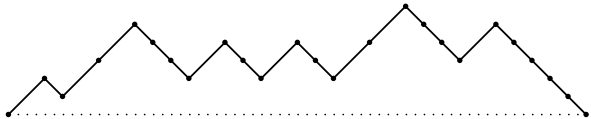
\scalebox{0.6} 
{
\begin{pspicture}(0,-1.28)(12.94,1.28)
\psdots[dotsize=0.12](0.06,-1.2)
\psdots[dotsize=0.12](0.86,-0.4)
\psdots[dotsize=0.12](1.26,-0.8)
\psdots[dotsize=0.12](2.06,0.0)
\psdots[dotsize=0.12](2.86,0.8)
\psdots[dotsize=0.12](3.26,0.4)
\psdots[dotsize=0.12](3.66,0.0)
\psdots[dotsize=0.12](4.06,-0.4)
\psdots[dotsize=0.12](4.86,0.4)
\psdots[dotsize=0.12](5.26,0.0)
\psdots[dotsize=0.12](5.66,-0.4)
\psdots[dotsize=0.12](6.46,0.4)
\psdots[dotsize=0.12](6.86,0.0)
\psdots[dotsize=0.12](7.26,-0.4)
\psdots[dotsize=0.12](8.06,0.4)
\psdots[dotsize=0.12](8.86,1.2)
\psdots[dotsize=0.12](9.26,0.8)
\psdots[dotsize=0.12](9.66,0.4)
\psdots[dotsize=0.12](10.06,0.0)
\psdots[dotsize=0.12](10.86,0.8)
\psdots[dotsize=0.12](11.26,0.4)
\psdots[dotsize=0.12](11.66,0.0)
\psdots[dotsize=0.12](12.06,-0.4)
\psdots[dotsize=0.12](12.46,-0.8)
\psdots[dotsize=0.12](12.86,-1.2)
\psline[linewidth=0.04cm,linestyle=dotted,dotsep=0.16cm](0.06,-1.2)(12.86,-1.2)
\psline[linewidth=0.04cm](0.06,-1.2)(0.86,-0.4)
\psline[linewidth=0.04cm](0.86,-0.4)(1.26,-0.8)
\psline[linewidth=0.04cm](1.26,-0.8)(2.86,0.8)
\psline[linewidth=0.04cm](2.86,0.8)(4.06,-0.4)
\psline[linewidth=0.04cm](4.06,-0.4)(4.86,0.4)
\psline[linewidth=0.04cm](4.86,0.4)(5.66,-0.4)
\psline[linewidth=0.04cm](5.66,-0.4)(6.46,0.4)
\psline[linewidth=0.04cm](6.46,0.4)(7.26,-0.4)
\psline[linewidth=0.04cm](7.26,-0.4)(8.86,1.2)
\psline[linewidth=0.04cm](8.86,1.2)(10.06,0.0)
\psline[linewidth=0.04cm](10.06,0.0)(10.86,0.8)
\psline[linewidth=0.04cm](10.86,0.8)(12.86,-1.2)
\end{pspicture} 
}

\end{example}
\medskip

\begin{definition} \label{defprime} An $m$-Dyck path $P$ is called {\it prime} if there does not exist a pair of $m$-Dyck paths of smaller size $Q$ and $R$ such that $P = Q\t_0 R$.\end{definition}
\medskip

\begin{remark} \label{remprime} For any $m$-Dyck path of size $n$ there exist a unique composition $(n_1,\dots ,n_r)$ of $n$ (with $n_i\geq 1$ for each $i$) and a unique family of prime Dyck paths $P_1\in \Dy_{n_1}^m,\ \dots , P_r\in \Dy_{n_r}^m$ such that $P = P_1\t_0 \dots \t_0 P_r$.

\end{remark}
\medskip

The proof of the following Lemma is immediate.

\begin{lemma} \label{lemma1} Let $P\in \Dy_{n_1}^m$ be a prime Dyck path and let $Q\in \Dy_{n_2}^m$ be another Dyck path. For any $1\leq j\leq L(P)$, the Dyck path $P\t _j Q$ is prime.\end{lemma}
\medskip

Define $\Dy _0^m :=\{ \bullet\}$, for $m\geq 1$. Any $m$-Dyck path $P$ of size $n$ may be written in two different ways, as:\begin{itemize}
\item $P= (((\rho _m\t _m P_0)\t _{m-1} P_1)\t _{m-2}\dots )\t _0 P_m,$ and 
\item $P = P_0\rq \t _0((( \rho _m\t _m P_1\rq)\t_{m-1}P_2\rq )\t _{m-2}\dots \t _2 P_{m-1}\rq)\t _1 P_m\rq ,$
for unique families of Dyck paths $P_0,\dots ,P_m$ and $P_0\rq, \dots ,P_m\rq$, with $P_j\in \Dy _{n_j}^m$ and $P_j\rq \in \Dy _{n_j\rq }^m$, $0\leq n_j, n_j\rq \leq n-1$ and ${\displaystyle \sum _{j=0}^m n_j =\sum _{j=0}^m n_j\rq = n-1}$.\end{itemize}
\medskip

For example, for the path $P$ in the example above, we get 
\medskip

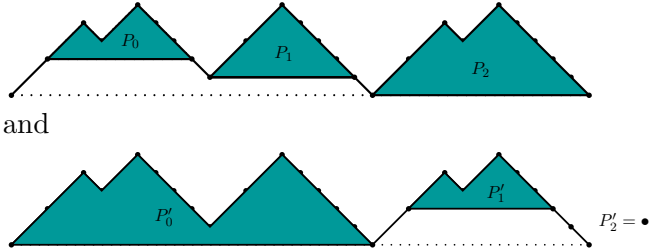
\scalebox{0.6} 
{
\begin{pspicture}(0,-1.084)(12.944,1.084)
\definecolor{color361b}{rgb}{0.0,0.6,0.6}
\psdots[dotsize=0.12](0.06,-1.0)
\psdots[dotsize=0.12](0.86,-0.2)
\psdots[dotsize=0.12](1.66,0.6)
\psdots[dotsize=0.12](2.06,0.2)
\psdots[dotsize=0.12](2.86,1.0)
\psdots[dotsize=0.12](3.66,0.2)
\psdots[dotsize=0.12](4.06,-0.2)
\psdots[dotsize=0.12](4.46,-0.6)
\psdots[dotsize=0.12](3.26,0.6)
\psdots[dotsize=0.12](5.26,0.2)
\psdots[dotsize=0.12](6.06,1.0)
\psdots[dotsize=0.12](6.46,0.6)
\psdots[dotsize=0.12](6.86,0.2)
\psdots[dotsize=0.12](7.26,-0.2)
\psdots[dotsize=0.12](7.66,-0.6)
\psdots[dotsize=0.12](8.06,-1.0)
\psdots[dotsize=0.12](8.86,-0.2)
\psdots[dotsize=0.12](9.66,0.6)
\psdots[dotsize=0.12](10.06,0.2)
\psdots[dotsize=0.12](10.86,1.0)
\psdots[dotsize=0.12](11.26,0.6)
\psdots[dotsize=0.12](11.66,0.2)
\psdots[dotsize=0.12](12.06,-0.2)
\psdots[dotsize=0.12](12.46,-0.6)
\psdots[dotsize=0.12](12.86,-1.0)
\psline[linewidth=0.048cm,fillcolor=color361b,linestyle=dotted,dotsep=0.16cm](0.06,-1.0)(12.86,-1.0)
\pspolygon[linewidth=0.048,linestyle=dotted,dotsep=0.16cm,fillstyle=solid,fillcolor=color361b](2.06,0.2)(1.66,0.6)(0.86,-0.2)(4.06,-0.2)(2.86,1.0)(2.86,1.0)
\pspolygon[linewidth=0.048,linestyle=dotted,dotsep=0.16cm,fillstyle=solid,fillcolor=color361b](2.86,1.0)(2.86,0.8)(2.86,1.0)
\pspolygon[linewidth=0.048,linestyle=dotted,dotsep=0.16cm,fillstyle=solid,fillcolor=color361b](2.86,1.0)(2.06,0.2)(2.06,0.2)(2.06,0.2)(2.06,0.2)
\pspolygon[linewidth=0.048,fillstyle=solid,fillcolor=color361b](2.06,0.2)(1.66,0.6)(0.86,-0.2)(4.06,-0.2)(2.86,1.0)
\pspolygon[linewidth=0.048,fillstyle=solid,fillcolor=color361b](4.46,-0.6)(6.06,1.0)(7.66,-0.6)(7.66,-0.6)(7.66,-0.6)(4.46,-0.6)(4.46,-0.6)(7.66,-0.6)(6.06,-0.6)
\psline[linewidth=0.048cm,fillcolor=color361b](0.86,-0.2)(0.06,-1.0)
\psline[linewidth=0.048cm,fillcolor=color361b](4.06,-0.2)(4.46,-0.6)
\psline[linewidth=0.048cm,fillcolor=color361b](7.66,-0.6)(8.06,-1.0)
\pspolygon[linewidth=0.048,fillstyle=solid,fillcolor=color361b](8.06,-1.0)(9.66,0.6)(10.06,0.2)(10.86,1.0)(12.86,-1.0)(8.06,-1.0)
\usefont{T1}{ptm}{m}{n}
\rput(2.7,0.205){$P_0$}
\usefont{T1}{ptm}{m}{n}
\rput(6.1,-0.075){$P_1$}
\usefont{T1}{ptm}{m}{n}
\rput(10.46,-0.435){$P_2$}
\end{pspicture} 
}

and
\medskip

\scalebox{0.6} 
{
\begin{pspicture}(0,-1.08)(15.001895,1.08)
\definecolor{color43b}{rgb}{0.0,0.6,0.6}
\psdots[dotsize=0.12](0.06,-1.0)
\psdots[dotsize=0.12](0.86,-0.2)
\psdots[dotsize=0.12](1.66,0.6)
\psdots[dotsize=0.12](2.06,0.2)
\psdots[dotsize=0.12](2.86,1.0)
\psdots[dotsize=0.12](3.26,0.6)
\psdots[dotsize=0.12](3.66,0.2)
\psdots[dotsize=0.12](4.06,-0.2)
\psdots[dotsize=0.12](4.46,-0.6)
\psdots[dotsize=0.12](5.26,0.2)
\psdots[dotsize=0.12](6.06,1.0)
\psdots[dotsize=0.12](6.46,0.6)
\psdots[dotsize=0.12](6.86,0.2)
\psdots[dotsize=0.12](7.26,-0.2)
\psdots[dotsize=0.12](7.66,-0.6)
\psdots[dotsize=0.12](8.06,-1.0)
\psdots[dotsize=0.12](8.86,-0.2)
\psdots[dotsize=0.12](9.66,0.6)
\psdots[dotsize=0.12](10.06,0.2)
\psdots[dotsize=0.12](10.86,1.0)
\psdots[dotsize=0.12](11.26,0.6)
\psdots[dotsize=0.12](11.66,0.2)
\psdots[dotsize=0.12](12.06,-0.2)
\psdots[dotsize=0.12](12.46,-0.6)
\psdots[dotsize=0.12](12.86,-1.0)
\psline[linewidth=0.04cm,linestyle=dotted,dotsep=0.16cm](0.06,-1.0)(12.86,-1.0)
\psline[linewidth=0.04cm](0.06,-1.0)(1.66,0.6)
\psline[linewidth=0.04cm](1.66,0.6)(2.06,0.2)
\psline[linewidth=0.04cm](2.06,0.2)(2.86,1.0)
\psline[linewidth=0.04cm](2.86,1.0)(4.46,-0.6)
\psline[linewidth=0.04cm](4.46,-0.6)(6.06,1.0)
\psline[linewidth=0.04cm](6.06,1.0)(8.06,-1.0)
\psline[linewidth=0.04cm](8.06,-1.0)(9.66,0.6)
\psline[linewidth=0.04cm](9.66,0.6)(10.06,0.2)
\psline[linewidth=0.04cm](10.06,0.2)(10.86,1.0)
\psline[linewidth=0.04cm](10.86,1.0)(12.86,-1.0)
\pspolygon[linewidth=0.04,fillstyle=solid,fillcolor=color43b](12.06,-0.2)(11.86,-0.2)(8.86,-0.2)(9.66,0.6)(10.06,0.2)(10.86,1.0)
\pspolygon[linewidth=0.04,fillstyle=solid,fillcolor=color43b](8.06,-1.0)(0.06,-1.0)(1.66,0.6)(2.06,0.2)(2.86,1.0)(4.46,-0.6)(6.06,1.0)(8.06,-1.0)
\usefont{T1}{ptm}{m}{n}
\rput(3.4414551,-0.375){$P_0\rq$}
\usefont{T1}{ptm}{m}{n}
\rput(10.8014555,0.145){$P_1\rq$}
\usefont{T1}{ptm}{m}{n}
\rput(13.641455,-0.495){$P_2\rq = \bullet$}
\end{pspicture} 
}
\bigskip

\begin{notation} \label{notrecursive} For any Dyck path $P=$
$$(((\rho _m\t _m P_0)\t _{m-1} P_1)\t _{m-2}\dots )\t _0 P_m = P_0\rq \t _0 ((\rho _m\t _m P_1\rq )\t _{m-1} \dots \t _2P_{m-1}\rq)\t _1P_m\rq,$$
 we denote it by $P = \bigvee_d (P_0,\dots ,P_m) = \bigvee _u(P_0\rq ,\dots ,P_m\rq )$.\end{notation}
 
 Note that, the $P_i$\rq s and the $P_i\rq $ may be just the point $\bullet$.
\medskip

\begin{remark}
\label{remark: generating series for Dyck paths}
Let $d_m(x)$ be the generating series of $K[\Dy^m]^+$, that is, 
\begin{center}
 $d_m(x) : =\displaystyle \sum_{n\geq 0}d_{m,n}x^n$
\end{center}
where $d_{m,n}$ is the dimension of $K[\Dy^m_n]$ and $d_{m,0}=1$. The preceding discussion implies that the series $d_m(x)$ 
satisfies the equation $x\cdot d_m(x)^{m+1} = d_m(x)-1$.
\end{remark}
\medskip

\begin{remark}\label{remnotations} A Dyck path $P$ is prime if, and only if, $P$ is of the form $P = \bigvee_d (P_0,\dots ,P_m) = \bigvee _u(P_0\rq ,\dots ,P_m\rq )$, with $P_r = P_0\rq = \bullet $.\end{remark}
\medskip

\begin{definition}\label{defcoloring} Let $P$ be a $m$-Dyck path of size $n$. The {\it standard coloring} of $P$ is a map $\alpha _P$ from the set of down steps $\dw(P)$ to the set $\{1,\dots ,n\}$, described recursively as follows: \begin{enumerate}
\item For $P = \rho _m\in \Dy_1^m$, $\alpha _{\rho _m}$ is the constant function $1$.
\item For $P = \bigvee_d (P_0,\dots ,P_m)$, with $P_j\in \Dy _{n_j}^m$, the set of down steps of $P$ is the disjoint union 
$$\dw (P) = \{ 1,\dots ,m\}\coprod \dw(P_0) \coprod \dots \coprod  \dw(P_m),$$ where the first subset $\{1,\dots ,m\}$ corresponds to the down steps of $\rho _m$.

The map $\alpha _P$ is defined by:
$$\alpha _P(d) = \begin{cases} 1,&{\rm for}\ d\in \{1,\dots ,m\},\\
\alpha _{P_j}(d)+n_0+\dots +n_{j-1}+1,&{\rm for}\ d\in \dw(P_j),\end{cases}$$
where $0\leq j\leq m$.
\end{enumerate}\end{definition}
\medskip

In our last example, we get the following coloring for $P$:
\medskip

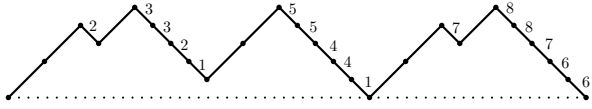
\scalebox{0.6} 
{
\begin{pspicture}(0,-1.124)(13.16,1.14)
\definecolor{color361b}{rgb}{0.0,0.6,0.6}
\psdots[dotsize=0.12](0.06,-1.04)
\psdots[dotsize=0.12](0.86,-0.24)
\psdots[dotsize=0.12](1.66,0.56)
\psdots[dotsize=0.12](2.06,0.16)
\psdots[dotsize=0.12](2.86,0.96)
\psdots[dotsize=0.12](3.66,0.16)
\psdots[dotsize=0.12](4.06,-0.24)
\psdots[dotsize=0.12](4.46,-0.64)
\psdots[dotsize=0.12](3.26,0.56)
\psdots[dotsize=0.12](5.26,0.16)
\psdots[dotsize=0.12](6.06,0.96)
\psdots[dotsize=0.12](6.46,0.56)
\psdots[dotsize=0.12](6.86,0.16)
\psdots[dotsize=0.12](7.26,-0.24)
\psdots[dotsize=0.12](7.66,-0.64)
\psdots[dotsize=0.12](8.06,-1.04)
\psdots[dotsize=0.12](8.86,-0.24)
\psdots[dotsize=0.12](9.66,0.56)
\psdots[dotsize=0.12](10.06,0.16)
\psdots[dotsize=0.12](10.86,0.96)
\psdots[dotsize=0.12](11.26,0.56)
\psdots[dotsize=0.12](11.66,0.16)
\psdots[dotsize=0.12](12.06,-0.24)
\psdots[dotsize=0.12](12.46,-0.64)
\psdots[dotsize=0.12](12.86,-1.04)
\psline[linewidth=0.048cm,fillcolor=color361b,linestyle=dotted,dotsep=0.16cm](0.06,-1.04)(12.86,-1.04)
\psline[linewidth=0.048cm,fillcolor=color361b](0.86,-0.24)(0.06,-1.04)
\psline[linewidth=0.048cm,fillcolor=color361b](4.06,-0.24)(4.46,-0.64)
\psline[linewidth=0.048cm,fillcolor=color361b](7.66,-0.64)(8.06,-1.04)
\psline[linewidth=0.048cm,fillcolor=color361b](0.86,-0.24)(1.66,0.56)
\psline[linewidth=0.048cm,fillcolor=color361b](1.66,0.56)(2.06,0.16)
\psline[linewidth=0.048cm,fillcolor=color361b](2.06,0.16)(2.86,0.96)
\psline[linewidth=0.048cm,fillcolor=color361b](2.86,0.96)(4.06,-0.24)
\psline[linewidth=0.048cm,fillcolor=color361b](4.46,-0.64)(6.06,0.96)
\psline[linewidth=0.048cm,fillcolor=color361b](6.06,0.96)(7.66,-0.64)
\psline[linewidth=0.048cm,fillcolor=color361b](8.06,-1.04)(9.66,0.56)
\psline[linewidth=0.048cm,fillcolor=color361b](9.66,0.56)(10.06,0.16)
\psline[linewidth=0.048cm,fillcolor=color361b](10.06,0.16)(10.86,0.96)
\psline[linewidth=0.048cm,fillcolor=color361b](10.86,0.96)(12.86,-1.04)
\usefont{T1}{ptm}{m}{n}
\rput(8.04,-0.695){\small $1$}
\usefont{T1}{ptm}{m}{n}
\rput(4.36,-0.315){\small $1$}
\usefont{T1}{ptm}{m}{n}
\rput(3.96,0.125){\small $2$}
\usefont{T1}{ptm}{m}{n}
\rput(1.94,0.545){\small $2$}
\usefont{T1}{ptm}{m}{n}
\rput(3.18,0.925){\small $3$}
\usefont{T1}{ptm}{m}{n}
\rput(3.58,0.525){\small $3$}
\usefont{T1}{ptm}{m}{n}
\rput(7.26,0.105){\small $4$}
\usefont{T1}{ptm}{m}{n}
\rput(7.58,-0.255){\small $4$}
\usefont{T1}{ptm}{m}{n}
\rput(6.82,0.525){\small $5$}
\usefont{T1}{ptm}{m}{n}
\rput(12.4,-0.275){\small $6$}
\usefont{T1}{ptm}{m}{n}
\rput(6.36,0.925){\small $5$}
\usefont{T1}{ptm}{m}{n}
\rput(12.86,-0.695){\small $6$}
\usefont{T1}{ptm}{m}{n}
\rput(9.98,0.505){\small $7$}
\usefont{T1}{ptm}{m}{n}
\rput(12.04,0.145){\small $7$}
\usefont{T1}{ptm}{m}{n}
\rput(11.18,0.945){\small $8$}
\usefont{T1}{ptm}{m}{n}
\rput(11.6,0.525){\small $8$}
\end{pspicture} 
}

\medskip

\begin{notation} \label{notcolor} For any path $P\in \Dy _n^m$ and any $1\leq k\leq n$, 
we denote by $\omega _k(P)$ the word $\omega _k^P:= \alpha _P(d_{k1})\dots  \alpha _P(d_{kL_k(P)})$, 
which is the image under $\alpha _P$ of the sequence of level $k$ down steps of $P$ (from left to right).\end{notation} 
\medskip

\begin{remark} \label{remcoloring} Let $P$ be an $m$-Dyck path of size $n$.
\begin{enumerate} 
\item For any down step $d\in \dw(P)$, the color of $d$ coincides with the rank of the up step $u\in \up(P)$ 
which is the first intersection of the  horizontal half-line beginning at the middle point of $d$ and going to the left side 
with the Dyck path $P$. In the example
\medskip

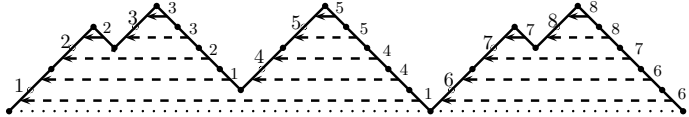
\scalebox{0.7} 
{
\begin{pspicture}(0,-1.114)(13.14,1.13)
\definecolor{color361b}{rgb}{0.0,0.6,0.6}
\psdots[dotsize=0.12](0.08,-1.03)
\psdots[dotsize=0.12](0.88,-0.23)
\psdots[dotsize=0.12](1.68,0.57)
\psdots[dotsize=0.12](2.08,0.17)
\psdots[dotsize=0.12](2.88,0.97)
\psdots[dotsize=0.12](3.68,0.17)
\psdots[dotsize=0.12](4.08,-0.23)
\psdots[dotsize=0.12](4.48,-0.63)
\psdots[dotsize=0.12](3.28,0.57)
\psdots[dotsize=0.12](5.28,0.17)
\psdots[dotsize=0.12](6.08,0.97)
\psdots[dotsize=0.12](6.48,0.57)
\psdots[dotsize=0.12](6.88,0.17)
\psdots[dotsize=0.12](7.28,-0.23)
\psdots[dotsize=0.12](7.68,-0.63)
\psdots[dotsize=0.12](8.08,-1.03)
\psdots[dotsize=0.12](8.88,-0.23)
\psdots[dotsize=0.12](9.68,0.57)
\psdots[dotsize=0.12](10.08,0.17)
\psdots[dotsize=0.12](10.88,0.97)
\psdots[dotsize=0.12](11.28,0.57)
\psdots[dotsize=0.12](11.68,0.17)
\psdots[dotsize=0.12](12.08,-0.23)
\psdots[dotsize=0.12](12.48,-0.63)
\psdots[dotsize=0.12](12.88,-1.03)
\psline[linewidth=0.048cm,fillcolor=color361b,linestyle=dotted,dotsep=0.16cm](0.08,-1.03)(12.88,-1.03)
\psdots[dotsize=0.12,fillstyle=solid,dotstyle=o](0.48,-0.63)
\psdots[dotsize=0.12,fillstyle=solid,dotstyle=o](1.28,0.17)
\psdots[dotsize=0.12,fillstyle=solid,dotstyle=o](2.48,0.57)
\psdots[dotsize=0.12,fillstyle=solid,dotstyle=o](4.88,-0.23)
\psdots[dotsize=0.12,fillstyle=solid,dotstyle=o](5.68,0.57)
\psdots[dotsize=0.12,fillstyle=solid,dotstyle=o](9.28,0.17)
\psdots[dotsize=0.12,fillstyle=solid,dotstyle=o](8.48,-0.63)
\psdots[dotsize=0.12,fillstyle=solid,dotstyle=o](10.48,0.57)
\psline[linewidth=0.048cm,fillcolor=color361b](0.88,-0.23)(0.08,-1.03)
\psline[linewidth=0.048cm,fillcolor=color361b](4.08,-0.23)(4.48,-0.63)
\psline[linewidth=0.048cm,fillcolor=color361b](7.68,-0.63)(8.08,-1.03)
\psline[linewidth=0.048cm,fillcolor=color361b](0.88,-0.23)(1.68,0.57)
\psline[linewidth=0.048cm,fillcolor=color361b](1.68,0.57)(2.08,0.17)
\psline[linewidth=0.048cm,fillcolor=color361b](2.08,0.17)(2.88,0.97)
\psline[linewidth=0.048cm,fillcolor=color361b](2.88,0.97)(4.08,-0.23)
\psline[linewidth=0.048cm,fillcolor=color361b](4.48,-0.63)(6.08,0.97)
\psline[linewidth=0.048cm,fillcolor=color361b](6.08,0.97)(7.68,-0.63)
\psline[linewidth=0.048cm,fillcolor=color361b](8.08,-1.03)(9.68,0.57)
\psline[linewidth=0.048cm,fillcolor=color361b](9.68,0.57)(10.08,0.17)
\psline[linewidth=0.048cm,fillcolor=color361b](10.08,0.17)(10.88,0.97)
\psline[linewidth=0.048cm,fillcolor=color361b](10.88,0.97)(12.88,-1.03)
\usefont{T1}{ptm}{m}{n}
\rput(8.04,-0.69){\footnotesize $1$}
\usefont{T1}{ptm}{m}{n}
\rput(4.36,-0.31){\footnotesize $1$}
\usefont{T1}{ptm}{m}{n}
\rput(3.96,0.13){\footnotesize $2$}
\usefont{T1}{ptm}{m}{n}
\rput(1.94,0.55){\footnotesize $2$}
\usefont{T1}{ptm}{m}{n}
\rput(3.18,0.93){\footnotesize $3$}
\usefont{T1}{ptm}{m}{n}
\rput(3.58,0.53){\footnotesize $3$}
\usefont{T1}{ptm}{m}{n}
\rput(7.26,0.11){\footnotesize $4$}
\usefont{T1}{ptm}{m}{n}
\rput(7.58,-0.25){\footnotesize $4$}
\usefont{T1}{ptm}{m}{n}
\rput(6.82,0.53){\footnotesize $5$}
\usefont{T1}{ptm}{m}{n}
\rput(12.4,-0.27){\footnotesize $6$}
\usefont{T1}{ptm}{m}{n}
\rput(6.36,0.93){\footnotesize $5$}
\usefont{T1}{ptm}{m}{n}
\rput(12.86,-0.69){\footnotesize $6$}
\usefont{T1}{ptm}{m}{n}
\rput(9.98,0.51){\footnotesize $7$}
\usefont{T1}{ptm}{m}{n}
\rput(12.04,0.15){\footnotesize $7$}
\usefont{T1}{ptm}{m}{n}
\rput(11.18,0.95){\footnotesize $8$}
\usefont{T1}{ptm}{m}{n}
\rput(11.6,0.53){\footnotesize $8$}
\usefont{T1}{ptm}{m}{n}
\rput(0.27,-0.545){$1$}
\usefont{T1}{ptm}{m}{n}
\rput(1.15,0.295){$2$}
\usefont{T1}{ptm}{m}{n}
\rput(2.43,0.755){$3$}
\usefont{T1}{ptm}{m}{n}
\rput(4.83,0.015){$4$}
\usefont{T1}{ptm}{m}{n}
\rput(5.53,0.655){$5$}
\usefont{T1}{ptm}{m}{n}
\rput(8.43,-0.425){$6$}
\usefont{T1}{ptm}{m}{n}
\rput(9.17,0.275){$7$}
\usefont{T1}{ptm}{m}{n}
\rput(10.41,0.695){$8$}
\psline[linewidth=0.048cm,fillcolor=color361b,linestyle=dashed,dash=0.16cm 0.16cm](2.1,0.11)(2.02,0.13)
\psline[linewidth=0.048cm,fillcolor=color361b,linestyle=dashed,dash=0.16cm 0.16cm,arrowsize=0.05291667cm 2.0,arrowlength=1.4,arrowinset=0.4]{->}(4.28,-0.43)(0.68,-0.43)
\psline[linewidth=0.048cm,fillcolor=color361b,linestyle=dashed,dash=0.16cm 0.16cm,arrowsize=0.05291667cm 2.0,arrowlength=1.4,arrowinset=0.4]{->}(3.88,-0.03)(1.08,-0.03)
\psline[linewidth=0.048cm,fillcolor=color361b,linestyle=dashed,dash=0.16cm 0.16cm,arrowsize=0.05291667cm 2.0,arrowlength=1.4,arrowinset=0.4]{->}(3.48,0.37)(2.28,0.37)
\psline[linewidth=0.048cm,fillcolor=color361b,linestyle=dashed,dash=0.16cm 0.16cm,arrowsize=0.05291667cm 2.0,arrowlength=1.4,arrowinset=0.4]{->}(3.08,0.77)(2.68,0.77)
\psline[linewidth=0.048cm,fillcolor=color361b,linestyle=dashed,dash=0.16cm 0.16cm,arrowsize=0.05291667cm 2.0,arrowlength=1.4,arrowinset=0.4]{->}(1.88,0.37)(1.48,0.37)
\psline[linewidth=0.048cm,fillcolor=color361b,linestyle=dashed,dash=0.16cm 0.16cm,arrowsize=0.05291667cm 2.0,arrowlength=1.4,arrowinset=0.4]{->}(7.48,-0.43)(4.68,-0.43)
\psline[linewidth=0.048cm,fillcolor=color361b,linestyle=dashed,dash=0.16cm 0.16cm,arrowsize=0.05291667cm 2.0,arrowlength=1.4,arrowinset=0.4]{->}(7.08,-0.03)(5.08,-0.03)
\psline[linewidth=0.048cm,fillcolor=color361b,linestyle=dashed,dash=0.16cm 0.16cm,arrowsize=0.05291667cm 2.0,arrowlength=1.4,arrowinset=0.4]{->}(6.68,0.37)(5.48,0.37)
\psline[linewidth=0.048cm,fillcolor=color361b,linestyle=dashed,dash=0.16cm 0.16cm,arrowsize=0.05291667cm 2.0,arrowlength=1.4,arrowinset=0.4]{->}(6.28,0.77)(5.88,0.77)
\psline[linewidth=0.048cm,fillcolor=color361b,linestyle=dashed,dash=0.16cm 0.16cm,arrowsize=0.05291667cm 2.0,arrowlength=1.4,arrowinset=0.4]{->}(12.68,-0.83)(8.28,-0.83)
\psline[linewidth=0.048cm,fillcolor=color361b,linestyle=dashed,dash=0.16cm 0.16cm,arrowsize=0.05291667cm 2.0,arrowlength=1.4,arrowinset=0.4]{->}(12.28,-0.43)(8.68,-0.43)
\psline[linewidth=0.048cm,fillcolor=color361b,linestyle=dashed,dash=0.16cm 0.16cm,arrowsize=0.05291667cm 2.0,arrowlength=1.4,arrowinset=0.4]{->}(11.88,-0.03)(9.08,-0.03)
\psline[linewidth=0.048cm,fillcolor=color361b,linestyle=dashed,dash=0.16cm 0.16cm,arrowsize=0.05291667cm 2.0,arrowlength=1.4,arrowinset=0.4]{->}(11.48,0.37)(10.28,0.37)
\psline[linewidth=0.048cm,fillcolor=color361b,linestyle=dashed,dash=0.16cm 0.16cm,arrowsize=0.05291667cm 2.0,arrowlength=1.4,arrowinset=0.4]{->}(11.08,0.77)(10.68,0.77)
\psline[linewidth=0.048cm,fillcolor=color361b,linestyle=dashed,dash=0.16cm 0.16cm,arrowsize=0.05291667cm 2.0,arrowlength=1.4,arrowinset=0.4]{->}(7.88,-0.83)(0.28,-0.83)
\psline[linewidth=0.048cm,fillcolor=color361b,linestyle=dashed,dash=0.16cm 0.16cm,arrowsize=0.05291667cm 2.0,arrowlength=1.4,arrowinset=0.4]{->}(9.88,0.37)(9.48,0.37)
\end{pspicture} 
}
\medskip

\item We have that $\vert \alpha _P^{-1}(i)\vert = m$, for any $1\leq i\leq n$.
\item For a fixed $1\leq k\leq n$, the word $\omega _k^P= \alpha _P(d_{1}^P)\dots  \alpha _P(d_{L_k(P)}^P)$ 
is decreasing for the usual order of the natural numbers. Moreover, the first $m$ digits of $\omega _n^P$ are $n$\rq s. 
\item If $Q$ is another $m$-Dyck path, then $\dw (P\t _i Q) = \dw (P)\coprod \dw (Q)$, and $\alpha _{P\t_iQ}$ is described by:
$$\alpha _{P\t_iQ}(d)= \begin{cases} \alpha _P(d),&{\rm for\ any}\ d\ {\rm which\ belongs\ initially\ to}\ P,\\
\alpha _Q(d) + n,&{\rm for\ any}\ d\ {\rm which\ belongs\ initially\ to}\ Q,\end{cases}$$
for any $0\leq i\leq L(P)$.
\end{enumerate}
\end{remark}

\bigskip

\bigskip

\section{Products on $m$-Dyck paths}
\medskip

\begin{definition} \label{defpartition} For any positive integer $n$, a {\it weak composition} of $n$ with $r+1$ parts is an ordered collection 
of non-negative integers ${\underline {\lambda}}= (\lambda_0,\dots ,\lambda _r)$ such that ${\displaystyle \sum _{i=0}^r\lambda _i = n}$. We say that the {\it length} of $\lam$ is $r+1$.\end{definition}
\medskip

\begin{notation} \label{notLambda} Given an $m$-Dyck path $P$ of size $n$, the set of all weak compositions of $L(P)$ of length $r+1$ is denoted $\Lambda _r(P)$.\end{notation}
\medskip

Let $P\in \Dy_{n_1}^m$ and $Q= Q_1\t_0 \dots \t_0 Q_r\in \Dy _{n_2}^m$ be two Dyck paths, where 
 $Q_j\in \Dy^m$ is prime, for $1\leq j\leq r$. 

Suppose that $\lam =(\lambda _0,\dots ,\lambda _r)$ is a weak composition of $L(P)$. 
Define a Dyck path $P*_{\lam} Q$ of size $n_1+n_2$ by the formula:
$$P*_{\lam} Q := ((((P\t _{\lambda _1+\dots +\lambda _r}Q_1)\t _{\lambda _2+\dots + \lambda _r} Q_2)\t _{\lambda _3+\dots +\lambda _r}\dots )\t _{\lambda _r}Q_r).$$

The product $*_{\lam}$ just divides the ordered set $\dw_{n_1}(P)$ of down steps of level $n_1$ of $P$ and glue, in order, the $i^{th}$ piece at the end of the path $Q_i$. If $\lambda _0 > 0$, the first $\lambda _0$ steps of $\dw_{n_1}(P)$ remain at the end of $P$.

\begin{example} Let $P = (2, 3, 1, 6)$ be a path in $\Dy _4^3$ and let $Q = (1,4,4,3, 2,3,4)$ be a $3$-Dyck path of size $7$, note that $Q = (1,4,4)\t _0 (3)\t_0 (2,3,4)$.
\medskip

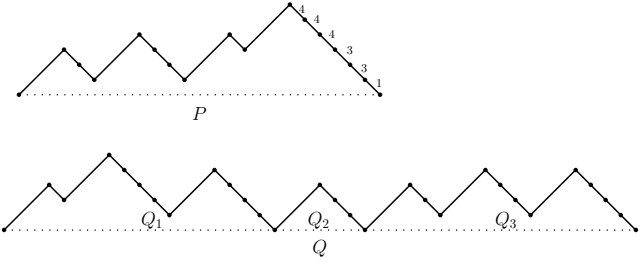
\scalebox{0.5} 
{
\begin{pspicture}(0,-3.08)(16.94,3.08)
\psdots[dotsize=0.12](0.46,0.6)
\psdots[dotsize=0.12](1.66,1.8)
\psdots[dotsize=0.12](2.06,1.4)
\psdots[dotsize=0.12](2.46,1.0)
\psdots[dotsize=0.12](3.66,2.2)
\psdots[dotsize=0.12](4.06,1.8)
\psdots[dotsize=0.12](4.46,1.4)
\psdots[dotsize=0.12](4.86,1.0)
\psdots[dotsize=0.12](6.06,2.2)
\psdots[dotsize=0.12](6.46,1.8)
\psdots[dotsize=0.12](7.66,3.0)
\psdots[dotsize=0.12](8.06,2.6)
\psdots[dotsize=0.12](8.46,2.2)
\psdots[dotsize=0.12](8.86,1.8)
\psdots[dotsize=0.12](9.26,1.4)
\psdots[dotsize=0.12](9.66,1.0)
\psdots[dotsize=0.12](10.06,0.6)
\psline[linewidth=0.04cm](0.46,0.6)(1.66,1.8)
\psline[linewidth=0.04cm](1.66,1.8)(2.46,1.0)
\psline[linewidth=0.04cm](2.46,1.0)(3.66,2.2)
\psline[linewidth=0.04cm](3.66,2.2)(4.86,1.0)
\psline[linewidth=0.04cm](4.86,1.0)(6.06,2.2)
\psline[linewidth=0.04cm](6.06,2.2)(6.46,1.8)
\psline[linewidth=0.04cm](6.46,1.8)(7.66,3.0)
\psline[linewidth=0.04cm](7.66,3.0)(10.06,0.6)
\psline[linewidth=0.04cm,linestyle=dotted,dotsep=0.16cm](0.46,0.6)(10.06,0.6)
\usefont{T1}{ptm}{m}{n}
\rput(10.04,0.92){\footnotesize $1$}
\usefont{T1}{ptm}{m}{n}
\rput(7.98,2.88){\footnotesize $4$}
\usefont{T1}{ptm}{m}{n}
\rput(8.38,2.6){\footnotesize $4$}
\usefont{T1}{ptm}{m}{n}
\rput(8.78,2.2){\footnotesize $4$}
\usefont{T1}{ptm}{m}{n}
\rput(9.64,1.32){\footnotesize $3$}
\usefont{T1}{ptm}{m}{n}
\rput(9.26,1.78){\footnotesize $3$}
\usefont{T1}{ptm}{m}{n}
\rput(5.26,0.1){\Large $P$}
\psdots[dotsize=0.12](0.06,-3.0)
\psdots[dotsize=0.12](1.26,-1.8)
\psdots[dotsize=0.12](1.66,-2.2)
\psdots[dotsize=0.12](2.86,-1.0)
\psdots[dotsize=0.12](3.26,-1.4)
\psdots[dotsize=0.12](3.66,-1.8)
\psdots[dotsize=0.12](4.06,-2.2)
\psdots[dotsize=0.12](4.46,-2.6)
\psdots[dotsize=0.12](5.66,-1.4)
\psdots[dotsize=0.12](6.06,-1.8)
\psdots[dotsize=0.12](6.46,-2.2)
\psdots[dotsize=0.12](6.86,-2.6)
\psdots[dotsize=0.12](7.26,-3.0)
\psdots[dotsize=0.12](8.46,-1.8)
\psdots[dotsize=0.12](8.86,-2.2)
\psdots[dotsize=0.12](9.26,-2.6)
\psdots[dotsize=0.12](9.66,-3.0)
\psdots[dotsize=0.12](10.86,-1.8)
\psdots[dotsize=0.12](11.26,-2.2)
\psdots[dotsize=0.12](11.66,-2.6)
\psdots[dotsize=0.12](12.86,-1.4)
\psdots[dotsize=0.12](13.26,-1.8)
\psdots[dotsize=0.12](13.66,-2.2)
\psdots[dotsize=0.12](14.06,-2.6)
\psdots[dotsize=0.12](15.26,-1.4)
\psdots[dotsize=0.12](15.66,-1.8)
\psdots[dotsize=0.12](16.06,-2.2)
\psdots[dotsize=0.12](16.46,-2.6)
\psdots[dotsize=0.12](16.86,-3.0)
\psline[linewidth=0.04cm](0.06,-3.0)(1.26,-1.8)
\psline[linewidth=0.04cm](1.26,-1.8)(1.66,-2.2)
\psline[linewidth=0.04cm](1.66,-2.2)(2.86,-1.0)
\psline[linewidth=0.04cm](2.86,-1.0)(4.46,-2.6)
\psline[linewidth=0.04cm](4.46,-2.6)(5.66,-1.4)
\psline[linewidth=0.04cm](5.66,-1.4)(7.26,-3.0)
\psline[linewidth=0.04cm](7.26,-3.0)(8.46,-1.8)
\psline[linewidth=0.04cm](8.46,-1.8)(9.66,-3.0)
\psline[linewidth=0.04cm](9.66,-3.0)(10.86,-1.8)
\psline[linewidth=0.04cm](10.86,-1.8)(11.66,-2.6)
\psline[linewidth=0.04cm](11.66,-2.6)(12.86,-1.4)
\psline[linewidth=0.04cm](12.86,-1.4)(14.06,-2.6)
\psline[linewidth=0.04cm](14.06,-2.6)(15.26,-1.4)
\psline[linewidth=0.04cm](15.26,-1.4)(16.86,-3.0)
\psline[linewidth=0.04cm,linestyle=dotted,dotsep=0.16cm](0.06,-3.0)(16.86,-3.0)
\usefont{T1}{ptm}{m}{n}
\rput(8.46,-3.5){\Large$Q$}
\usefont{T1}{ptm}{m}{n}
\rput(4.0,-2.755){\Large $Q_1$}
\usefont{T1}{ptm}{m}{n}
\rput(8.43,-2.755){\Large $Q_2$}
\usefont{T1}{ptm}{m}{n}
\rput(13.41,-2.755){\Large $Q_3$}
\end{pspicture} 
}

\bigskip

Consider the weak composition $\lam = (1,2,2,1)$ of $L(P)= 6$ of length $4$. 
The word on the top level of $P$ is $\omega _4^P:= 444331$. The path $P\t_{(1,2,2,1)}Q$ is:
\medskip

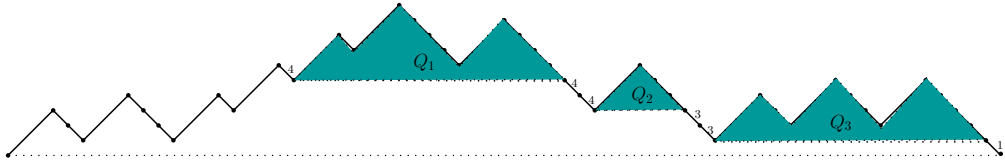
\scalebox{0.5} 
{
\begin{pspicture}(0,-2.08)(26.72,2.08)
\definecolor{color732b}{rgb}{0.0,0.6,0.6}
\psdots[dotsize=0.12](0.06,-2.0)
\psdots[dotsize=0.12](1.26,-0.8)
\psdots[dotsize=0.12](1.66,-1.2)
\psdots[dotsize=0.12](2.06,-1.6)
\psdots[dotsize=0.12](3.26,-0.4)
\psdots[dotsize=0.12](3.66,-0.8)
\psdots[dotsize=0.12](4.06,-1.2)
\psdots[dotsize=0.12](4.46,-1.6)
\psdots[dotsize=0.12](5.66,-0.4)
\psdots[dotsize=0.12](6.06,-0.8)
\psdots[dotsize=0.12](7.26,0.4)
\psdots[dotsize=0.12](7.66,0.0)
\psdots[dotsize=0.12](8.86,1.2)
\psdots[dotsize=0.12](9.26,0.8)
\psdots[dotsize=0.12](10.46,2.0)
\psdots[dotsize=0.12](10.86,1.6)
\psdots[dotsize=0.12](11.26,1.2)
\psdots[dotsize=0.12](11.66,0.8)
\psdots[dotsize=0.12](12.06,0.4)
\psdots[dotsize=0.12](13.26,1.6)
\psdots[dotsize=0.12](13.66,1.2)
\psdots[dotsize=0.12](14.06,0.8)
\psdots[dotsize=0.12](14.46,0.4)
\psdots[dotsize=0.12](14.86,0.0)
\psdots[dotsize=0.12](15.26,-0.4)
\psdots[dotsize=0.12](15.66,-0.8)
\psdots[dotsize=0.12](16.86,0.4)
\psdots[dotsize=0.12](17.26,0.0)
\psdots[dotsize=0.12](17.66,-0.4)
\psdots[dotsize=0.12](18.06,-0.8)
\psdots[dotsize=0.12](18.46,-1.2)
\psdots[dotsize=0.12](18.86,-1.6)
\psdots[dotsize=0.12](20.06,-0.4)
\psdots[dotsize=0.12](20.46,-0.8)
\psdots[dotsize=0.12](20.86,-1.2)
\psdots[dotsize=0.12](22.06,0.0)
\psdots[dotsize=0.12](22.46,-0.4)
\psdots[dotsize=0.12](22.86,-0.8)
\psdots[dotsize=0.12](23.26,-1.2)
\psdots[dotsize=0.12](24.46,0.0)
\psdots[dotsize=0.12](24.86,-0.4)
\psdots[dotsize=0.12](25.26,-0.8)
\psdots[dotsize=0.12](25.66,-1.2)
\psdots[dotsize=0.12](26.06,-1.6)
\psline[linewidth=0.04cm](0.06,-2.0)(1.26,-0.8)
\psline[linewidth=0.04cm](1.26,-0.8)(2.06,-1.6)
\psline[linewidth=0.04cm](2.06,-1.6)(3.26,-0.4)
\psline[linewidth=0.04cm](3.26,-0.4)(4.46,-1.6)
\psline[linewidth=0.04cm](4.46,-1.6)(5.66,-0.4)
\psline[linewidth=0.04cm](5.66,-0.4)(6.06,-0.8)
\psline[linewidth=0.04cm](6.06,-0.8)(7.26,0.4)
\psline[linewidth=0.04cm](7.26,0.4)(7.66,0.0)
\psline[linewidth=0.04cm](7.66,0.0)(8.86,1.2)
\psline[linewidth=0.04cm](8.86,1.2)(9.26,0.8)
\psline[linewidth=0.04cm](9.26,0.8)(10.46,2.0)
\psline[linewidth=0.04cm](10.46,2.0)(12.06,0.4)
\psline[linewidth=0.04cm](12.06,0.4)(13.26,1.6)
\psline[linewidth=0.04cm](13.26,1.6)(15.66,-0.8)
\psline[linewidth=0.04cm](15.66,-0.8)(16.86,0.4)
\psline[linewidth=0.04cm](16.86,0.4)(18.86,-1.6)
\psline[linewidth=0.04cm](18.86,-1.6)(20.06,-0.4)
\psline[linewidth=0.04cm](20.06,-0.4)(20.86,-1.2)
\psline[linewidth=0.04cm](20.86,-1.2)(22.06,0.0)
\psline[linewidth=0.04cm](22.06,0.0)(23.26,-1.2)
\psline[linewidth=0.04cm](23.26,-1.2)(24.46,0.0)
\psline[linewidth=0.04cm](24.46,0.0)(26.46,-2.0)
\psline[linewidth=0.04cm,linestyle=dotted,dotsep=0.16cm](0.06,-2.0)(26.46,-2.0)
\usefont{T1}{ptm}{m}{n}
\rput(7.58,0.3){\footnotesize $4$}
\psline[linewidth=0.048cm,linestyle=dotted,dotsep=0.16cm](7.68,0.0)(14.9,-0.02)
\psline[linewidth=0.048cm,linestyle=dotted,dotsep=0.16cm](15.7,-0.8)(17.98,-0.78)
\psline[linewidth=0.048cm,linestyle=dotted,dotsep=0.16cm](18.86,-1.64)(26.08,-1.64)
\usefont{T1}{ptm}{m}{n}
\rput(15.12,-0.08){\footnotesize $4$}
\usefont{T1}{ptm}{m}{n}
\rput(15.58,-0.52){\footnotesize $4$}
\usefont{T1}{ptm}{m}{n}
\rput(18.74,-1.32){\footnotesize $3$}
\usefont{T1}{ptm}{m}{n}
\rput(18.4,-0.9){\footnotesize $3$}
\psdots[dotsize=0.12](26.46,-1.96)
\usefont{T1}{ptm}{m}{n}
\rput(26.44,-1.72){\footnotesize $1$}
\pspolygon[linewidth=0.04,linestyle=dotted,dotsep=0.16cm,fillstyle=solid,fillcolor=color732b](14.86,0.02)(7.66,0.0)(8.84,1.2)(9.24,0.74)(10.48,2.0)(12.06,0.36)(13.22,1.64)(14.78,0.1)(14.88,0.0)
\pspolygon[linewidth=0.04,linestyle=dotted,dotsep=0.16cm,fillstyle=solid,fillcolor=color732b](18.0,-0.78)(15.68,-0.8)(16.88,0.4)
\pspolygon[linewidth=0.04,linestyle=dotted,dotsep=0.16cm,fillstyle=solid,fillcolor=color732b](26.0,-1.58)(18.84,-1.62)(20.06,-0.38)(20.9,-1.2)(22.08,0.06)(23.28,-1.28)(24.48,0.06)(26.06,-1.62)
\usefont{T1}{ptm}{m}{n}
\rput(11.13,0.485){\Large $Q_1$}
\usefont{T1}{ptm}{m}{n}
\rput(16.93,-0.395){\Large $Q_2$}
\usefont{T1}{ptm}{m}{n}
\rput(22.21,-1.135){\Large $Q_3$}
\usefont{T1}{ptm}{m}{n}

\end{pspicture} 
}

\end{example}
\medskip

The last point of Remark \ref{remcoloring} implies that for any $P\in \Dy_{n_1}^m$, any $Q = Q_1\t_0\dots \t_0 Q_r$ and any $\lam \in \Lambda _r(P)$, the set of down steps of 
$P*_{\lam}Q$ is:
$$\dw (P*_{\lam}Q) = \dw (P)\ \coprod\ \dw (Q),$$
and the standard coloring $\alpha _{P*_{\lam}Q}$ is described by:
$$\alpha _{P*_{\lam}Q}(d) = \begin{cases} \alpha _P(d),& {\rm for}\ d\in \dw(P),\\
\alpha _Q(d) + n_1,& {\rm for}\ d\in \dw(Q).\end{cases} \leqno (3.1)$$
\medskip

\begin{notation} \label{notparti} Let $P$ be a Dyck path with $\dw _{n}(P)=(d_{1}^P,\dots ,d_{L(P)}^P)$ , and let $\lam = (\lambda _0,\dots ,\lambda _r)$ be a weak composition of $L(P)$. For $0\leq i\leq m$, we denote by $\Lambda _r ^i(P)$ the set of all weak compositions $\lam$ of length $r+1$ such that the restriction $\alpha_P (d_{L(P){-}\lambda _r+1}^P),\dots ,\alpha_P(d_{L(P)}^P)$ of the word $\omega _{n}^P$ to its last $\lambda_r$ letters satisfies the following conditions:\begin{enumerate}
\item any digit in the word $\alpha_P (d_{L(P)-\lambda _r+1}^P),\dots ,\alpha_P(d_{L(P)}^P)$ appears at most $i$ times,
\item there exists at least one integer $1\leq i_0\leq n$ such that $i_0$ appears exactly $i$ times in $\alpha_P (d_{L(P){-}\lambda _r+1}^P),\dots ,\alpha_P(d_{L(P)}^P)$.\end{enumerate}\end{notation}
\medskip

For example, for $P = (0,2,1,3,4)\in \Dy _5^2$, 
\medskip

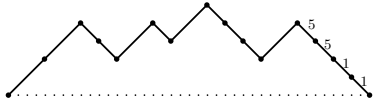
\scalebox{0.6} 
{
\begin{pspicture}(0,-1.08)(8.22,1.08)
\psdots[dotsize=0.12](0.06,-1.0)
\psdots[dotsize=0.12](0.86,-0.2)
\psdots[dotsize=0.12](1.66,0.6)
\psdots[dotsize=0.12](2.06,0.2)
\psdots[dotsize=0.12](2.46,-0.2)
\psdots[dotsize=0.12](3.26,0.6)
\psdots[dotsize=0.12](3.66,0.2)
\psdots[dotsize=0.12](4.46,1.0)
\psdots[dotsize=0.12](4.86,0.6)
\psdots[dotsize=0.12](5.26,0.2)
\psdots[dotsize=0.12](5.66,-0.2)
\psdots[dotsize=0.12](6.46,0.6)
\psdots[dotsize=0.12](6.86,0.2)
\psdots[dotsize=0.12](7.26,-0.2)
\psdots[dotsize=0.12](7.66,-0.6)
\psdots[dotsize=0.12](8.06,-1.0)
\psline[linewidth=0.04cm](0.06,-1.0)(1.66,0.6)
\psline[linewidth=0.04cm](1.66,0.6)(2.46,-0.2)
\psline[linewidth=0.04cm](2.46,-0.2)(3.26,0.6)
\psline[linewidth=0.04cm](3.26,0.6)(3.66,0.2)
\psline[linewidth=0.04cm](3.66,0.2)(4.46,1.0)
\psline[linewidth=0.04cm](4.46,1.0)(5.66,-0.2)
\psline[linewidth=0.04cm](5.66,-0.2)(6.46,0.6)
\psline[linewidth=0.04cm](6.46,0.6)(8.06,-1.0)
\psline[linewidth=0.04cm,linestyle=dotted,dotsep=0.16cm](0.06,-1.0)(8.06,-1.0)
\usefont{T1}{ptm}{m}{n}
\rput(7.54,-0.3){\footnotesize $1$}
\usefont{T1}{ptm}{m}{n}
\rput(7.94,-0.7){\footnotesize $1$}
\usefont{T1}{ptm}{m}{n}
\rput(6.78,0.56){\footnotesize $5$}
\usefont{T1}{ptm}{m}{n}
\rput(7.14,0.12){\footnotesize $5$}
\end{pspicture} 
}
\medskip

 we get that $\lam _1= (1, 1, 2)$ belongs to $\Lambda _2^2(P)$, while $\lam _2 = (0, 3,1)$ belongs to $\Lambda _2^1(P)$. 
 \medskip
 
 Observe that
 $$\Lambda _r^0(P) = \{ (\lambda _0,\dots ,\lambda _{r-1}, 0)\mid \sum _{i=0}^{r-1}\lambda _i = L(P)\ {\rm and}\ r\geq 1\}.$$
 \medskip
 
The set of all weak compositions of $L(P)$ is the disjoint union $\coprod _{r\geq 0}\bigl (\coprod _{i=0}^m\Lambda _r^i (P)\bigr )$, for any $m$-Dyck path $P$ of size $n$. 
 \medskip
 
 The following result is a straightforward consequence of Lemma \ref{lemma1} and the definition of $*_{\lambda}$.

\begin{lemma} \label{lemma2} Let $P = P_1\t _0\dots \t_0 P_s$ in $\Dy _{n_1}^m$ and $Q= Q_1\t_0\dots \t_0Q_r$ in $\Dy _{n_2}^m$ be two Dyck paths, where $P_1,\dots ,P_s,Q_1,\dots ,Q_r$ are prime, and let $\lam \in \Lambda _r^i(P)$ be a weak composition. We have that:\begin{enumerate}
\item if $i > 0$, then
$$P*_{\lam}Q = P_1\t _0\dots \t_0 P_{s-1}\t_0 (P_s*_{\lam } Q),$$
where $P_s*_{\lam } Q$ is prime.
\item if $i=0$, then $\lam = (\lambda _0,\dots ,\lambda _{r-1},0)$ and
$$P*_{\lam}Q = P_1\t _0\dots \t_0 P_{s-1}\t_0 (P_s*_{\lam } (Q_1\t _0 \dots \t_0 Q_{j_0}))\t _0 Q_{j_0+1}\t_0\dots \t_0 Q_r,$$
where $j_0$ is the maximal element of $\{0,\dots , r{-}1\}$ such that $\lambda _{j_0}\neq 0$.
\end{enumerate}
\end{lemma}
\bigskip
 
The product on the graded vector space $\K[\Dy^m]$, spanned by the set of all $m$-Dyck paths, is defined as follows.
 
 \begin{definition} \label{defprodi} Let $P\in \Dy _{n_1}^m$ and $Q\in \Dy _{n_2}^m$ be two Dyck paths, such that $Q = Q_1\t _0\dots \t_0 Q_r$ with $Q_i$ prime, $1\leq i\leq r$. For any integer $0\leq j\leq m$, define
$$P*_jQ= \sum _{\lam\in \Lambda _r^j(P)} P*_{\lam}Q.$$
The product extends in a unique way to a linear map from $\K[\Dy^m]\ot \K[\Dy^m]$ to $\K[\Dy^m]$.\end{definition}
\medskip

\begin{example} Let $P = (1,3)$ be the $2$-Dyck path 
\medskip

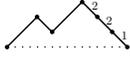
\scalebox{0.5} 
{
\begin{pspicture}(0,-0.68)(3.46,0.68)
\psdots[dotsize=0.12](0.06,-0.6)
\psdots[dotsize=0.12](0.86,0.2)
\psdots[dotsize=0.12](1.26,-0.2)
\psdots[dotsize=0.12](2.06,0.6)
\psdots[dotsize=0.12](2.46,0.2)
\psdots[dotsize=0.12](2.86,-0.2)
\psdots[dotsize=0.12](3.26,-0.6)
\psline[linewidth=0.04cm](0.06,-0.6)(0.86,0.2)
\psline[linewidth=0.04cm](0.86,0.2)(1.26,-0.2)
\psline[linewidth=0.04cm](1.26,-0.2)(2.06,0.6)
\psline[linewidth=0.04cm](2.06,0.6)(3.26,-0.6)
\psline[linewidth=0.04cm,linestyle=dotted,dotsep=0.16cm](0.06,-0.6)(3.26,-0.6)
\usefont{T1}{ptm}{m}{n}
\rput(2.4,0.5){\footnotesize $2$}
\usefont{T1}{ptm}{m}{n}
\rput(2.78,0.1){\footnotesize $2$}
\usefont{T1}{ptm}{m}{n}
\rput(3.18,-0.3){\footnotesize $1$}
\end{pspicture} 
}
\medskip

and let $Q = (0, 2, 4,2) = (0,2,4)\t_0 (2)$ in $\Dy_4^2$, 
\medskip

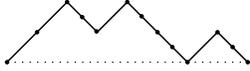
\scalebox{0.5} 
{
\begin{pspicture}(0,-0.88)(6.54,0.88)
\definecolor{color766b}{rgb}{0.0,0.6,0.6}
\psdots[dotsize=0.12](1.66,0.8)
\psdots[dotsize=0.12](2.06,0.4)
\psdots[dotsize=0.12](2.44,0.02)
\psdots[dotsize=0.12](3.64,0.42)
\psline[linewidth=0.04cm,fillcolor=color766b](1.66,0.8)(2.46,0.0)
\psdots[dotsize=0.12](3.26,0.8)
\psdots[dotsize=0.12](4.06,0.0)
\psdots[dotsize=0.12](4.46,-0.4)
\psdots[dotsize=0.12](4.86,-0.8)
\psdots[dotsize=0.12](5.66,0.0)
\psdots[dotsize=0.12](6.06,-0.4)
\psdots[dotsize=0.12](6.46,-0.8)
\psline[linewidth=0.04cm,fillcolor=color766b](2.46,0.0)(3.26,0.8)
\psline[linewidth=0.04cm,fillcolor=color766b](3.26,0.8)(4.86,-0.8)
\psline[linewidth=0.04cm,fillcolor=color766b](4.86,-0.8)(5.66,0.0)
\psline[linewidth=0.04cm,fillcolor=color766b](5.66,0.0)(6.46,-0.8)
\psdots[dotsize=0.12](0.86,0.0)
\psdots[dotsize=0.12](0.06,-0.8)
\psline[linewidth=0.04cm,fillcolor=color766b,linestyle=dotted,dotsep=0.16cm](0.06,-0.8)(6.46,-0.8)
\psline[linewidth=0.04cm,fillcolor=color766b](0.06,-0.8)(1.66,0.8)
\end{pspicture} 
}
\medskip

we get that $P*_0Q = P*_{(3,0,0)}Q + P*_{(2,1,0)}Q + P*_{(1,2,0)}Q + P*_{(0,3,0)}Q =$
$$(1,3,0,2,4,2) + (1,2,0,2,5,2) + (1,1,0,2,6,2) + (1,0,0,2,7,2) =$$

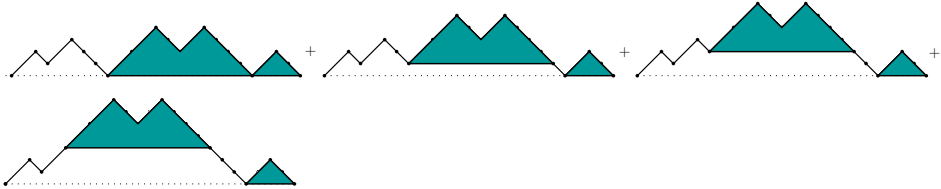
\scalebox{0.4} 
{
\begin{pspicture}(0,-3.08)(31.26,3.08)
\definecolor{color507b}{rgb}{0.0,0.6,0.6}
\psdots[dotsize=0.12](0.26,0.6)
\psdots[dotsize=0.12](1.06,1.4)
\psdots[dotsize=0.12](1.46,1.0)
\psdots[dotsize=0.12](2.26,1.8)
\psdots[dotsize=0.12](2.66,1.4)
\psdots[dotsize=0.12](3.06,1.0)
\psdots[dotsize=0.12](3.46,0.6)
\psdots[dotsize=0.12](4.26,1.4)
\psdots[dotsize=0.12](5.06,2.2)
\psdots[dotsize=0.12](5.46,1.8)
\psdots[dotsize=0.12](5.86,1.4)
\psdots[dotsize=0.12](6.66,2.2)
\psdots[dotsize=0.12](7.06,1.8)
\psdots[dotsize=0.12](7.46,1.4)
\psdots[dotsize=0.12](7.86,1.0)
\psdots[dotsize=0.12](8.26,0.6)
\psdots[dotsize=0.12](9.06,1.4)
\psdots[dotsize=0.12](9.46,1.0)
\psdots[dotsize=0.12](9.86,0.6)
\psline[linewidth=0.04cm,linestyle=dotted,dotsep=0.16cm](0.26,0.6)(0.06,0.6)
\psline[linewidth=0.04cm,linestyle=dotted,dotsep=0.16cm](0.26,0.6)(9.86,0.6)
\psline[linewidth=0.04cm](0.26,0.6)(1.06,1.4)
\psline[linewidth=0.04cm](1.06,1.4)(1.46,1.0)
\psline[linewidth=0.04cm](1.46,1.0)(2.26,1.8)
\psline[linewidth=0.04cm](2.26,1.8)(3.46,0.6)
\psline[linewidth=0.04cm](3.46,0.6)(5.06,2.2)
\psline[linewidth=0.04cm](5.06,2.2)(5.86,1.4)
\psline[linewidth=0.04cm](5.86,1.4)(6.66,2.2)
\psline[linewidth=0.04cm](6.66,2.2)(8.26,0.6)
\psline[linewidth=0.04cm](8.26,0.6)(9.06,1.4)
\psline[linewidth=0.04cm](9.06,1.4)(9.86,0.6)
\usefont{T1}{ptm}{m}{n}
\rput(10.2,1.405){\Large $+$}
\psdots[dotsize=0.12](10.66,0.6)
\psdots[dotsize=0.12](11.46,1.4)
\psdots[dotsize=0.12](11.86,1.0)
\psdots[dotsize=0.12](12.66,1.8)
\psdots[dotsize=0.12](13.06,1.4)
\psdots[dotsize=0.12](13.46,1.0)
\psdots[dotsize=0.12](14.26,1.8)
\psdots[dotsize=0.12](15.06,2.6)
\psdots[dotsize=0.12](15.46,2.2)
\psdots[dotsize=0.12](15.86,1.8)
\psdots[dotsize=0.12](16.66,2.6)
\psdots[dotsize=0.12](17.06,2.2)
\psdots[dotsize=0.12](17.46,1.8)
\psdots[dotsize=0.12](17.86,1.4)
\psdots[dotsize=0.12](18.26,1.0)
\psdots[dotsize=0.12](18.66,0.6)
\psdots[dotsize=0.12](19.46,1.4)
\psdots[dotsize=0.12](19.86,1.0)
\psdots[dotsize=0.12](20.26,0.6)
\psline[linewidth=0.04cm](10.66,0.6)(11.46,1.4)
\psline[linewidth=0.04cm](11.46,1.4)(11.86,1.0)
\psline[linewidth=0.04cm](11.86,1.0)(12.66,1.8)
\psline[linewidth=0.04cm](12.66,1.8)(13.46,1.0)
\psline[linewidth=0.04cm](13.46,1.0)(15.06,2.6)
\psline[linewidth=0.04cm](15.06,2.6)(15.86,1.8)
\psline[linewidth=0.04cm](15.86,1.8)(16.66,2.6)
\psline[linewidth=0.04cm](17.06,2.2)(18.66,0.6)
\psline[linewidth=0.04cm](18.66,0.6)(19.46,1.4)
\psline[linewidth=0.04cm](19.46,1.4)(20.26,0.6)
\usefont{T1}{ptm}{m}{n}
\rput(20.64,1.365){\Large $+$}
\psdots[dotsize=0.12](21.06,0.6)
\psdots[dotsize=0.12](21.86,1.4)
\psdots[dotsize=0.12](22.26,1.0)
\psdots[dotsize=0.12](23.06,1.8)
\psdots[dotsize=0.12](23.46,1.4)
\psdots[dotsize=0.12](24.26,2.2)
\psdots[dotsize=0.12](25.06,3.0)
\psdots[dotsize=0.12](25.46,2.6)
\psdots[dotsize=0.12](25.86,2.2)
\psdots[dotsize=0.12](26.66,3.0)
\psline[linewidth=0.04cm](16.66,2.6)(17.06,2.2)
\psdots[dotsize=0.12](27.06,2.6)
\psdots[dotsize=0.12](27.46,2.2)
\psdots[dotsize=0.12](27.86,1.8)
\psdots[dotsize=0.12](28.26,1.4)
\psdots[dotsize=0.12](28.66,1.0)
\psdots[dotsize=0.12](29.06,0.6)
\psdots[dotsize=0.12](29.86,1.4)
\psdots[dotsize=0.12](30.26,1.0)
\psdots[dotsize=0.12](30.66,0.6)
\psline[linewidth=0.04cm](21.06,0.6)(21.86,1.4)
\psline[linewidth=0.04cm](21.86,1.4)(22.26,1.0)
\psline[linewidth=0.04cm](22.26,1.0)(23.06,1.8)
\psline[linewidth=0.04cm](23.06,1.8)(23.46,1.4)
\psline[linewidth=0.04cm](23.46,1.4)(25.06,3.0)
\psline[linewidth=0.04cm](25.06,3.0)(25.86,2.2)
\psline[linewidth=0.04cm](25.86,2.2)(26.66,3.0)
\psline[linewidth=0.04cm](26.66,3.0)(29.06,0.6)
\psline[linewidth=0.04cm](29.06,0.6)(29.86,1.4)
\psline[linewidth=0.04cm](29.86,1.4)(30.66,0.6)
\psline[linewidth=0.04cm,linestyle=dotted,dotsep=0.16cm](10.66,0.6)(20.26,0.6)
\psline[linewidth=0.04cm,linestyle=dotted,dotsep=0.16cm](21.06,0.6)(30.66,0.6)
\usefont{T1}{ptm}{m}{n}
\rput(30.94,1.345){\Large $+$}
\psdots[dotsize=0.12](0.06,-3.0)
\psdots[dotsize=0.12](0.86,-2.2)
\psdots[dotsize=0.12](1.26,-2.6)
\psdots[dotsize=0.12](2.06,-1.8)
\psdots[dotsize=0.12](2.86,-1.0)
\psdots[dotsize=0.12](3.66,-0.2)
\psdots[dotsize=0.12](4.06,-0.6)
\psdots[dotsize=0.12](4.46,-1.0)
\psdots[dotsize=0.12](5.26,-0.2)
\psdots[dotsize=0.12](5.66,-0.6)
\psdots[dotsize=0.12](6.06,-1.0)
\psdots[dotsize=0.12](6.46,-1.4)
\psdots[dotsize=0.12](6.86,-1.8)
\psdots[dotsize=0.12](7.26,-2.2)
\psdots[dotsize=0.12](7.66,-2.6)
\psdots[dotsize=0.12](0.06,-3.0)
\psline[linewidth=0.04cm](0.06,-3.0)(0.86,-2.2)
\psline[linewidth=0.04cm](0.86,-2.2)(1.26,-2.6)
\psline[linewidth=0.04cm](1.26,-2.6)(3.66,-0.2)
\psline[linewidth=0.04cm](3.66,-0.2)(4.46,-1.0)
\psline[linewidth=0.04cm,linestyle=dotted,dotsep=0.16cm](0.06,-3.0)(9.86,-3.0)
\psline[linewidth=0.04cm](4.46,-1.0)(5.26,-0.2)
\psdots[dotsize=0.12](8.06,-3.0)
\psdots[dotsize=0.12](8.86,-2.2)
\psdots[dotsize=0.12](9.26,-2.6)
\psdots[dotsize=0.12](9.66,-3.0)
\psline[linewidth=0.04cm](5.26,-0.2)(8.06,-3.0)
\psline[linewidth=0.04cm](8.06,-3.0)(8.86,-2.2)
\psline[linewidth=0.04cm](8.86,-2.2)(9.66,-3.0)
\psdots[dotsize=0.12,fillstyle=solid,dotstyle=o](4.86,-0.6)
\psdots[dotsize=0.12,fillstyle=solid,dotstyle=o](8.46,-2.6)
\pspolygon[linewidth=0.048,fillstyle=solid,fillcolor=color507b](18.26,1.0)(18.26,1.0)(13.46,1.0)(15.06,2.6)(15.86,1.8)(16.66,2.6)(18.26,1.0)
\pspolygon[linewidth=0.048,fillstyle=solid,fillcolor=color507b](18.66,0.6)(19.46,1.4)(20.26,0.6)(20.26,0.6)
\pspolygon[linewidth=0.048,fillstyle=solid,fillcolor=color507b](23.46,1.4)(28.26,1.4)(26.66,3.0)(25.86,2.2)(25.06,3.0)
\pspolygon[linewidth=0.048,fillstyle=solid,fillcolor=color507b](29.06,0.6)(30.66,0.6)(29.86,1.4)
\pspolygon[linewidth=0.048,fillstyle=solid,fillcolor=color507b](2.06,-1.8)(6.86,-1.8)(5.26,-0.2)(4.46,-1.0)(3.66,-0.2)(3.66,-0.2)
\pspolygon[linewidth=0.048,fillstyle=solid,fillcolor=color507b](9.66,-3.0)(8.06,-3.0)(8.86,-2.2)
\pspolygon[linewidth=0.048,fillstyle=solid,fillcolor=color507b](4.66,0.6)(3.46,0.6)(5.06,2.2)(5.86,1.4)(6.66,2.2)(8.26,0.6)
\pspolygon[linewidth=0.048,fillstyle=solid,fillcolor=color507b](9.06,0.6)(8.26,0.6)(9.06,1.4)(9.86,0.6)(9.86,0.6)
\end{pspicture} 
}
\medskip

and 

$$\displaylines {
P*_1Q = P*_{(2,0,1)}Q + P*_{(1,1,1)} Q +P*_{(1,0,2)} Q + P*_{(0,2,1)} Q + P*_{(0,1,2)} Q =\hfill \cr
(1,2,0,2,4,3) + (1,1, 0,2,5,3) + (1,1,0,2,4,4) + (1,0,0,2,6,3) + (1,0,0,2,5,4)= \cr}$$

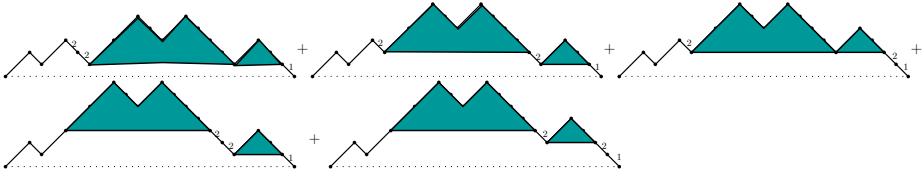
\scalebox{0.4} 
{
\begin{pspicture}(0,-2.78)(30.66,2.78)
\definecolor{color1277b}{rgb}{0.0,0.6,0.6}
\psdots[dotsize=0.12](0.06,0.3)
\psdots[dotsize=0.12](0.86,1.1)
\psdots[dotsize=0.12](1.26,0.7)
\psdots[dotsize=0.12](2.06,1.5)
\psdots[dotsize=0.12](2.46,1.1)
\psdots[dotsize=0.12](2.86,0.7)
\psdots[dotsize=0.12](3.66,1.5)
\psdots[dotsize=0.12](4.46,2.3)
\psdots[dotsize=0.12](4.86,1.9)
\psdots[dotsize=0.12](5.26,1.5)
\psdots[dotsize=0.12](6.06,2.3)
\psdots[dotsize=0.12](6.46,1.9)
\psdots[dotsize=0.12](6.86,1.5)
\psdots[dotsize=0.12](7.26,1.1)
\psdots[dotsize=0.12](7.66,0.7)
\psdots[dotsize=0.12](8.46,1.5)
\psdots[dotsize=0.12](8.86,1.1)
\psdots[dotsize=0.12](9.26,0.7)
\psdots[dotsize=0.12](9.66,0.3)
\psline[linewidth=0.04cm,linestyle=dotted,dotsep=0.16cm](0.06,0.3)(9.66,0.3)
\psline[linewidth=0.04cm](0.06,0.3)(0.86,1.1)
\psline[linewidth=0.04cm](0.86,1.1)(1.26,0.7)
\psline[linewidth=0.04cm](1.26,0.7)(2.06,1.5)
\psline[linewidth=0.04cm](2.06,1.5)(2.86,0.7)
\psline[linewidth=0.04cm](2.86,0.7)(4.46,2.3)
\psline[linewidth=0.04cm](4.46,2.3)(5.26,1.5)
\psline[linewidth=0.04cm](5.26,1.5)(6.06,2.3)
\psline[linewidth=0.04cm](6.06,2.3)(7.66,0.7)
\psline[linewidth=0.04cm](7.66,0.7)(8.46,1.5)
\psline[linewidth=0.04cm](8.46,1.5)(9.66,0.3)
\usefont{T1}{ptm}{m}{n}
\rput(9.94,1.205){\Large $+$}
\psdots[dotsize=0.12](10.26,0.3)
\psdots[dotsize=0.12](11.06,1.1)
\psdots[dotsize=0.12](11.46,0.7)
\psdots[dotsize=0.12](12.26,1.5)
\psdots[dotsize=0.12](12.66,1.1)
\psdots[dotsize=0.12](13.46,1.9)
\psdots[dotsize=0.12](14.26,2.7)
\psdots[dotsize=0.12](14.66,2.3)
\psdots[dotsize=0.12](15.06,1.9)
\psdots[dotsize=0.12](15.86,2.7)
\psdots[dotsize=0.12](16.26,2.3)
\psdots[dotsize=0.12](16.66,1.9)
\psdots[dotsize=0.12](17.06,1.5)
\psdots[dotsize=0.12](17.46,1.1)
\psdots[dotsize=0.12](17.86,0.7)
\psdots[dotsize=0.12](18.66,1.5)
\psdots[dotsize=0.12](19.06,1.1)
\psdots[dotsize=0.12](19.46,0.7)
\psdots[dotsize=0.12](19.86,0.3)
\psline[linewidth=0.04cm](10.26,0.3)(11.06,1.1)
\psline[linewidth=0.04cm](11.06,1.1)(11.46,0.7)
\psline[linewidth=0.04cm](11.46,0.7)(12.26,1.5)
\psline[linewidth=0.04cm](12.26,1.5)(12.66,1.1)
\psline[linewidth=0.04cm](12.66,1.1)(14.26,2.7)
\psline[linewidth=0.04cm](14.26,2.7)(15.06,1.9)
\psline[linewidth=0.04cm](15.06,1.9)(15.86,2.7)
\psline[linewidth=0.04cm](15.86,2.7)(17.86,0.7)
\psline[linewidth=0.04cm](17.86,0.7)(18.66,1.5)
\psline[linewidth=0.04cm](18.66,1.5)(19.86,0.3)
\psline[linewidth=0.04cm,linestyle=dotted,dotsep=0.16cm](10.26,0.3)(19.86,0.3)
\usefont{T1}{ptm}{m}{n}
\rput(20.14,1.205){\Large $+$}
\psdots[dotsize=0.12](20.46,0.3)
\psdots[dotsize=0.12](21.26,1.1)
\psdots[dotsize=0.12](21.66,0.7)
\psdots[dotsize=0.12](22.46,1.5)
\psdots[dotsize=0.12](22.86,1.1)
\psdots[dotsize=0.12](23.66,1.9)
\psdots[dotsize=0.12](24.46,2.7)
\psdots[dotsize=0.12](24.86,2.3)
\psdots[dotsize=0.12](25.26,1.9)
\psdots[dotsize=0.12](26.06,2.7)
\psdots[dotsize=0.12](26.46,2.3)
\psdots[dotsize=0.12](26.86,1.9)
\psdots[dotsize=0.12](27.26,1.5)
\psdots[dotsize=0.12](27.66,1.1)
\psdots[dotsize=0.12](28.46,1.9)
\psdots[dotsize=0.12](28.86,1.5)
\psdots[dotsize=0.12](29.26,1.1)
\psdots[dotsize=0.12](29.66,0.7)
\psdots[dotsize=0.12](30.06,0.3)
\psline[linewidth=0.04cm,linestyle=dotted,dotsep=0.16cm](20.46,0.3)(30.06,0.3)
\psline[linewidth=0.04cm](20.46,0.3)(21.26,1.1)
\psline[linewidth=0.04cm](21.26,1.1)(21.66,0.7)
\psline[linewidth=0.04cm](21.66,0.7)(22.46,1.5)
\psline[linewidth=0.04cm](22.46,1.5)(22.86,1.1)
\psline[linewidth=0.04cm](22.86,1.1)(24.46,2.7)
\psline[linewidth=0.04cm](24.46,2.7)(25.26,1.9)
\psline[linewidth=0.04cm](25.26,1.9)(26.06,2.7)
\psline[linewidth=0.04cm](26.06,2.7)(27.66,1.1)
\psline[linewidth=0.04cm](27.66,1.1)(28.46,1.9)
\psline[linewidth=0.04cm](28.46,1.9)(30.06,0.3)
\usefont{T1}{ptm}{m}{n}
\rput(30.34,1.205){\Large $+$}
\psdots[dotsize=0.12](0.06,-2.7)
\psdots[dotsize=0.12](0.86,-1.9)
\psdots[dotsize=0.12](4.06,-0.3)
\psdots[dotsize=0.12](2.86,-0.7)
\psdots[dotsize=0.12](3.66,0.1)
\psdots[dotsize=0.12](6.86,-1.5)
\psdots[dotsize=0.12](4.46,-0.7)
\psdots[dotsize=0.12](5.26,0.1)
\psdots[dotsize=0.12](5.66,-0.3)
\psdots[dotsize=0.12](6.06,-0.7)
\psdots[dotsize=0.12](6.46,-1.1)
\psdots[dotsize=0.12](1.26,-2.3)
\psdots[dotsize=0.12](2.06,-1.5)
\psdots[dotsize=0.12](7.26,-1.9)
\psdots[dotsize=0.12](7.66,-2.3)
\psdots[dotsize=0.12](8.46,-1.5)
\psdots[dotsize=0.12](8.86,-1.9)
\psdots[dotsize=0.12](9.26,-2.3)
\psdots[dotsize=0.12](9.66,-2.7)
\psline[linewidth=0.04cm](0.06,-2.7)(0.86,-1.9)
\psline[linewidth=0.04cm](0.86,-1.9)(1.26,-2.3)
\psline[linewidth=0.04cm](1.26,-2.3)(3.66,0.1)
\psline[linewidth=0.04cm](3.66,0.1)(4.46,-0.7)
\psline[linewidth=0.04cm](4.46,-0.7)(5.26,0.1)
\psline[linewidth=0.04cm](5.26,0.1)(7.66,-2.3)
\psline[linewidth=0.04cm](7.66,-2.3)(8.46,-1.5)
\psline[linewidth=0.04cm](8.46,-1.5)(9.66,-2.7)
\psline[linewidth=0.04cm,linestyle=dotted,dotsep=0.16cm](0.06,-2.7)(9.66,-2.7)
\usefont{T1}{ptm}{m}{n}
\rput(10.34,-1.795){\Large $+$}
\psdots[dotsize=0.12](10.86,-2.7)
\psdots[dotsize=0.12](11.66,-1.9)
\psdots[dotsize=0.12](12.06,-2.3)
\psdots[dotsize=0.12](12.86,-1.5)
\psdots[dotsize=0.12](13.66,-0.7)
\psdots[dotsize=0.12](14.46,0.1)
\psdots[dotsize=0.12](14.86,-0.3)
\psdots[dotsize=0.12](15.26,-0.7)
\psdots[dotsize=0.12](16.06,0.1)
\psdots[dotsize=0.12](16.46,-0.3)
\psdots[dotsize=0.12](16.86,-0.7)
\psdots[dotsize=0.12](17.26,-1.1)
\psdots[dotsize=0.12](17.66,-1.5)
\psdots[dotsize=0.12](18.06,-1.9)
\psdots[dotsize=0.12](18.86,-1.1)
\psdots[dotsize=0.12](19.26,-1.5)
\psdots[dotsize=0.12](19.66,-1.9)
\psdots[dotsize=0.12](20.06,-2.3)
\psdots[dotsize=0.12](20.46,-2.7)
\psline[linewidth=0.04cm,linestyle=dotted,dotsep=0.16cm](10.86,-2.7)(20.46,-2.7)
\psline[linewidth=0.04cm](10.86,-2.7)(11.66,-1.9)
\psline[linewidth=0.04cm](11.66,-1.9)(12.06,-2.3)
\psline[linewidth=0.04cm](12.06,-2.3)(14.46,0.1)
\psline[linewidth=0.04cm](14.46,0.1)(15.26,-0.7)
\psline[linewidth=0.04cm](15.26,-0.7)(16.06,0.1)
\psline[linewidth=0.04cm](16.06,0.1)(18.06,-1.9)
\psline[linewidth=0.04cm](18.06,-1.9)(18.86,-1.1)
\psline[linewidth=0.04cm](18.86,-1.1)(20.46,-2.7)
\usefont{T1}{ptm}{m}{n}
\rput(2.76,1.0){\footnotesize $2$}
\usefont{T1}{ptm}{m}{n}
\rput(2.34,1.38){\footnotesize $2$}
\usefont{T1}{ptm}{m}{n}
\rput(9.56,0.6){\footnotesize $1$}
\usefont{T1}{ptm}{m}{n}
\rput(12.52,1.42){\footnotesize $2$}
\usefont{T1}{ptm}{m}{n}
\rput(17.76,0.96){\footnotesize $2$}
\usefont{T1}{ptm}{m}{n}
\rput(19.74,0.62){\footnotesize $1$}
\usefont{T1}{ptm}{m}{n}
\rput(29.62,0.98){\footnotesize $2$}
\usefont{T1}{ptm}{m}{n}
\rput(22.78,1.4){\footnotesize $2$}
\usefont{T1}{ptm}{m}{n}
\rput(29.98,0.6){\footnotesize $1$}
\usefont{T1}{ptm}{m}{n}
\rput(7.1,-1.6){\footnotesize $2$}
\usefont{T1}{ptm}{m}{n}
\rput(7.56,-2.0){\footnotesize $2$}
\usefont{T1}{ptm}{m}{n}
\rput(9.54,-2.4){\footnotesize $1$}
\usefont{T1}{ptm}{m}{n}
\rput(18.0,-1.62){\footnotesize $2$}
\usefont{T1}{ptm}{m}{n}
\rput(19.96,-2.0){\footnotesize $2$}
\usefont{T1}{ptm}{m}{n}
\rput(20.46,-2.38){\footnotesize $1$}
\pspolygon[linewidth=0.048,fillstyle=solid,fillcolor=color1277b](5.28,0.76)(7.68,0.7)(6.04,2.3)(5.28,1.46)(4.46,2.24)(2.84,0.68)(2.84,0.68)
\pspolygon[linewidth=0.048,fillstyle=solid,fillcolor=color1277b](9.22,0.7)(7.7,0.66)(8.46,1.48)(8.46,1.48)(8.5,1.5)
\pspolygon[linewidth=0.048,fillstyle=solid,fillcolor=color1277b](12.66,1.12)(17.48,1.1)(15.9,2.66)(15.1,1.88)(14.28,2.7)(14.28,2.68)(14.28,2.7)
\pspolygon[linewidth=0.048,fillstyle=solid,fillcolor=color1277b](18.66,1.5)(17.86,0.7)(19.46,0.7)
\pspolygon[linewidth=0.048,fillstyle=solid,fillcolor=color1277b](25.86,1.1)(22.86,1.1)(24.46,2.7)(25.26,1.9)(26.06,2.7)(27.66,1.1)
\pspolygon[linewidth=0.048,fillstyle=solid,fillcolor=color1277b](28.26,1.1)(27.66,1.1)(28.46,1.9)(29.26,1.1)
\pspolygon[linewidth=0.048,fillstyle=solid,fillcolor=color1277b](2.06,-1.5)(6.86,-1.5)(5.26,0.1)(4.46,-0.7)(3.66,0.1)(3.66,0.1)(3.66,0.1)
\pspolygon[linewidth=0.048,fillstyle=solid,fillcolor=color1277b](8.46,-1.5)(9.26,-2.3)(7.66,-2.3)(7.66,-2.3)(7.66,-2.3)
\pspolygon[linewidth=0.048,fillstyle=solid,fillcolor=color1277b](17.66,-1.5)(12.86,-1.5)(14.46,0.1)(15.26,-0.7)(16.06,0.1)
\pspolygon[linewidth=0.048,fillstyle=solid,fillcolor=color1277b](18.06,-1.9)(18.86,-1.1)(19.66,-1.9)(19.66,-1.9)
\end{pspicture} 
}

\end{example}
\medskip

\begin{proposition} \label{proprelac} Let $P\in \Dy _{n_1}^m$ and $Q = Q_1\t _0\dots \t _0 Q_r\in \Dy _{n_2}^m$ be two Dyck paths, with $Q_j\in \Dy _{n_{2j}}^m$ prime for $1\leq j\leq r$. \begin{enumerate}
\item For nonnegative integers $s\geq 1$ and $0\leq i < j\leq m$, the map 
$$\psi_{ij}(P,Q):\{ \Lambda _r^i(P)\times \Lambda _s^j(Q) \longrightarrow \{ (\lam ,{\underline {\delta}})\mid \lam \in \Lambda_r^i(P)\ {\rm and}\ {\underline {\delta}}\in \Lambda _s^j(P*_{\lam}Q)\},$$ 
which sends $(\lam , {\underline {\tau}})\mapsto (\lam , {\underline{\delta}} :=(\tau_0,\dots ,\tau _{s-1},\tau_s +\lambda_r))$ is bijective.
\item For any integer $0\leq i\leq m$, the map $\psi _i^1(P,Q)(\lam ,{\underline {\tau}}): =$
$$ ((\lambda _0,\dots ,\lambda _{r-1},\lambda _r + \dots + \lambda _{r+s{-}j_{\underline {\tau}}}) , (\tau_0,\dots ,\tau_{j_{\underline {\tau}}}+\lambda _r,\lambda _{r+1},\dots ,\lambda _{r+s{-}j_{\underline {\tau}}}))$$ defines a bijection from the set $\{(\lam ,{\underline {\tau}})\mid {\underline {\tau}}\in \Lambda _s^0(Q)\ {\rm and}\ \lam \in \Lambda _{r+s{-}j_{\underline {\tau}}}^i (P)\}$ to the set 
$$\{ ({\underline {\gamma }}, {\underline {\delta }})\mid {\underline {\gamma }} \in {\displaystyle \coprod_{j=i}^m\Lambda _r^j(P)}\ {\rm and}\  {\underline {\delta }}\in \Lambda _s^i(P*_{\underline {\gamma }}Q)\ {\rm such\ that}\ \delta _s \leq \gamma _r\},$$  
where $j_{\underline {\tau}}$ is the maximal integer $0\leq j\leq s-1$ such that $\tau _j > 0$, and $\coprod $ denotes the disjoint union. 
\item For any integer $0\leq i\leq m$, the map
$$\psi _i^2(P,Q) (\lam , {\underline {\tau}}) := (\lam , {\underline{\delta}}:= (\tau_0,\dots ,\tau _{s-1}, \tau _s+\lambda _r),$$ from 
$\Lambda _r^i(P)\t ({\displaystyle \coprod _{j=1}^i \Lambda _s^j(Q)})$ to the set $$\{ ({\underline{\gamma}} , {\underline {\delta}})\mid {\underline{\gamma}} \in \Lambda _r^i(P)\ {\rm and}\ {\underline {\delta}}\in \Lambda _s^i(P*_{\underline {\gamma}}Q)\ {\rm such\ that}\ \gamma _r < \delta _s\},$$
is bijective.
\end{enumerate}
\end{proposition}
\medskip

\begin{proo} $(1)$ For the first point, let $\lam \in \Lambda _r^i(P)$ and ${\underline {\tau}}\in  \Lambda _s^j(Q)$ be two weak compositions. 

If $\dw_{n_1}(P) = (d_{1}^P,\dots ,d_{L(P)}^P)$ and $ \dw_{n_2}(Q)= (d_{1}^Q,\dots ,d_{L(Q)}^Q)$, then:
$$\dw _{n_1+n_2}(P*_{\lam}Q) = (d_{1}^Q,\dots ,d_{L(Q)}^Q, d_{L(P){-}\lambda_r+1}^P,\dots ,d_{L(P)}^P).$$

The map $\psi _{ij}$ is defined by the formula:
$$\psi _{ij}(\lam, {\underline {\tau}}) := (\lam , {\underline{\delta}} :=(\tau_0,\dots ,\tau _{s{-}1},\tau_s +\lambda_r)).$$

Clearly, $\lam $ belongs to $ \Lambda_r^i(P)$. On the other hand, 
$$\dw _{n_1+n_2}(P*_{\lam}Q) = (d_1^Q,\dots ,d_{L(Q)}^Q,d_{L(P){-}\lambda _r+1}^P,\dots ,d_{L(P)}^P),$$
 which implies that the subset of the last 
$\tau_s +\lambda_r$ down steps of $P*_{\lam}Q$ is
$(d_{L(Q){-}\tau _s+1}^Q, \dots ,d_{L(Q)}^Q, d_{L(P){-}\lambda _r+1}^P,\dots ,d_{L(P)}^P)$.
\medskip

Note that: \begin{enumerate}
\item $\alpha _{P*_{\lam}Q}(d_{L(P){-}\lambda _r+1}^P)\dots \alpha_{P*_{\lam}Q}(d_{L(P)}^P)$ is a sequence of elements in the set $\{1,\dots n_1\}$ such that any digit appears at most $i$ times.
\item $\alpha _{P*_{\lam}Q}(d_{L(Q){-}\tau _s+1}^Q)\dots \alpha _{P*_{\lam}Q}(d_{L(Q)}^Q)$ is a sequence of elements in the set $\{n_1+1,\dots ,n_1+n_2\}$ where there exists at least one digit that appears $j$ times, and no digit appears more than $j$ times.\end{enumerate}
So, ${\underline{\delta}}$ belongs to $\Lambda _s^j(P*_{\lam}Q)$.
\medskip

For any pair of weak compositions $\lam \in \Lambda _r^i(P)$ and ${\underline {\delta}}\in \Lambda _s^j(P*_{\lam}Q)$, we get:
\begin{multline}\nonumber \omega_{n_1+n_2}^{P*_{\lam}Q}=\\
\qquad\alpha _Q(d_{L(Q){-}\delta_s+1}^Q)+n_1,\dots ,\alpha _Q(d_{L(Q)}^Q)+n_1,\alpha _P(d_{L(P){-}\lambda _r+1}^P),\dots ,\alpha _P(d_{L(P)}^P).\end{multline}

As the expression $\alpha_P(d_{L(P){-}\lambda_r+1})\dots \alpha _P(d_{L(P)})$ is a word in the alphabet $\{ 1,\dots ,n_1\}$ such that no digit appears more than $i$ times, and $i < j$, then ${\underline {\tau}} := (\delta _0,\dots ,\delta _s{-}\lambda _r)$ must belong to $\Lambda _s^j(Q)$. 
 \medskip
 
It is immediate to prove that the map $(\lam ,{\underline {\delta}})\mapsto (\lam ,{\underline {\tau}})$ is inverse to $\psi_{ij}(P,Q)$, which ends the proof of $(1)$.
\bigskip

$(2)$ If $\lam\in \Lambda _{r+s{-}j_{\underline {\tau}}}^i (P)$ and ${\underline {\tau}}\in \Lambda _s^0(Q)$, then it is immediate to verify that \begin{enumerate}[(i)]
\item ${\underline {\gamma }} = (\lambda _0,\dots ,\lambda _{r-1},\lambda _r + \dots + \lambda _{r+s{-}j_{\underline {\tau}}})$ belongs to $\Lambda _r^j(P)$, for $i\leq j\leq m$,
\item $ {\underline {\delta }} = (\tau_0,\dots ,\tau_{j_{\underline {\tau}}}+\lambda _r,\lambda _{r+1},\dots ,\lambda _{r+s{-}j_{\underline {\tau}}})$ belongs to $\Lambda _s^i(P*_{\underline {\gamma }}Q)$,
\item $\delta _s = \lambda _{r+s{-}j_{\underline {\tau}}}\leq \gamma _r = \lambda _r + \dots + \lambda _{r+s{-}j_{\underline {\tau}}}$.\end{enumerate}

Assume that we have two weak compositions ${\underline {\gamma }} = (\gamma _0,\dots ,\gamma _r)\in {\displaystyle \coprod _{j=i}^m\Lambda _r^j(P)}$ and $ {\underline {\delta }} = (\delta _0, \dots ,\delta _s) \in \Lambda _s^i(P*_{\underline {\gamma }}Q)$ such that $\delta _s \leq \gamma _r$.
\medskip

Let $j_0$ be the maximal integer $0\leq j_0\leq s{-}1$, such that $\delta _{j_0}+\dots +\delta _s > \gamma _r$. 
\medskip

Define
\begin{enumerate}[(a)]
\item $\lam := (\gamma _0,\dots , \gamma_{r-1}, \gamma_r{-}\delta_{j_0+1}{-}\dots {-}\delta_s , \delta_{j_0+1},\dots, \delta_s)$,
\item ${\underline {\tau}}:= (\delta _0,\dots ,\delta_{j_0{-}1}, \delta _{j_0}+\dots +\delta _s -\gamma _r, 0,\dots ,0)$.\end{enumerate}
It is clear that $\lam \in \Lambda _{r+s-j_0}^i(P)$, 
${\underline {\tau}}\in \Lambda _s^0(Q)$ and $\psi _i^1(P,Q)(\lam ,{\underline {\tau}}) = ({\underline {\gamma }},  {\underline {\delta }})$, which shows that $\psi _i^1$ is bijective, ending the proof of $(2)$. 
\bigskip

$(3)$ For $\lam \in \Lambda _r^i(P)$ and ${\underline {\tau}}\in \Lambda _s^j(Q)$, for $1\leq j\leq i$, we have that the weak composition $\psi _i^2(P,Q) (\lam , {\underline {\tau}}) = ({\underline{\gamma}}, {\underline{\delta}})$ satisfies the following conditions:\begin{enumerate}[(i)]
\item ${\underline{\gamma}} = \lam$ belongs to $\Lambda _r^i(P)$,
\item the weak composition ${\underline{\tau}}$ belongs to $\Lambda _s^j(Q)$ for some $1\leq j\leq i$. So, the sequence $ \alpha _Q(d_{L(Q){-}\tau _s+1}^Q)+n_1\dots \alpha _Q(d_{L(Q)}^Q)+n_1$ is a word in the digits of $\{n_1+1,\dots ,n_1+n_2\}$ such that each sequence appears at most $j$ times.

On the other hand, the sequence $\alpha_P(d_{L(P)-\gamma _r+1}^P)\dots \alpha _P(d_{L(P)}^P)$ 
is a word in $\{1,\dots ,n_1\}$ such that some digit appears exactly $i$ times in it and no digit appears more than $i$ times. 

The sequence of level $n_1+n_2$ of $P*_{\underline {\gamma}}Q$ is 
\begin{multline}\nonumber\omega _{n_1+n_2}^{P*_{\underline {\gamma}}Q} = \\
 \alpha _Q(d_{L(Q){-}\tau _s+1}^Q)+n_1\dots \alpha _Q(d_{L(Q)}^Q)+n_1\alpha_P(d_{L(P)-\gamma _r+1}^P)
\dots \alpha _P(d_{L(P)}^P),\end{multline}
which shows that ${\underline {\delta}} = (\tau _0,\dots ,\tau_{s-1},\tau _s+\lambda _r)$ belongs to $\Lambda _s^i(P*_{\gamma}Q)$.

\item As ${\underline {\gamma }} = \lam$ and ${\underline {\delta}} = (\tau_0,\dots ,\tau_{s-1}, \tau_s + \lambda _r)$, with $\tau _s > 0$, we get that $\gamma _r < \delta _s$.\end{enumerate}
\medskip

The map $({\underline {\gamma}}, {\underline {\delta}})\mapsto ({\underline {\gamma}}, (\delta_0,\dots ,\delta _{s-1}, \delta _s- \gamma _r))$ is the inverse map of $\psi _i^2(P,Q)$.
\end{proo}
\medskip

\begin{theorem} \label{theorem1} The binary operations $\{ *_j\}_{0\leq j\leq m}$ defined on $\K[\Dy^m]$ satisfy the following relations:\begin{enumerate}
\item $x*_i (y*_j z) = (x*_iy)*_jz$, for $0\leq i < j\leq m$,
\item $x*_i (y*_0z + \dots + y*_iz) = (x*_iy + \dots +x*_my)*_iz$, for $0\leq i\leq m$,\end{enumerate}
where $x,y,z$ are arbitrary elements of $\K[\Dy^m]$.\end{theorem}
\medskip

\begin{proo} Clearly, it suffices to prove the relations for any Dyck paths $P$, $Q$ and $Z$. Suppose that $P\in \Dy _{n_1}^m$, $Q = Q_1\t_0\dots Q_r\in \Dy _{n_2}^m$ and $Z = Z_1\t_0\dots \t_0Z_s\in \Dy _{n_3}^m$, where $Q_1,\dots ,Q_r, Z_1,\dots , Z_s$ are prime Dyck paths.

$(1)$ For $0\leq i < j\leq m$, applying a recursive argument on $s$ and Lemma \ref{lemma2} it is easy to see that, for any pair $(\lam ,{\underline {\tau}})\in \Lambda _r^i(P)\t \Lambda _s^j(Q)$, we get:
\begin{multline}\nonumber
P*_{\lam} (Q*_{\underline {\tau}}Z) = ((P\t _{\lambda _1+\dots +\lambda _r}Q_1)\t _{\lambda _2+\dots +\lambda _r} \dots )\t _{\lambda _r} (Q_r*_{\underline {\tau}} Z) =\\
(((P\t _{\lambda _1+\dots +\lambda _r}Q_1)\t _{\lambda _2+\dots +\lambda _r} \dots )\t _{\lambda _r} Q_r)*_{\underline {\delta }} Z,\end{multline}
where ${\underline {\delta }}= (\tau _0,\dots ,\tau _{s-1}, \tau_s + \lambda _r)$.

Applying the same notation than in Proposition \ref{proprelac}, we get that 

\noi $P*_{\lam} (Q*_{\underline {\tau}} Z) = (P*_{\lam} Q)*_{\underline {\delta }} Z$ if, and only if, $\psi _{ij}(P,Q) (\lam ,{\underline {\tau}}) = (\lam , {\underline {\delta}})$. 
The result follows applying point $(1)$ of Proposition \ref{proprelac}.
\medskip

$(2)$ We write ${\displaystyle \sum _{j = 0}^i P*_i (Q*_jZ) = P*_i (Q*_0Z) + \sum _{j=1}^i P*_i (Q*_j Z)}$ and 
we work the terms on the right hand side separately.
\medskip

$a)$ Suppose that ${\underline {\tau}}\in \Lambda _s^0(Q)$, by Lemma \ref{lemma2} we get that:
\begin{multline}\nonumber Q*_{\underline {\tau}}Z = Q_1\t_0\dots \t _0Q_{r-1}\t_0 (Q_r*_{\underline {\tau}\rq} (Z_1\t _0\dots \t _0 Z_{j_{\underline {\tau}}}))\t_0 Z_{j_{\underline {\tau}}+1}\t _0\dots \t_0 Z_s,\end{multline}
where $\underline {\tau}\rq = (\tau _0,\dots ,\tau _{j_{\underline {\tau}}})$ and $Q_r*_{\underline {\tau}\rq} (Z_1\t _0\dots \t _0 Z_{j_{\underline {\tau}}})$ is prime.

Applying $P*_{\lam}$, we obtain that:
\begin{multline}\nonumber
P*_{\lam} (Q*_{\underline {\tau}}Z) =\\
(P*_{\lam^1} (Q_1\t _0\dots \t _0 Q_{r-1}\t_0 (Q_r*_{\underline {\tau}\rq} (Z_1\t _0\dots \t _0 Z_{j_{\underline {\tau}}}))))*_{\lam^2} (Z_{j_{\underline{\tau}}+1}\t_0\dots \t _0 Z_s)=\\
((P*_{\lam ^1}Q)*_{\underline {\tau}^2}(Z_1\t_0 \dots \t _0 Z_{j_{\underline {\tau}}}))*_{\lam^2} (Z_{j_{\underline{\tau}}+1}\t_0\dots \t _0 Z_s)=\\
\hfill (P*_{\lam ^1}Q)*_{\underline {\delta}} Z,\end{multline}
for the weak compositions ${\lam^1} = (\lambda _0,\dots ,\lambda _{r-1},\lambda _r+\dots +\lambda _{r+s{-}j_{\underline {\tau}}})$, 

\noi $\lam^2 = (\lambda _r,\dots ,\lambda  _{r+s{-}j_{\underline {\tau}}})$, ${\underline {\tau}^2}= (\tau_0,\dots ,\tau_{j_{\underline {\tau}}-1}, \tau_{j_{\underline {\tau}}}+ \lambda _r+\dots +\lambda  _{r+s{-}j_{\underline {\tau}}})$ and
${\underline {\delta}} = (\tau_0,\dots ,\tau_{j_{\underline {\tau}}-1}, \tau_{j_{\underline {\tau}}}+ \lambda _r,\lambda _{r+1},\dots ,\lambda  _{r+s{-}j_{\underline {\tau}}})$.
\medskip

The formula above implies that for any pair $(\lam , {\underline {\tau}})\in \Lambda _{r+s{-}j_{\underline {\tau}}}^i(P)\t \Lambda _s^0(Q)$, the elements $P*_{\lam}(Q*_{\underline {\tau}}Z)$ and $(P*_{\underline {\gamma}} Q)*_{\underline {\delta}}Z$ are equal whenever 
$$\psi _i^1(P,Q)(\lam , {\underline {\tau}}) = ({\underline {\gamma}},{\underline {\delta}}).$$ 
So, we have proved that 
$$P*_i(Q*_0 Z) =\sum _{({\underline {\gamma}},{\underline {\delta}})} (P*_{\underline {\gamma}} Q) *_{\underline {\delta}} Z,$$
where the sum is taken over all 
${\underline {\gamma}}\in \displaystyle\coprod_{j=i}^m \Lambda^j_r(P)$ and 
${\underline {\delta}}\in\Lambda _s^i(P*_{\underline {\gamma}}Q)$ such that $\delta _s \leq \gamma _r$.
\medskip

$b)$ Suppose now that $(\lam ,{\underline {\tau}})$ belongs to $\Lambda _r^i(P)\t ({\displaystyle \coprod _{j=1}^i\Lambda _s^j(Q)})$. We have that:
$$Q*_{\underline {\tau}}Z = Q_1\t _0\dots \t _0 Q_{r-1}\t _0 (Q_r*_{\underline {\tau}} Z),$$
with $Q_1,\dots ,Q_{r-1}, Q_r*_{\underline {\tau}} Z$ prime. Let us compute
\begin{multline}\nonumber
 P*_{\lam} (Q*_{\underline {\tau}} Z)=\\
  (P*_{\lam^1}(Q_1\t_0\dots \t _0 Q_{r-1}))\t _{\lambda _r} (Q_r*_{\underline {\tau}} Z)=\\
 \hfill (P*_{\lam} Q)*_{\underline {\delta}} Z,\end{multline}
 where $\lam ^1 = (\lambda _0,\dots ,\lambda _{r-2},\lambda _{r-1}+\lambda _r)$ and ${\underline {\delta}} = (\tau _0,\dots, \tau_{s-1},\tau_s+\lambda _r).$
 \medskip
 
Using the notation of Proposition \ref{proprelac}, we have proved that 
$$P*_{\lam} (Q*_{\underline {\tau}} Z) = (P*_{\underline {\gamma}} Q)*_{\underline {\delta}} Z),$$ whenever 
 $\psi _i^2(P,Q)(\lam ,{\underline {\tau}}) = ({\underline {\gamma}}, {\underline {\delta}})$. So, we get:
 $$\sum _{j = 1}^i P*_i(Q*_j Z) = \sum _{({\underline {\gamma}},{\underline {\delta}})} (P*_{\underline {\gamma}} Q) *_{\underline {\delta}} Z,$$
where the sum is taken over all $({\underline {\gamma}},{\underline {\delta}})\in \Lambda _r^i(P)\t \Lambda _s^i(P*_{\underline {\gamma}}Q)$ such that $\delta _s > \gamma _r$.
\medskip

Adding up $a)$ and $b)$, we get that:
$$\sum _{j = 0}^i P*_i (Q*_j Z) = \sum _{j=i}^m (P*_jQ)*_i Z,$$
which ends the proof.
\end{proo}
\bigskip

\bigskip

\section{Connection with the $m$-Tamari lattice}
\medskip

For $n\geq 1$, let ${\mathcal Y}_n$ denotes the set of planar rooted binary trees with $n+1$ leaves.

\begin{notation} \label{trees} Define binary operations $\vee , /$ and $\backslash $ on the set of trees as follows: \begin{enumerate}
\item $\vee $ is the map which sends an ordered pair of trees $(t, w)$ to the tree obtained by joining the roots of $t$ and $w$ to a new root.
\item The element $t / w $ is the tree obtained by joining the root of $t$ to the first leaf of $w$.
\item The element $t\backslash w$ is the tree obtained by joining the root of $w$ to the last leaf of $t$.
\end{enumerate}
for  any $t$ and $w$ in ${\mathcal Y}_{\infty}:=\bigcup_{n\geq 0}{\mathcal Y}_n$. \end{notation}

The diagrams below show a more graphical description of the previous definitions,

\scalebox{0.9} 
{
\begin{pspicture}(0,-0.59919924)(11.354668,0.6191992)
\psline[linewidth=0.04cm](1.5967969,0.02080078)(1.9967968,-0.3791992)
\psline[linewidth=0.04cm](1.9967968,-0.3791992)(2.396797,0.02080078)
\psline[linewidth=0.04cm](1.9967968,-0.3791992)(1.9967968,-0.5791992)
\psline[linewidth=0.04cm](5.996797,0.02080078)(6.3967967,-0.3791992)
\psline[linewidth=0.04cm](6.3967967,-0.3791992)(6.796797,0.02080078)
\psline[linewidth=0.04cm](5.996797,0.22080079)(5.996797,0.02080078)
\psline[linewidth=0.04cm](6.3967967,-0.3791992)(6.3967967,-0.5791992)
\psline[linewidth=0.04cm](10.196796,0.02080078)(10.596797,-0.3791992)
\psline[linewidth=0.04cm](10.596797,-0.3791992)(10.996797,0.02080078)
\psline[linewidth=0.04cm](10.996797,0.22080079)(10.996797,0.02080078)
\psline[linewidth=0.04cm](10.596797,-0.3791992)(10.596797,-0.5791992)
\usefont{T1}{ptm}{m}{n}
\rput(1.0041308,-0.24919921){\large $t\vee w =$}
\usefont{T1}{ptm}{m}{n}
\rput(1.568252,0.24580078){$t$}
\usefont{T1}{ptm}{m}{n}
\rput(2.408252,0.28580078){$w$}
\usefont{T1}{ptm}{m}{n}
\rput(5.4741306,-0.24919921){\large $t/w =$}
\usefont{T1}{ptm}{m}{n}
\rput(5.968252,0.42580077){$t$}
\usefont{T1}{ptm}{m}{n}
\rput(6.388252,-0.05419922){$w$}
\usefont{T1}{ptm}{m}{n}
\rput(9.704131,-0.26919922){\large $t\backslash w =$}
\usefont{T1}{ptm}{m}{n}
\rput(10.588252,-0.07419922){$t$}
\usefont{T1}{ptm}{m}{n}
\rput(10.988252,0.4058008){$w$}
\end{pspicture} 
}
\medskip

Note that, adding ${\mathcal Y}_0 :=\{\vert \}$, for any $t\in {\mathcal Y}_n$ 
there exist unique trees $t^l\in {\mathcal Y}_{n_1}$ and $t^r\in {\mathcal Y}_{n_2}$ such that $t=t^l\vee t^r$.
\medskip

\begin{definition} \label{Tamariorder}  The {\it Tamari order } (see \cite{Tam}) on ${\mathcal Y}_n$, $n\geq 1$, is the partial order transitively spanned by the following relations:\begin{enumerate} 
\item $(t\vee w)\vee z < t\vee (w\vee z)$,
\item if $t < w$, then $t\vee z < w\vee z$,
\item if $w< z$, then $t\vee w < t\vee z$,\end{enumerate}
for $t, w, z \in {\mathcal Y}_{\infty}$. \end{definition}
\medskip

It is well-known that the set $\Dy_n^1$ of paths of size $n$ has the same cardinal that the set of planar binary rooted trees ${\mathcal Y}_n$.

Consider the map  $\Gamma _n: \Dy _n^1\longrightarrow{\mathcal Y}_n$, 
$n\geq 0$, defined by:\begin{enumerate}
\item $\Gamma _0(\bullet ) : = \vert$, is the unique element of ${\mathcal Y}_0$,
\item $\Gamma _n (P\t _0 Q) := \Gamma _{n_1}(P)/ \Gamma _{n_2}(Q)$, 
\item $\Gamma_{n_1+1}(\rho_1 \t _{1} P) = \vert \vee \Gamma_{n_1}(P)$,
\end{enumerate}
for any pair of Dyck paths $P\in \Dy_{n_1}^1$ and $Q\in \Dy_{n_2}^1$. The inverse application is defined recursively on $n$ by:\begin{enumerate}
\item $\Gamma _0^{-1}(\vert ) = \bullet $,
\item $\Gamma _n^{-1}(t^l\vee t^r) = \bigvee_u(\Gamma _{n_1}^{-1}(t^l),\Gamma _{n_2}^{-1}(t^r)$),\end{enumerate}
for any $t^l\in {\mathcal Y}_{n_1}$ and $t^r\in {\mathcal Y}_{n_2}$.
\medskip

So, the Tamari order is defined on $\Dy_n^1$, via the bijective map $\Gamma _n$, for $n\geq 1$. 
\medskip

F. Bergeron extended the Tamari order to the sets $\Dy _n^m$ of Dyck paths (see \cite{BerPre}) . Let us describe briefly the $m$-Tamari lattice $\Dy _n^m$.
\medskip

Let $P$ be an $m$-Dyck path. 
For any down step $d_0\in \dw (P)$ which is followed by an up step $u\in \up (P)$, consider the excursion $P_u$ of $u$ in $P$ and its matching down step $w_u$ as described in Definition \ref{excursion}. Let $P_{(d_0)}$ be the Dyck path obtained by removing $d_0$
and gluing the initial vertex of $u$ to the end of the step preceding $d_0$, and attaching $d_0$ at the final point of $w_u$. For example
\medskip

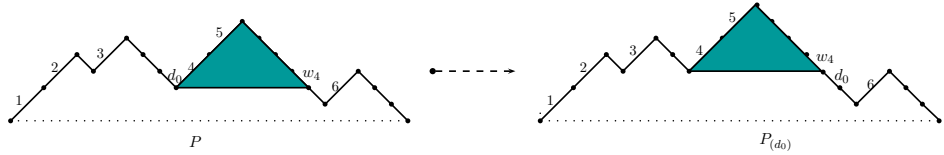
\scalebox{0.55} 
{
\begin{pspicture}(0,-1.8188477)(22.58,1.7988477)
\definecolor{color111b}{rgb}{0.0,0.6,0.6}
\psdots[dotsize=0.12](0.06,-1.0811523)
\psdots[dotsize=0.12](0.86,-0.28115234)
\psdots[dotsize=0.12](1.66,0.51884764)
\psdots[dotsize=0.12](2.06,0.11884765)
\psdots[dotsize=0.12](2.86,0.9188477)
\psdots[dotsize=0.12](3.26,0.51884764)
\psdots[dotsize=0.12](3.66,0.11884765)
\psdots[dotsize=0.12](4.06,-0.28115234)
\psdots[dotsize=0.12](4.86,0.51884764)
\psdots[dotsize=0.12](5.66,1.3188477)
\psdots[dotsize=0.12](6.06,0.9188477)
\psdots[dotsize=0.12](6.46,0.51884764)
\psdots[dotsize=0.12](6.86,0.11884765)
\psdots[dotsize=0.12](7.26,-0.28115234)
\psdots[dotsize=0.12](7.66,-0.68115234)
\psdots[dotsize=0.12](8.46,0.11884765)
\psdots[dotsize=0.12](8.86,-0.28115234)
\psdots[dotsize=0.12](9.26,-0.68115234)
\psdots[dotsize=0.12](9.66,-1.0811523)
\psline[linewidth=0.04cm](0.06,-1.0811523)(1.66,0.51884764)
\psline[linewidth=0.04cm](1.66,0.51884764)(2.06,0.11884765)
\psline[linewidth=0.04cm](2.06,0.11884765)(2.86,0.9188477)
\psline[linewidth=0.04cm](2.86,0.9188477)(3.26,0.51884764)
\psline[linewidth=0.04cm](3.26,0.51884764)(3.66,0.11884765)
\psline[linewidth=0.04cm](3.66,0.11884765)(4.06,-0.28115234)
\psline[linewidth=0.04cm](4.06,-0.28115234)(4.86,0.51884764)
\psline[linewidth=0.04cm](4.86,0.51884764)(5.66,1.3188477)
\psline[linewidth=0.04cm](5.66,1.3188477)(6.06,0.9188477)
\psline[linewidth=0.04cm](6.06,0.9188477)(6.46,0.51884764)
\psline[linewidth=0.04cm](6.46,0.51884764)(6.86,0.11884765)
\psline[linewidth=0.04cm](6.86,0.11884765)(7.26,-0.28115234)
\psline[linewidth=0.04cm](7.26,-0.28115234)(7.66,-0.68115234)
\psline[linewidth=0.04cm](7.66,-0.68115234)(8.46,0.11884765)
\psline[linewidth=0.04cm](8.46,0.11884765)(8.86,-0.28115234)
\psline[linewidth=0.04cm](8.86,-0.28115234)(9.26,-0.68115234)
\psline[linewidth=0.04cm](9.26,-0.68115234)(9.66,-1.0811523)
\psline[linewidth=0.04cm,linestyle=dotted,dotsep=0.16cm](9.66,-1.0811523)(0.06,-1.0811523)
\usefont{T1}{ptm}{m}{n}
\rput(0.26748535,-0.59615237){\small $1$}
\usefont{T1}{ptm}{m}{n}
\rput(1.1274854,0.22384766){\small $2$}
\usefont{T1}{ptm}{m}{n}
\rput(2.2274854,0.5438477){\small $3$}
\usefont{T1}{ptm}{m}{n}
\rput(4.4274855,0.22384766){\small $4$}
\usefont{T1}{ptm}{m}{n}
\rput(5.1074853,1.0638477){\small $5$}
\usefont{T1}{ptm}{m}{n}
\rput(7.9074855,-0.25615233){\small $6$}
\psline[linewidth=0.04cm,linestyle=dashed,dash=0.16cm 0.16cm,dotsize=0.07055555cm 2.0,arrowsize=0.05291667cm 2.0,arrowlength=1.4,arrowinset=0.4]{*->}(10.26,0.11884765)(12.26,0.11884765)
\usefont{T1}{ptm}{m}{n}
\rput(4.031455,0.0038476563){$d_0$}
\pspolygon[linewidth=0.04,fillstyle=solid,fillcolor=color111b](7.26,-0.28115234)(5.66,1.3188477)(4.06,-0.28115234)(7.26,-0.28115234)(7.26,-0.28115234)
\usefont{T1}{ptm}{m}{n}
\rput(7.331455,0.023847656){$w_4$}
\psdots[dotsize=0.12](12.86,-1.0811523)
\psdots[dotsize=0.12](13.66,-0.28115234)
\psdots[dotsize=0.12](14.46,0.51884764)
\psdots[dotsize=0.12](14.86,0.11884765)
\psdots[dotsize=0.12](15.66,0.9188477)
\psdots[dotsize=0.12](16.06,0.51884764)
\psdots[dotsize=0.12](16.46,0.11884765)
\psdots[dotsize=0.12](18.46,1.3188477)
\psdots[dotsize=0.12](17.26,0.9188477)
\psdots[dotsize=0.12](19.3,0.51884764)
\psdots[dotsize=0.12](18.86,0.9188477)
\psdots[dotsize=0.12](19.7,0.09884766)
\psdots[dotsize=0.12](20.1,-0.28115234)
\psdots[dotsize=0.12](20.5,-0.68115234)
\psdots[dotsize=0.12](21.3,0.11884765)
\psdots[dotsize=0.12](21.7,-0.28115234)
\psdots[dotsize=0.12](18.1,1.7188476)
\psdots[dotsize=0.12](22.5,-1.0811523)
\psdots[dotsize=0.12](22.1,-0.68115234)
\psline[linewidth=0.04cm,fillcolor=color111b,linestyle=dotted,dotsep=0.16cm](12.86,-1.0811523)(22.46,-1.0811523)
\psline[linewidth=0.04cm,fillcolor=color111b,linestyle=dotted,dotsep=0.16cm](12.86,-1.0811523)(12.86,-0.88115233)
\psline[linewidth=0.04cm,fillcolor=color111b](12.86,-1.0811523)(13.66,-0.28115234)
\psline[linewidth=0.04cm,fillcolor=color111b](13.66,-0.28115234)(14.46,0.51884764)
\psline[linewidth=0.04cm,fillcolor=color111b](14.46,0.51884764)(14.86,0.11884765)
\psline[linewidth=0.04cm,fillcolor=color111b](14.86,0.11884765)(15.66,0.9188477)
\psline[linewidth=0.04cm,fillcolor=color111b](15.66,0.9188477)(16.06,0.51884764)
\psline[linewidth=0.04cm,fillcolor=color111b](16.06,0.51884764)(16.46,0.11884765)
\psline[linewidth=0.04cm,fillcolor=color111b](16.46,0.11884765)(17.26,0.9188477)
\psline[linewidth=0.04cm,fillcolor=color111b](17.26,0.9188477)(18.06,1.7188476)
\usefont{T1}{ptm}{m}{n}
\rput(4.5214553,-1.5761523){$P$}
\usefont{T1}{ptm}{m}{n}
\rput(18.561455,-1.5961523){$P_{(d_0)}$}
\psline[linewidth=0.04cm,fillcolor=color111b](18.46,1.3188477)(18.06,1.7188476)
\psline[linewidth=0.04cm,fillcolor=color111b](18.86,0.9188477)(18.46,1.3188477)
\psline[linewidth=0.04cm,fillcolor=color111b](18.86,0.9188477)(19.26,0.51884764)
\psline[linewidth=0.04cm,fillcolor=color111b](19.74,0.058847655)(20.14,-0.34115234)
\psline[linewidth=0.04cm,fillcolor=color111b](20.1,-0.28115234)(20.5,-0.68115234)
\psline[linewidth=0.04cm,fillcolor=color111b](20.5,-0.68115234)(21.3,0.11884765)
\psline[linewidth=0.04cm,fillcolor=color111b](21.3,0.11884765)(21.7,-0.28115234)
\psline[linewidth=0.04cm,fillcolor=color111b](21.7,-0.28115234)(22.1,-0.68115234)
\psline[linewidth=0.04cm,fillcolor=color111b](22.1,-0.68115234)(22.5,-1.0611523)
\usefont{T1}{ptm}{m}{n}
\rput(13.107486,-0.5761523){\small $1$}
\usefont{T1}{ptm}{m}{n}
\rput(13.907485,0.22384766){\small $2$}
\usefont{T1}{ptm}{m}{n}
\rput(15.107486,0.62384766){\small $3$}
\usefont{T1}{ptm}{m}{n}
\rput(17.507484,1.4238477){\small $5$}
\usefont{T1}{ptm}{m}{n}
\rput(16.707485,0.62384766){\small $4$}
\usefont{T1}{ptm}{m}{n}
\rput(20.847485,-0.19615234){\small $6$}
\usefont{T1}{ptm}{m}{n}
\rput(19.751455,0.42384765){$w_4$}
\usefont{T1}{ptm}{m}{n}
\rput(20.151455,0.0038476563){$d_0$}
\pspolygon[linewidth=0.04,fillstyle=solid,fillcolor=color111b](19.66,0.11884765)(16.46,0.11884765)(18.06,1.7188476)
\end{pspicture} 
}

It is immediate to see that $\alpha _{P_{d_0}}(d) = \alpha _P(d)$, for any $d\in \dw (P)$.
\medskip

\begin{definition} \label{defmTamari}  The $m$-Tamari order on $\Dy ^m_n$ is the transitive relation spanned by the covering relation:
$$P \lessdot P_{(d)},$$
for any $d\in \dw (P)$ such that the final vertex of $d$ is the initial point of an up step $u\in \up (P)$. We use the symbol $\lessdot $ for a covering relation. \end{definition}
\medskip

The Hasse diagrams for $m = 2$ and $n = 1,2$ are:
\medskip

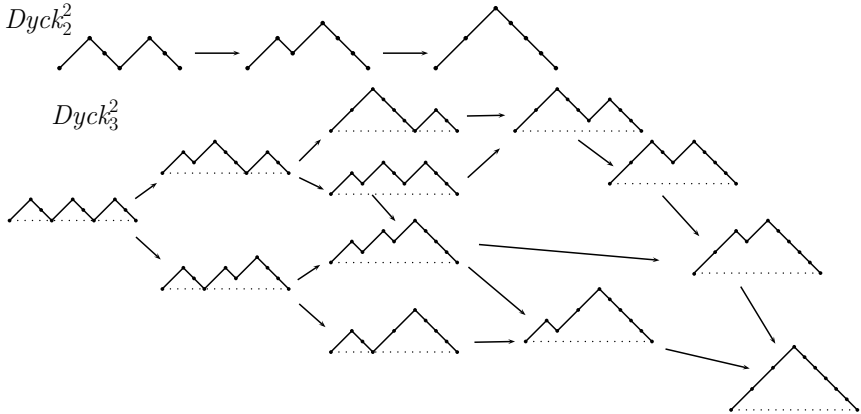
\scalebox{0.5} 
{
\begin{pspicture}(0,-0.88)(14.6575,0.88)
\psdots[dotsize=0.12](1.3775,-0.8)
\psdots[dotsize=0.12](2.1775,0.0)
\psdots[dotsize=0.12](2.5775,-0.4)
\psdots[dotsize=0.12](2.9775,-0.8)
\psdots[dotsize=0.12](3.7775,0.0)
\psdots[dotsize=0.12](4.1775,-0.4)
\psdots[dotsize=0.12](4.5775,-0.8)
\psdots[dotsize=0.12](6.3775,-0.8)
\psdots[dotsize=0.12](7.1775,0.0)
\psdots[dotsize=0.12](7.5775,-0.4)
\psdots[dotsize=0.12](8.3775,0.4)
\psdots[dotsize=0.12](8.7775,0.0)
\psdots[dotsize=0.12](9.1775,-0.4)
\psdots[dotsize=0.12](9.5775,-0.8)
\psdots[dotsize=0.12](11.3775,-0.8)
\psdots[dotsize=0.12](12.1775,0.0)
\psdots[dotsize=0.12](12.9775,0.8)
\psdots[dotsize=0.12](13.3775,0.4)
\psdots[dotsize=0.12](13.7775,0.0)
\psdots[dotsize=0.12](14.1775,-0.4)
\psdots[dotsize=0.12](14.5775,-0.8)
\psline[linewidth=0.04cm](1.3775,-0.8)(2.1775,0.0)
\psline[linewidth=0.04cm](2.1775,0.0)(2.5775,-0.4)
\psline[linewidth=0.04cm](2.5775,-0.4)(2.9775,-0.8)
\psline[linewidth=0.04cm](2.9775,-0.8)(3.7775,0.0)
\psline[linewidth=0.04cm](3.7775,0.0)(4.1775,-0.4)
\psline[linewidth=0.04cm](4.1775,-0.4)(4.5775,-0.8)
\psline[linewidth=0.04cm](6.3775,-0.8)(7.1775,0.0)
\psline[linewidth=0.04cm](7.1775,0.0)(7.5775,-0.4)
\psline[linewidth=0.04cm](7.5775,-0.4)(8.3775,0.4)
\psline[linewidth=0.04cm](8.3775,0.4)(8.7775,0.0)
\psline[linewidth=0.04cm](8.7775,0.0)(9.1775,-0.4)
\psline[linewidth=0.04cm](9.1775,-0.4)(9.5775,-0.8)
\psline[linewidth=0.04cm](11.3775,-0.8)(12.1775,0.0)
\psline[linewidth=0.04cm](12.1775,0.0)(12.9775,0.8)
\psline[linewidth=0.04cm](12.9775,0.8)(13.3775,0.4)
\psline[linewidth=0.04cm](13.3775,0.4)(13.7775,0.0)
\psline[linewidth=0.04cm](13.7775,0.0)(14.1775,-0.4)
\psline[linewidth=0.04cm](14.1775,-0.4)(14.5775,-0.8)
\usefont{T1}{ptm}{m}{n}
\rput(0.8171972,0.5){\huge $\Dy _2^2$}
\psline[linewidth=0.04cm,arrowsize=0.05291667cm 2.0,arrowlength=1.4,arrowinset=0.4]{->}(4.9775,-0.4)(6.1775,-0.4)
\psline[linewidth=0.04cm,arrowsize=0.05291667cm 2.0,arrowlength=1.4,arrowinset=0.4]{->}(9.9775,-0.4)(11.1775,-0.4)
\end{pspicture} 
}
\medskip

\scalebox{0.5} 
{
\begin{pspicture}(0,-4.34)(22.66,4.34)
\usefont{T1}{ptm}{m}{n}
\rput(2.061455,3.575){\huge $\Dy _3^2$}
\psdots[dotsize=0.1](0.05,0.77)
\psdots[dotsize=0.1](0.61,1.33)
\psdots[dotsize=0.1](0.89,1.05)
\psdots[dotsize=0.1](1.17,0.77)
\psdots[dotsize=0.1](1.73,1.33)
\psdots[dotsize=0.1](2.01,1.05)
\psdots[dotsize=0.1](2.29,0.77)
\psdots[dotsize=0.1](2.85,1.33)
\psdots[dotsize=0.1](3.13,1.05)
\psdots[dotsize=0.1](3.41,0.77)
\psline[linewidth=0.04cm](0.05,0.77)(0.61,1.33)
\psline[linewidth=0.04cm](0.61,1.33)(1.17,0.77)
\psline[linewidth=0.04cm](1.17,0.77)(1.73,1.33)
\psline[linewidth=0.04cm](1.73,1.33)(2.29,0.77)
\psline[linewidth=0.04cm](2.29,0.77)(2.85,1.33)
\psline[linewidth=0.04cm](2.85,1.33)(3.41,0.77)
\psdots[dotsize=0.1](4.11,2.03)
\psdots[dotsize=0.1](4.67,2.59)
\psdots[dotsize=0.1](4.95,2.31)
\psdots[dotsize=0.1](5.51,2.87)
\psdots[dotsize=0.1](5.79,2.59)
\psdots[dotsize=0.1](6.07,2.31)
\psdots[dotsize=0.1](6.35,2.03)
\psdots[dotsize=0.1](6.91,2.59)
\psdots[dotsize=0.1](7.19,2.31)
\psdots[dotsize=0.1](7.47,2.03)
\psline[linewidth=0.04cm](4.11,2.03)(4.67,2.59)
\psline[linewidth=0.04cm](4.67,2.59)(4.95,2.31)
\psline[linewidth=0.04cm](4.95,2.31)(5.51,2.87)
\psline[linewidth=0.04cm](5.51,2.87)(6.35,2.03)
\psline[linewidth=0.04cm](6.35,2.03)(6.91,2.59)
\psline[linewidth=0.04cm](6.91,2.59)(7.47,2.03)
\psdots[dotsize=0.1](4.11,-1.05)
\psdots[dotsize=0.1](4.67,-0.49)
\psdots[dotsize=0.1](4.95,-0.77)
\psdots[dotsize=0.1](5.23,-1.05)
\psdots[dotsize=0.1](5.79,-0.49)
\psdots[dotsize=0.1](6.07,-0.77)
\psdots[dotsize=0.1](6.63,-0.21)
\psdots[dotsize=0.1](6.91,-0.49)
\psdots[dotsize=0.1](7.19,-0.77)
\psdots[dotsize=0.1](7.47,-1.05)
\psline[linewidth=0.04cm](4.11,-1.05)(4.67,-0.49)
\psline[linewidth=0.04cm](4.67,-0.49)(5.23,-1.05)
\psline[linewidth=0.04cm](5.23,-1.05)(5.79,-0.49)
\psline[linewidth=0.04cm](5.79,-0.49)(6.07,-0.77)
\psline[linewidth=0.04cm](6.07,-0.77)(6.63,-0.21)
\psline[linewidth=0.04cm](6.63,-0.21)(7.47,-1.05)
\psdots[dotsize=0.1](8.59,3.15)
\psdots[dotsize=0.1](9.15,3.71)
\psdots[dotsize=0.1](9.71,4.27)
\psdots[dotsize=0.1](9.99,3.99)
\psdots[dotsize=0.1](10.27,3.71)
\psdots[dotsize=0.1](10.55,3.43)
\psdots[dotsize=0.1](10.83,3.15)
\psdots[dotsize=0.1](11.39,3.71)
\psdots[dotsize=0.1](11.67,3.43)
\psdots[dotsize=0.1](11.95,3.15)
\psline[linewidth=0.04cm](8.59,3.15)(9.71,4.27)
\psline[linewidth=0.04cm](9.71,4.27)(10.83,3.15)
\psline[linewidth=0.04cm](10.83,3.15)(11.39,3.71)
\psline[linewidth=0.04cm](11.39,3.71)(11.95,3.15)
\psdots[dotsize=0.1](8.59,1.47)
\psdots[dotsize=0.1](9.15,2.03)
\psdots[dotsize=0.1](9.43,1.75)
\psdots[dotsize=0.1](9.99,2.31)
\psdots[dotsize=0.1](10.27,2.03)
\psdots[dotsize=0.1](10.55,1.75)
\psdots[dotsize=0.1](11.11,2.31)
\psdots[dotsize=0.1](11.39,2.03)
\psdots[dotsize=0.1](11.67,1.75)
\psdots[dotsize=0.1](11.95,1.47)
\psline[linewidth=0.04cm](8.59,1.47)(9.15,2.03)
\psline[linewidth=0.04cm](9.15,2.03)(9.43,1.75)
\psline[linewidth=0.04cm](9.43,1.75)(9.99,2.31)
\psline[linewidth=0.04cm](9.99,2.31)(10.55,1.75)
\psline[linewidth=0.04cm](10.55,1.75)(11.11,2.31)
\psline[linewidth=0.04cm](11.11,2.31)(11.95,1.47)
\psdots[dotsize=0.1](8.59,-0.35)
\psdots[dotsize=0.1](9.15,0.21)
\psdots[dotsize=0.1](9.43,-0.07)
\psdots[dotsize=0.1](9.99,0.49)
\psdots[dotsize=0.1](10.27,0.21)
\psdots[dotsize=0.1](10.83,0.77)
\psdots[dotsize=0.1](11.11,0.49)
\psdots[dotsize=0.1](11.39,0.21)
\psdots[dotsize=0.1](11.67,-0.07)
\psdots[dotsize=0.1](11.95,-0.35)
\psline[linewidth=0.04cm](8.59,-0.35)(9.15,0.21)
\psline[linewidth=0.04cm](9.15,0.21)(9.43,-0.07)
\psline[linewidth=0.04cm](9.43,-0.07)(9.99,0.49)
\psline[linewidth=0.04cm](9.99,0.49)(10.27,0.21)
\psline[linewidth=0.04cm](10.27,0.21)(10.83,0.77)
\psline[linewidth=0.04cm](10.83,0.77)(11.95,-0.35)
\psdots[dotsize=0.1](8.59,-2.73)
\psdots[dotsize=0.1](9.15,-2.17)
\psdots[dotsize=0.1](9.43,-2.45)
\psdots[dotsize=0.1](9.71,-2.73)
\psdots[dotsize=0.1](10.27,-2.17)
\psdots[dotsize=0.1](10.83,-1.61)
\psdots[dotsize=0.1](11.11,-1.89)
\psdots[dotsize=0.1](11.39,-2.17)
\psdots[dotsize=0.1](11.67,-2.45)
\psdots[dotsize=0.1](11.95,-2.73)
\psline[linewidth=0.04cm](8.59,-2.73)(9.15,-2.17)
\psline[linewidth=0.04cm](9.15,-2.17)(9.71,-2.73)
\psline[linewidth=0.04cm](9.71,-2.73)(10.83,-1.61)
\psline[linewidth=0.04cm](10.83,-1.61)(11.95,-2.73)
\psdots[dotsize=0.1](13.49,3.15)
\psdots[dotsize=0.1](14.05,3.71)
\psdots[dotsize=0.1](14.61,4.27)
\psdots[dotsize=0.1](14.89,3.99)
\psdots[dotsize=0.1](15.17,3.71)
\psdots[dotsize=0.1](15.45,3.43)
\psdots[dotsize=0.1](16.01,3.99)
\psdots[dotsize=0.1](16.29,3.71)
\psdots[dotsize=0.1](16.57,3.43)
\psdots[dotsize=0.1](16.85,3.15)
\psline[linewidth=0.04cm](13.49,3.15)(14.61,4.27)
\psline[linewidth=0.04cm](14.61,4.27)(15.45,3.43)
\psline[linewidth=0.04cm](15.45,3.43)(16.01,3.99)
\psline[linewidth=0.04cm](16.01,3.99)(16.85,3.15)
\psline[linewidth=0.04cm,arrowsize=0.05291667cm 2.0,arrowlength=1.4,arrowinset=0.4]{->}(3.39,1.35)(3.97,1.77)
\psline[linewidth=0.04cm,arrowsize=0.05291667cm 2.0,arrowlength=1.4,arrowinset=0.4]{->}(3.41,0.33)(4.07,-0.27)
\psline[linewidth=0.04cm,arrowsize=0.05291667cm 2.0,arrowlength=1.4,arrowinset=0.4]{->}(7.75,2.35)(8.33,2.93)
\psline[linewidth=0.04cm,arrowsize=0.05291667cm 2.0,arrowlength=1.4,arrowinset=0.4]{->}(7.75,1.91)(8.41,1.59)
\psline[linewidth=0.04cm,arrowsize=0.05291667cm 2.0,arrowlength=1.4,arrowinset=0.4]{->}(7.71,-0.81)(8.31,-0.41)
\psline[linewidth=0.04cm,arrowsize=0.05291667cm 2.0,arrowlength=1.4,arrowinset=0.4]{->}(7.75,-1.45)(8.45,-2.07)
\psline[linewidth=0.04cm,arrowsize=0.05291667cm 2.0,arrowlength=1.4,arrowinset=0.4]{->}(12.21,3.55)(13.23,3.57)
\psline[linewidth=0.04cm,arrowsize=0.05291667cm 2.0,arrowlength=1.4,arrowinset=0.4]{->}(12.23,1.91)(13.09,2.69)
\psline[linewidth=0.04cm,arrowsize=0.05291667cm 2.0,arrowlength=1.4,arrowinset=0.4]{->}(9.69,1.45)(10.31,0.75)
\psdots[dotsize=0.1](14.33,-1.89)
\psdots[dotsize=0.1](14.61,-2.17)
\psdots[dotsize=0.1](15.17,-1.61)
\psdots[dotsize=0.1](15.73,-1.05)
\psdots[dotsize=0.1](16.01,-1.33)
\psdots[dotsize=0.1](16.29,-1.61)
\psdots[dotsize=0.1](16.57,-1.89)
\psdots[dotsize=0.1](16.85,-2.17)
\psdots[dotsize=0.1](17.13,-2.45)
\psline[linewidth=0.04cm](13.77,-2.45)(14.33,-1.89)
\psline[linewidth=0.04cm](14.33,-1.89)(14.61,-2.17)
\psline[linewidth=0.04cm](14.61,-2.17)(15.73,-1.05)
\psline[linewidth=0.04cm](15.73,-1.05)(17.13,-2.45)
\psline[linewidth=0.04cm,arrowsize=0.05291667cm 2.0,arrowlength=1.4,arrowinset=0.4]{->}(12.25,-0.47)(13.77,-1.75)
\psline[linewidth=0.04cm,arrowsize=0.05291667cm 2.0,arrowlength=1.4,arrowinset=0.4]{->}(12.41,-2.47)(13.49,-2.45)
\psdots[dotsize=0.1](16.01,1.75)
\psdots[dotsize=0.1](16.57,2.31)
\psdots[dotsize=0.1](17.13,2.87)
\psdots[dotsize=0.1](17.41,2.59)
\psdots[dotsize=0.1](17.69,2.31)
\psdots[dotsize=0.1](18.25,2.87)
\psdots[dotsize=0.1](18.53,2.59)
\psdots[dotsize=0.1](18.81,2.31)
\psdots[dotsize=0.1](19.09,2.03)
\psdots[dotsize=0.1](19.37,1.75)
\psline[linewidth=0.04cm](16.01,1.75)(17.13,2.87)
\psline[linewidth=0.04cm](17.13,2.87)(17.69,2.31)
\psline[linewidth=0.04cm](17.69,2.31)(18.25,2.87)
\psline[linewidth=0.04cm](18.25,2.87)(19.37,1.75)
\psdots[dotsize=0.1](18.25,-0.65)
\psdots[dotsize=0.1](18.81,-0.07)
\psdots[dotsize=0.1](19.65,0.21)
\psdots[dotsize=0.1](20.21,0.77)
\psdots[dotsize=0.1](20.49,0.49)
\psdots[dotsize=0.1](20.77,0.21)
\psdots[dotsize=0.1](21.05,-0.07)
\psdots[dotsize=0.1](21.33,-0.35)
\psdots[dotsize=0.1](19.79,-3.71)
\psdots[dotsize=0.1](20.35,-3.15)
\psdots[dotsize=0.1](20.91,-2.59)
\psdots[dotsize=0.1](19.23,-4.27)
\psdots[dotsize=0.1](21.19,-2.87)
\psdots[dotsize=0.1](21.47,-3.15)
\psdots[dotsize=0.1](21.75,-3.43)
\psdots[dotsize=0.1](22.03,-3.71)
\psdots[dotsize=0.1](22.31,-3.99)
\psdots[dotsize=0.1](22.59,-4.27)
\psline[linewidth=0.04cm](19.23,-4.27)(20.91,-2.59)
\psline[linewidth=0.04cm](20.91,-2.59)(22.59,-4.27)
\psline[linewidth=0.04cm,linestyle=dotted,dotsep=0.16cm](0.05,0.77)(3.41,0.77)
\psline[linewidth=0.04cm,linestyle=dotted,dotsep=0.16cm](4.11,-1.05)(7.47,-1.05)
\psline[linewidth=0.04cm,linestyle=dotted,dotsep=0.16cm](4.11,2.03)(7.47,2.03)
\psline[linewidth=0.04cm,linestyle=dotted,dotsep=0.16cm](8.59,3.15)(11.95,3.15)
\psline[linewidth=0.04cm,linestyle=dotted,dotsep=0.16cm](8.59,1.47)(11.95,1.47)
\psline[linewidth=0.04cm,linestyle=dotted,dotsep=0.16cm](8.59,-0.35)(11.95,-0.35)
\psline[linewidth=0.04cm,linestyle=dotted,dotsep=0.16cm](8.59,-2.73)(11.95,-2.73)
\psline[linewidth=0.04cm,linestyle=dotted,dotsep=0.16cm](13.77,-2.45)(17.13,-2.45)
\psline[linewidth=0.04cm,linestyle=dotted,dotsep=0.16cm](13.49,3.15)(16.85,3.15)
\psline[linewidth=0.04cm,linestyle=dotted,dotsep=0.16cm](16.01,1.75)(19.37,1.75)
\psline[linewidth=0.04cm,linestyle=dotted,dotsep=0.16cm](19.23,-4.27)(22.59,-4.27)
\psdots[dotsize=0.1](13.77,-2.45)
\psline[linewidth=0.04cm,arrowsize=0.05291667cm 2.0,arrowlength=1.4,arrowinset=0.4]{->}(15.17,2.91)(16.11,2.21)
\psdots[dotsize=0.1](19.37,0.49)
\psdots[dotsize=0.1](21.61,-0.63)
\psline[linewidth=0.04cm](18.25,-0.63)(19.37,0.49)
\psline[linewidth=0.04cm](19.37,0.49)(19.65,0.21)
\psline[linewidth=0.04cm](19.65,0.21)(20.21,0.77)
\psline[linewidth=0.04cm](20.21,0.77)(21.61,-0.63)
\psline[linewidth=0.04cm,linestyle=dotted,dotsep=0.16cm](18.25,-0.63)(21.61,-0.63)
\psline[linewidth=0.04cm,arrowsize=0.05291667cm 2.0,arrowlength=1.4,arrowinset=0.4]{->}(12.53,0.13)(17.33,-0.29)
\psline[linewidth=0.04cm,arrowsize=0.05291667cm 2.0,arrowlength=1.4,arrowinset=0.4]{->}(17.41,1.43)(18.39,0.35)
\psline[linewidth=0.04cm,arrowsize=0.05291667cm 2.0,arrowlength=1.4,arrowinset=0.4]{->}(19.49,-0.99)(20.41,-2.51)
\psline[linewidth=0.04cm,arrowsize=0.05291667cm 2.0,arrowlength=1.4,arrowinset=0.4]{->}(17.49,-2.65)(19.73,-3.17)
\end{pspicture} 
}

\bigskip

For $m=1$, it is easy to see that the order defined on $\Dy _n^1$ in \cite{BerPre} 
is the order induced by the Tamari order on ${\mathcal Y}_n$ via the map $\Gamma _n^{-1}$.
 That is, $\Gamma _n$ is an isomorphism of partially ordered sets, for $n\geq 1$.
\bigskip

The goal of the present section is to show that the binary operations $*_i: \K[\Dy _n^m]\ot \K[\Dy_r^m]\longrightarrow \K[\Dy _{n+r}^m]$ 
are described in terms of the $m$-Tamari order. Let us begin by describing the situation in the case $m = 1$.
\medskip

\begin{definition} \label{defdendriform} (see \cite{Lod}) A {\it dendriform algebra} over $\K$ is a vector space $A$ equipped with binary operations $\succ $ and $\prec $ satisfying the following conditions\begin{enumerate}
\item $x\succ (y\succ z) = (x\succ y + x\prec y) \succ z$,
\item $x\succ (y\prec z) = (x\succ y)\prec z $, 
\item $x\prec (y\succ z + y \prec z)=(x\prec y)\prec z$,\end{enumerate}
for $x,y,z \in A$.\end{definition}

In \cite{LodRon}, J.-L. Loday and the third author showed that the vector space $\K[{\mathcal Y}_{\infty}]$, spanned by $\bigcup _{n\geq 1} {\mathcal Y}_n$, may be endowed with a natural dendriform structure, in such a way that $\K[{\mathcal Y}_{\infty}]$ is the free dendriform algebra on one generator. 
\medskip

The dendriform structure on $\K[{\mathcal Y}_{\infty}]$ is described in terms of the Tamari order and the binary operations $/$ and $\backslash $ (see \cite{LodRon1}) as follows: \begin{enumerate}
\item ${\displaystyle t \succ w = \sum _{t/w \leq z\leq (t\backslash w^l)\vee w^r} z}$,
\item ${\displaystyle t\prec w = \sum _{t^l\vee (t^r/w)\leq z \leq t\backslash w} z }$.\end{enumerate}
\bigskip

It is not difficult to see that, for $m=1$, we have:\begin{enumerate}
\item $\Gamma _n(x) \succ \Gamma _r(y) = \Gamma _n(x) *_0\Gamma _r(y)$,
\item $\Gamma _n(x) \prec \Gamma _r(y) = \Gamma _n(x) *_1\Gamma _r(y)$,\end{enumerate}
 for any pair of elements $x\in \Dy_n^1$ and $y\in \Dy_r^1$.
\medskip

\begin{remark} \label{remmatching} \begin{enumerate} \item Let $Q$ be a prime Dyck path, for any pair of $m$-Dyck path $P$, we get:
$$ P\t _0 Q < P\t _1 Q <\dots < P\t _{L(P)}Q,$$
in the $m$-Tamari lattice.
\item If $P < P\rq $ in $\Dy _{n_1}^m$ are such that $L(P) = L(P\rq)$, and $Q < Q\rq $ in $\Dy_{n_2}^m$, then \begin{enumerate}
\item $P\t _k Q < P\t _k Q\rq $, for any $0\leq k\leq L(P)$, 
\item $P\t _k Q < P\rq \t_k Q$, for any $0\leq k\leq L(P)$.\end{enumerate}
\end{enumerate}
\end{remark}
\medskip

For the rest of the section, the $m$-Dyck path $Q$ is supposed to be a product $Q=Q_0\t _0\dots \t_0 Q_r$, where all the $Q_j$\rq s are prime Dyck paths.

\begin{lemma} \label{lemmalambdaordre} Let $P\in \Dy _{n_1}^m$ and $Q\in \Dy _{n_2}^m$ be two Dyck paths. Two weak compositions $\lam$ and ${\underline {\gamma}}$ in $\Lambda_r(P)$ satisfy that 
$$\lambda _j + \dots +\lambda _r \leq \gamma _j + \dots + \gamma _r,$$
for $1\leq j\leq r$, if, and only if, $P*_{\lam}Q \leq P*_{\underline {\gamma}} Q$.\end{lemma}
\medskip

\begin{proo} If $Q$ is prime, the result follows from point $(1)$ of Remark \ref{remmatching}. 
Suppose that $Q = Q_1\t _0\dots \t_0 Q_r$, for $r >1$. 
A recursive argument shows that, for any pair of elements $\lam\rq$ and ${\underline {\gamma }\rq}$ in $\Lambda _{r-1}(P)$,  we have that $$P*_{\lam\rq} (Q_1\t_0 \dots \t_0 Q_{r-1}) \leq P*_{\underline {\gamma}\rq} (Q_1\t_0 \dots \t_0 Q_{r-1}),$$
 whenever $\lambda\rq _j + \dots +\lambda\rq  _{r-1} \leq \gamma\rq  _j + \dots + \lambda\rq _{r-1},$ for $1\leq j\leq r-1$.
\medskip 
 
We have \begin{enumerate}
\item $P*_{\lam}Q = (P*_{\lam\rq} (Q_1\t _0\dots \t_0 Q_{r-1}))\t_{\lambda _r} Q_r$, 

\item $P*_{\underline {\gamma}} Q = (P*_{\underline {\gamma}\rq} (Q_1\t _0\dots \t_0 Q_{r-1}))\t_{\gamma _r} Q_r$, 

\end{enumerate}
where $\lam \rq = (\lambda _0, \dots ,\lambda _{r-1}, \lambda _{r-1} + \lambda _r)$ and 
$\underline {\gamma}\rq = (\gamma _0, \dots ,\gamma _{r-1}, \gamma _{r-1} + \gamma_r)$.
By the recursive hypothesis, we get that 
$$P*_{\lam\rq} (Q_1\t _0\dots \t_0 Q_{r-1})\leq P*_{\underline {\gamma}\rq} (Q_1\t _0\dots \t_0 Q_{r-1}),$$
 and using that $\lambda _r \leq \gamma _r$ we finally obtain $P*_{\lam}Q \leq P*_{\underline {\gamma}} Q$.
\medskip

Conversely, suppose that $P*_{\lam}Q \leq P*_{\underline {\gamma}} Q$. Point $(3)$ of Remark \ref{remmatching} implies that $$\lambda _j + \dots +\lambda _r \leq \gamma _j + \dots + \lambda _r,$$
for $1\leq j\leq r$, which ends the proof. \end{proo}
\medskip

\begin{notation} \label{notTamari} For any $m$-Dyck path $P$ of size $n$ and any $0\leq i\leq m$, let \begin{enumerate}
\item $c_i(P)$ be the minimal number of elements such that the word 
$$\alpha _P(d_{L(P){-}c_i(P)+1})\dots \alpha _P(d_{L(P)})$$ contains $i$ times an integer in $\{1,\dots ,n\}$ and no integer more than $i$ times,
\item $C_i(P)$ be the maximal integer such that the word $$\alpha _P(d_{L(P){-}C_i(P)+1})\dots \alpha _P(d_{L(P)})$$ contains at least one integer repeated $i$ times and no integer repeated $i+1$ times.\end{enumerate}
\medskip

Let $P\in \Dy _{n_1}^m$ and $Q\in \Dy _{n_2}^m$ be two Dyck paths. For any integer $0\leq i\leq m$, let $P /_i Q$ and $P \backslash _i Q$ be the Dyck paths defined as follows:\begin{enumerate}
\item $P/_i Q := P\t _{c_i(P)} Q$,
\item $P\backslash _i Q := (P\t _{L(P)} (Q_1\t_0\dots \t _0 Q_{r-1}))\t _{C_i(P)} Q_r$.\end{enumerate}
\end{notation}
\medskip

\begin{proposition} \label{propTamarim} For any pair of Dyck paths $P\in \Dy _{n_1}^m$ and $Q\in \Dy _{n_2}^m$ and any integer $0\leq i\leq m$, the product $*_i$ is given in terms of the $m$-Tamari order by the following formula:
$$P*_iQ = \sum _{P/_iQ\leq Z\leq P\backslash _i Q} Z.$$\end{proposition}
\medskip

\begin{proo} Suppose that $Q=Q_1\t_0\dots \t_0Q_r$, with all the $Q_i$\rq s prime and that $\lam \in \Lambda _r^i(P)$. 

The weak composition $\lam = (\lambda _0,\dots ,\lambda _r)$ satisfies that $c_i(P)\leq \lambda _r\leq C_i(P)$ and ${\displaystyle \sum _{j=0}^r\lambda _i = L(P)}$. 
\medskip

As \begin{itemize}
\item $P/_iQ = P*_{(L(P){-}c_i(P), 0, \dots ,0, c_i(P))} Q$, and 
\item $P\backslash _i Q = P*_{(0,\dots , 0, L(P){-}C_i(P), C_i(P))}Q$,\end{itemize}
 applying Lemma \ref{lemmalambdaordre}, it is easily seen that $P/_i Q\leq P*_{\lam}Q\leq P\backslash _i Q$.
\medskip

Recall that, whenever $R < S$ in the Tamari lattice, the set $\dw (R)$ of down steps of $R$ is identified with the set $\dw  (S)$. For any $d\in \dw (P)$ the levels of $d$ in $R$ and in $S$ are different but $\alpha _R(d) = \alpha _S(d)$. 

Note that the unique down steps which have different levels in the Dyck paths $P/_iQ$ and $P\backslash _i Q$ are colored by the set of integers $\{ 1,\dots ,n_1\}$. So, for any $P/_iQ\leq Z\leq P\backslash _i Q$ and any $1\leq l\leq r$, we get that 
$$L_j(Z) =L_j(Q_l),\ {\rm for}\ n_1+n_{21}+\dots +n_{2(l-1)} < j < n_1+n_{21}+\dots + n_{2l}. \leqno (*)$$

Define $$\lambda _j = \begin{cases} L_{n_1+n_{21}+\dots +n_{2j}} (Z) -L(Q_j),&\ {\rm for}\ 1\leq j\leq r,\\
L_{n_1}(Z) - L(P),&\ {\rm for}\ j=0.\end{cases} $$

The arguments above show that \begin{enumerate}
\item $c_i\leq \lambda _r\leq C_i$,
\item $c_i \leq  \lambda _j+ \dots +\lambda _r \leq L(P)$, for $1\leq j\leq r-1$,
\item $0\leq L_{n_1}(Z)\leq L(P) - c_i.$\end{enumerate}

From $(*)$, we get that $Z = P*_{\lam} Q$.
\medskip

Lemma \ref{lemmalambdaordre} and $P/_iQ\leq P*_{\lam}Q\leq P\backslash _i Q$ imply that $\lam \in \Lambda _r^i(P)$.\end{proo}
\medskip

Let us define the product $*$ on $\K[\Dy ^m]$ as the sum $* := {\displaystyle \sum _{i=0}^m *_i}$. It is not difficult to see, using Proposition \ref{propTamarim}, that 
$$P*Q = \sum _{P/_0 Q \leq Z\leq P\backslash _mQ} Z.$$
\medskip

\begin{example} Consider the Dyck paths $P = (1,3)$ and $Q=(2,2)$ in $\Dy _2^2$, the following diagram describes the Tamari interval $I_{P*Q}$ of all $Z\in \Dy _4^2$ such that 
$P * Q = {\displaystyle \sum _{Z\in I_{P*Q}} Z}$.

The Dyck paths in red are the terms of $P*_0Q$, the ones in green are the terms of $P*_1Q$, 
and the ones in blue are the terms appearing in $P*_2Q$.
\medskip

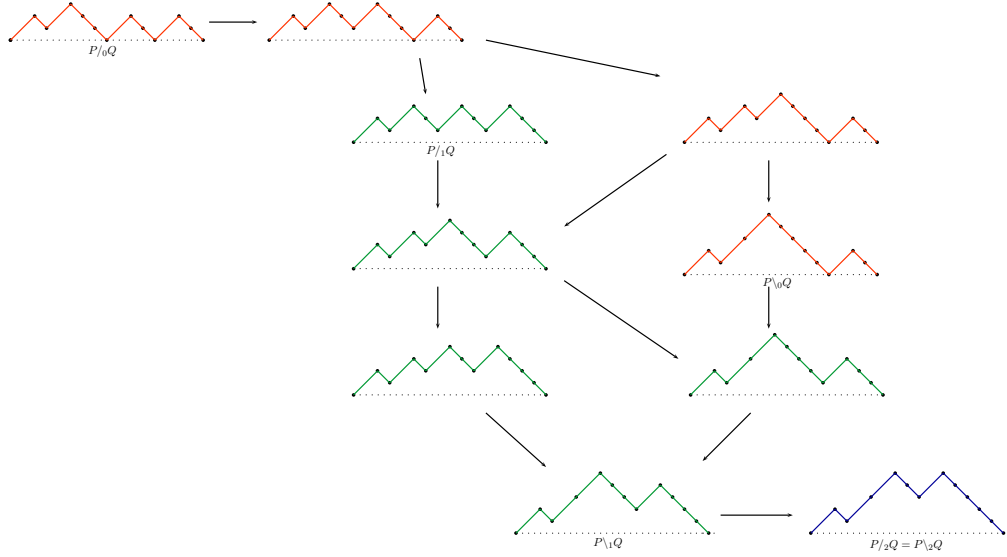
\scalebox{0.4} 
{
\begin{pspicture}(0,-9.158848)(33.14,9.138847)
\definecolor{color228}{rgb}{1.0,0.2,0.0}
\definecolor{color258}{rgb}{0.0,0.6,0.2}
\definecolor{color294}{rgb}{0.0,0.0,0.6}
\psdots[dotsize=0.12](0.06,7.8588476)
\psdots[dotsize=0.12](0.86,8.658848)
\psdots[dotsize=0.12](1.26,8.258847)
\psdots[dotsize=0.12](2.06,9.058847)
\psdots[dotsize=0.12](2.46,8.658848)
\psdots[dotsize=0.12](2.86,8.258847)
\psdots[dotsize=0.12](3.26,7.8588476)
\psdots[dotsize=0.12](4.06,8.658848)
\psdots[dotsize=0.12](4.46,8.258847)
\psdots[dotsize=0.12](4.86,7.8588476)
\psdots[dotsize=0.12](5.66,8.658848)
\psdots[dotsize=0.12](6.06,8.258847)
\psdots[dotsize=0.12](6.46,7.8588476)
\psdots[dotsize=0.12](8.66,7.8588476)
\psdots[dotsize=0.12](9.46,8.658848)
\psdots[dotsize=0.12](9.86,8.258847)
\psdots[dotsize=0.12](10.66,9.058847)
\psdots[dotsize=0.12](11.06,8.658848)
\psdots[dotsize=0.12](11.46,8.258847)
\psdots[dotsize=0.12](12.26,9.058847)
\psdots[dotsize=0.12](12.66,8.658848)
\psdots[dotsize=0.12](13.06,8.258847)
\psdots[dotsize=0.12](13.46,7.8588476)
\psdots[dotsize=0.12](14.26,8.658848)
\psdots[dotsize=0.12](14.66,8.258847)
\psdots[dotsize=0.12](15.06,7.8588476)
\psdots[dotsize=0.12](12.26,5.2588477)
\psdots[dotsize=0.12](12.66,4.8588476)
\psdots[dotsize=0.12](13.46,5.658848)
\psdots[dotsize=0.12](13.86,5.2588477)
\psdots[dotsize=0.12](14.26,4.8588476)
\psdots[dotsize=0.12](15.06,5.658848)
\psdots[dotsize=0.12](15.46,5.2588477)
\psdots[dotsize=0.12](15.86,4.8588476)
\psdots[dotsize=0.12](16.66,5.658848)
\psdots[dotsize=0.12](17.06,5.2588477)
\psdots[dotsize=0.12](17.46,4.8588476)
\psdots[dotsize=0.12](17.86,4.4588475)
\psdots[dotsize=0.12](22.46,4.4588475)
\psdots[dotsize=0.12](23.26,5.2588477)
\psdots[dotsize=0.12](23.66,4.8588476)
\psdots[dotsize=0.12](24.46,5.658848)
\psdots[dotsize=0.12](24.86,5.2588477)
\psdots[dotsize=0.12](25.66,6.0588474)
\psdots[dotsize=0.12](26.06,5.658848)
\psdots[dotsize=0.12](26.46,5.2588477)
\psdots[dotsize=0.12](26.86,4.8588476)
\psdots[dotsize=0.12](27.26,4.4588475)
\psdots[dotsize=0.12](28.06,5.2588477)
\psdots[dotsize=0.12](28.46,4.8588476)
\psdots[dotsize=0.12](28.86,4.4588475)
\psdots[dotsize=0.12](11.46,0.25884765)
\psdots[dotsize=0.12](12.26,1.0588477)
\psdots[dotsize=0.12](12.66,0.65884763)
\psdots[dotsize=0.12](13.46,1.4588476)
\psdots[dotsize=0.12](13.86,1.0588477)
\psdots[dotsize=0.12](14.66,1.8588476)
\psdots[dotsize=0.12](15.06,1.4588476)
\psdots[dotsize=0.12](15.46,1.0588477)
\psdots[dotsize=0.12](15.86,0.65884763)
\psdots[dotsize=0.12](16.66,1.4588476)
\psdots[dotsize=0.12](17.06,1.0588477)
\psdots[dotsize=0.12](17.46,0.65884763)
\psdots[dotsize=0.12](17.86,0.25884765)
\psdots[dotsize=0.12](22.46,0.058847655)
\psdots[dotsize=0.12](23.26,0.8588477)
\psdots[dotsize=0.12](23.66,0.45884764)
\psdots[dotsize=0.12](24.46,1.2588477)
\psdots[dotsize=0.12](25.26,2.0588477)
\psdots[dotsize=0.12](25.66,1.6588477)
\psdots[dotsize=0.12](26.06,1.2588477)
\psdots[dotsize=0.12](26.46,0.8588477)
\psdots[dotsize=0.12](26.86,0.45884764)
\psdots[dotsize=0.12](27.26,0.058847655)
\psdots[dotsize=0.12](28.06,0.8588477)
\psdots[dotsize=0.12](28.46,0.45884764)
\psdots[dotsize=0.12](28.86,0.058847655)
\psdots[dotsize=0.12](11.46,4.4588475)
\psdots[dotsize=0.12](11.46,-3.9411523)
\psdots[dotsize=0.12](12.26,-3.1411524)
\psdots[dotsize=0.12](12.66,-3.5411522)
\psdots[dotsize=0.12](13.46,-2.7411523)
\psdots[dotsize=0.12](13.86,-3.1411524)
\psdots[dotsize=0.12](14.66,-2.3411524)
\psdots[dotsize=0.12](15.06,-2.7411523)
\psdots[dotsize=0.12](15.46,-3.1411524)
\psdots[dotsize=0.12](16.26,-2.3411524)
\psdots[dotsize=0.12](16.66,-2.7411523)
\psdots[dotsize=0.12](17.06,-3.1411524)
\psdots[dotsize=0.12](17.46,-3.5411522)
\psdots[dotsize=0.12](17.86,-3.9411523)
\psdots[dotsize=0.12](22.66,-3.9411523)
\psdots[dotsize=0.12](23.46,-3.1411524)
\psdots[dotsize=0.12](23.86,-3.5411522)
\psdots[dotsize=0.12](24.66,-2.7411523)
\psdots[dotsize=0.12](25.46,-1.9411523)
\psdots[dotsize=0.12](25.86,-2.3411524)
\psdots[dotsize=0.12](26.26,-2.7411523)
\psdots[dotsize=0.12](26.66,-3.1411524)
\psdots[dotsize=0.12](27.06,-3.5411522)
\psdots[dotsize=0.12](27.86,-2.7411523)
\psdots[dotsize=0.12](28.26,-3.1411524)
\psdots[dotsize=0.12](28.66,-3.5411522)
\psdots[dotsize=0.12](29.06,-3.9411523)
\psdots[dotsize=0.12](16.86,-8.541152)
\psdots[dotsize=0.12](17.66,-7.7411523)
\psdots[dotsize=0.12](18.06,-8.141152)
\psdots[dotsize=0.12](18.86,-7.341152)
\psdots[dotsize=0.12](19.66,-6.5411525)
\psdots[dotsize=0.12](20.06,-6.9411526)
\psdots[dotsize=0.12](20.46,-7.341152)
\psdots[dotsize=0.12](20.86,-7.7411523)
\psdots[dotsize=0.12](21.66,-6.9411526)
\psdots[dotsize=0.12](22.06,-7.341152)
\psdots[dotsize=0.12](22.46,-7.7411523)
\psdots[dotsize=0.12](22.86,-8.141152)
\psdots[dotsize=0.12](23.26,-8.541152)
\psdots[dotsize=0.12](26.66,-8.541152)
\psdots[dotsize=0.12](27.46,-7.7411523)
\psdots[dotsize=0.12](27.86,-8.141152)
\psdots[dotsize=0.12](28.66,-7.341152)
\psdots[dotsize=0.12](29.46,-6.5411525)
\psdots[dotsize=0.12](29.86,-6.9411526)
\psdots[dotsize=0.12](30.26,-7.341152)
\psdots[dotsize=0.12](31.06,-6.5411525)
\psdots[dotsize=0.12](31.46,-6.9411526)
\psdots[dotsize=0.12](31.86,-7.341152)
\psdots[dotsize=0.12](32.26,-7.7411523)
\psdots[dotsize=0.12](32.66,-8.141152)
\psdots[dotsize=0.12](33.06,-8.541152)
\psdots[dotsize=0.12](0.06,7.8588476)
\psline[linewidth=0.04cm,linestyle=dotted,dotsep=0.16cm](0.06,7.8588476)(6.46,7.8588476)
\psline[linewidth=0.04cm,linestyle=dotted,dotsep=0.16cm](8.66,7.8588476)(15.06,7.8588476)
\psline[linewidth=0.04cm,linestyle=dotted,dotsep=0.16cm](22.46,4.4588475)(28.86,4.4588475)
\psline[linewidth=0.04cm,linestyle=dotted,dotsep=0.16cm](11.46,4.4588475)(18.06,4.4588475)
\psline[linewidth=0.04cm,linestyle=dotted,dotsep=0.16cm](11.46,0.25884765)(17.86,0.25884765)
\psline[linewidth=0.04cm,linestyle=dotted,dotsep=0.16cm](22.46,0.058847655)(29.06,0.058847655)
\psline[linewidth=0.04cm,linestyle=dotted,dotsep=0.16cm](11.46,-3.9411523)(18.06,-3.9411523)
\psline[linewidth=0.04cm,linestyle=dotted,dotsep=0.16cm](22.66,-3.9411523)(29.06,-3.9411523)
\psline[linewidth=0.04cm,linestyle=dotted,dotsep=0.16cm](16.86,-8.541152)(22.46,-8.541152)
\psline[linewidth=0.04cm,linestyle=dotted,dotsep=0.16cm](22.46,-8.541152)(23.06,-8.541152)
\psline[linewidth=0.04cm,linestyle=dotted,dotsep=0.16cm](23.06,-8.541152)(23.46,-8.541152)
\psline[linewidth=0.04cm,linestyle=dotted,dotsep=0.16cm](26.66,-8.541152)(33.06,-8.541152)
\psline[linewidth=0.04cm,arrowsize=0.05291667cm 2.0,arrowlength=1.4,arrowinset=0.4]{->}(6.66,8.458848)(8.26,8.458848)
\psline[linewidth=0.04cm,arrowsize=0.05291667cm 2.0,arrowlength=1.4,arrowinset=0.4]{->}(13.66,7.2588477)(13.86,6.0588474)
\psline[linewidth=0.04cm,arrowsize=0.05291667cm 2.0,arrowlength=1.4,arrowinset=0.4]{->}(15.86,7.8588476)(21.66,6.658848)
\psline[linewidth=0.04cm,arrowsize=0.05291667cm 2.0,arrowlength=1.4,arrowinset=0.4]{->}(14.26,3.8588476)(14.26,2.2588477)
\psline[linewidth=0.04cm,arrowsize=0.05291667cm 2.0,arrowlength=1.4,arrowinset=0.4]{->}(25.26,3.8588476)(25.26,2.4588478)
\psline[linewidth=0.04cm,arrowsize=0.05291667cm 2.0,arrowlength=1.4,arrowinset=0.4]{->}(21.86,4.0588474)(18.46,1.6588477)
\psline[linewidth=0.04cm,arrowsize=0.05291667cm 2.0,arrowlength=1.4,arrowinset=0.4]{->}(14.26,-0.34115234)(14.26,-1.7411523)
\psline[linewidth=0.04cm,arrowsize=0.05291667cm 2.0,arrowlength=1.4,arrowinset=0.4]{->}(18.46,-0.14115234)(22.26,-2.7411523)
\psline[linewidth=0.04cm,arrowsize=0.05291667cm 2.0,arrowlength=1.4,arrowinset=0.4]{->}(25.26,-0.34115234)(25.26,-1.7411523)
\psline[linewidth=0.04cm,arrowsize=0.05291667cm 2.0,arrowlength=1.4,arrowinset=0.4]{->}(24.66,-4.5411525)(23.06,-6.1411524)
\psline[linewidth=0.04cm,arrowsize=0.05291667cm 2.0,arrowlength=1.4,arrowinset=0.4]{->}(15.86,-4.5411525)(17.86,-6.341152)
\psline[linewidth=0.04cm,arrowsize=0.05291667cm 2.0,arrowlength=1.4,arrowinset=0.4]{->}(23.66,-7.9411526)(26.06,-7.9411526)
\psline[linewidth=0.04cm,linecolor=color228](0.06,7.8588476)(0.86,8.658848)
\psline[linewidth=0.04cm,linecolor=color228](0.86,8.658848)(1.26,8.258847)
\psline[linewidth=0.04cm,linecolor=color228](1.26,8.258847)(2.06,9.058847)
\psline[linewidth=0.04cm,linecolor=color228](2.06,9.058847)(3.26,7.8588476)
\psline[linewidth=0.04cm,linecolor=color228](3.26,7.8588476)(4.06,8.658848)
\psline[linewidth=0.04cm,linecolor=color228](4.06,8.658848)(4.86,7.8588476)
\psline[linewidth=0.04cm,linecolor=color228](4.86,7.8588476)(5.66,8.658848)
\psline[linewidth=0.04cm,linecolor=color228](5.66,8.658848)(6.46,7.8588476)
\psline[linewidth=0.04cm,linecolor=color228](8.66,7.8588476)(9.46,8.658848)
\psline[linewidth=0.04cm,linecolor=color228](9.46,8.658848)(9.86,8.258847)
\psline[linewidth=0.04cm,linecolor=color228](9.86,8.258847)(10.66,9.058847)
\psline[linewidth=0.04cm,linecolor=color228](10.66,9.058847)(11.46,8.258847)
\psline[linewidth=0.04cm,linecolor=color228](11.46,8.258847)(12.26,9.058847)
\psline[linewidth=0.04cm,linecolor=color228](12.26,9.058847)(13.46,7.8588476)
\psline[linewidth=0.04cm,linecolor=color228](13.46,7.8588476)(14.26,8.658848)
\psline[linewidth=0.04cm,linecolor=color228](14.26,8.658848)(15.06,7.8588476)
\psline[linewidth=0.04cm,linecolor=color228](22.46,4.4588475)(23.26,5.2588477)
\psline[linewidth=0.04cm,linecolor=color228](23.26,5.2588477)(23.66,4.8588476)
\psline[linewidth=0.04cm,linecolor=color228](23.66,4.8588476)(24.46,5.658848)
\psline[linewidth=0.04cm,linecolor=color228](24.46,5.658848)(24.86,5.2588477)
\psline[linewidth=0.04cm,linecolor=color228](24.86,5.2588477)(25.66,6.0588474)
\psline[linewidth=0.04cm,linecolor=color228](25.66,6.0588474)(27.26,4.4588475)
\psline[linewidth=0.04cm,linecolor=color228](27.26,4.4588475)(28.06,5.2588477)
\psline[linewidth=0.04cm,linecolor=color228](28.06,5.2588477)(28.86,4.4588475)
\psline[linewidth=0.04cm,linecolor=color228](22.46,0.058847655)(23.26,0.8588477)
\psline[linewidth=0.04cm,linecolor=color228](23.26,0.8588477)(23.66,0.45884764)
\psline[linewidth=0.04cm,linecolor=color228](23.66,0.45884764)(25.26,2.0588477)
\psline[linewidth=0.04cm,linecolor=color228](25.26,2.0588477)(27.26,0.058847655)
\psline[linewidth=0.04cm,linecolor=color228](27.26,0.058847655)(28.06,0.8588477)
\psline[linewidth=0.04cm,linecolor=color228](28.06,0.8588477)(28.86,0.058847655)
\psline[linewidth=0.04cm,linecolor=color258](11.46,4.4588475)(12.26,5.2588477)
\psline[linewidth=0.04cm,linecolor=color258](12.26,5.2588477)(12.66,4.8588476)
\psline[linewidth=0.04cm,linecolor=color258](12.66,4.8588476)(13.46,5.658848)
\psline[linewidth=0.04cm,linecolor=color258](13.46,5.658848)(14.26,4.8588476)
\psline[linewidth=0.04cm,linecolor=color258](14.26,4.8588476)(15.06,5.658848)
\psline[linewidth=0.04cm,linecolor=color258](15.06,5.658848)(15.86,4.8588476)
\psline[linewidth=0.04cm,linecolor=color258](15.86,4.8588476)(16.66,5.658848)
\psline[linewidth=0.04cm,linecolor=color258](16.66,5.658848)(17.86,4.4588475)
\psline[linewidth=0.04cm,linecolor=color258](11.46,0.25884765)(12.26,1.0588477)
\psline[linewidth=0.04cm,linecolor=color258](12.26,1.0588477)(12.66,0.65884763)
\psline[linewidth=0.04cm,linecolor=color258](12.66,0.65884763)(13.46,1.4588476)
\psline[linewidth=0.04cm,linecolor=color258](13.46,1.4588476)(13.86,1.0588477)
\psline[linewidth=0.04cm,linecolor=color258](13.86,1.0588477)(14.66,1.8588476)
\psline[linewidth=0.04cm,linecolor=color258](14.66,1.8588476)(15.86,0.65884763)
\psline[linewidth=0.04cm,linecolor=color258](15.86,0.65884763)(16.66,1.4588476)
\psline[linewidth=0.04cm,linecolor=color258](16.66,1.4588476)(17.86,0.25884765)
\psline[linewidth=0.04cm,linecolor=color258](11.46,-3.9411523)(12.26,-3.1411524)
\psline[linewidth=0.04cm,linecolor=color258](12.26,-3.1411524)(12.66,-3.5411522)
\psline[linewidth=0.04cm,linecolor=color258](12.66,-3.5411522)(13.46,-2.7411523)
\psline[linewidth=0.04cm,linecolor=color258](13.46,-2.7411523)(13.86,-3.1411524)
\psline[linewidth=0.04cm,linecolor=color258](13.86,-3.1411524)(14.66,-2.3411524)
\psline[linewidth=0.04cm,linecolor=color258](14.66,-2.3411524)(15.46,-3.1411524)
\psline[linewidth=0.04cm,linecolor=color258](15.46,-3.1411524)(16.26,-2.3411524)
\psline[linewidth=0.04cm,linecolor=color258](16.26,-2.3411524)(17.86,-3.9411523)
\psline[linewidth=0.04cm,linecolor=color258](22.66,-3.9411523)(23.46,-3.1411524)
\psline[linewidth=0.04cm,linecolor=color258](23.46,-3.1411524)(23.86,-3.5411522)
\psline[linewidth=0.04cm,linecolor=color258](23.86,-3.5411522)(25.46,-1.9411523)
\psline[linewidth=0.04cm,linecolor=color258](25.46,-1.9411523)(27.06,-3.5411522)
\psline[linewidth=0.04cm,linecolor=color258](27.06,-3.5411522)(27.86,-2.7411523)
\psline[linewidth=0.04cm,linecolor=color258](27.86,-2.7411523)(29.06,-3.9411523)
\psline[linewidth=0.04cm,linecolor=color258](16.86,-8.541152)(17.66,-7.7411523)
\psline[linewidth=0.04cm,linecolor=color258](17.66,-7.7411523)(18.06,-8.141152)
\psline[linewidth=0.04cm,linecolor=color258](18.06,-8.141152)(19.66,-6.5411525)
\psline[linewidth=0.04cm,linecolor=color258](19.66,-6.5411525)(20.86,-7.7411523)
\psline[linewidth=0.04cm,linecolor=color258](20.86,-7.7411523)(21.66,-6.9411526)
\psline[linewidth=0.04cm,linecolor=color258](21.66,-6.9411526)(23.26,-8.541152)
\psline[linewidth=0.04cm,linecolor=color294](26.66,-8.541152)(27.46,-7.7411523)
\psline[linewidth=0.04cm,linecolor=color294](27.46,-7.7411523)(27.86,-8.141152)
\psline[linewidth=0.04cm,linecolor=color294](27.86,-8.141152)(29.46,-6.5411525)
\psline[linewidth=0.04cm,linecolor=color294](29.46,-6.5411525)(30.26,-7.341152)
\psline[linewidth=0.04cm,linecolor=color294](30.26,-7.341152)(31.06,-6.5411525)
\psline[linewidth=0.04cm,linecolor=color294](31.06,-6.5411525)(33.06,-8.541152)
\usefont{T1}{ptm}{m}{n}
\rput(3.1214552,7.5038476){$P/_0Q$}
\usefont{T1}{ptm}{m}{n}
\rput(25.541454,-0.23615235){$P\backslash_0Q$}
\usefont{T1}{ptm}{m}{n}
\rput(14.361455,4.1638474){$P/_1Q$}
\usefont{T1}{ptm}{m}{n}
\rput(19.921455,-8.876152){$P\backslash_1Q$}
\usefont{T1}{ptm}{m}{n}
\rput(29.831455,-8.936152){$P/_2Q = P\backslash _2Q$}
\end{pspicture} 
}

\end{example}

\bigskip

\bigskip

\section{$\Dy^m$ algebras}
\medskip

We apply Theorem \ref{theorem1} to introduce the notion of $\Dy ^m$ algebra, for $m\geq 1$. When $m=1$, we recover J.-L. Loday\rq s dendriform algebras.
\medskip

The present section contains two main results: \begin{enumerate}
\item We prove that the vector space generated by all $m$-Dyck paths, with the products $*_i$, $0\leq i\leq m$, is the free $\Dy ^m$ algebra on one generator.
\item We define, for $m\geq 2$ and any composition ${\underline r}$ of $m+1$ in $l+1$ parts, a functor ${\mathbb F}_{\underline r}$ from the category of $\Dy ^m$ algebras into the category of $\Dy ^l$ algebras,  which sends free objects into free objects.
\end{enumerate}
\medskip

\begin{definition} \label{defDyckmalgebra} A {\it $\Dy ^m$ algebra} over $\K$ is a vector space $D$ equipped with $m+1$ binary operations $*_i: D\ot D\longrightarrow D$, for $0\leq i\leq m$, satisfying the following relations:\begin{enumerate}
\item $x*_i(y*_j z) = (x*_i y)*_j z$, for $0\leq i < j\leq m$,
\item ${\displaystyle \sum _{j=0}^i x*_i (y*_j z) = \sum _{k=i}^m (x*_k y)*_i z}$,\end{enumerate}
for any elements $x, y $ and $z$ in $D$.\end{definition}
\medskip

Clearly, a $\Dy^1$ algebra is a dendriform algebra, as described in Definition \ref{defdendriform}. 

\begin{remark} \label{rempropDyckalg} Let $D$ be a $\Dy ^m$ algebra. The relations of Definition \ref{defDyckmalgebra} imply that,\begin{enumerate}
\item the underlying vector space $D$, with the product $* := {\displaystyle\sum _{i=0}^m *_i}$, is an associative algebra.
\item for $1\leq l\leq m$ and any composition ${\underline r}=(r_0,\dots ,r_l)$ of $m+1$ of length $l+1$, the vector space $D$ equipped with 
 the binary operations
$$x {\overline {*}}_i y := \sum _{j = r_0+\dots +r_{i-1}+1}^{r_o+\dots +r_i} x*_j y,$$
for $0\leq i\leq l$, where $r_{-1} = -1$, is a $\Dy ^l$ algebra. 
So, they define a functor ${\mathbb F}_{\underline r}$ from the category of $\Dy ^m$ algebras into the category of $\Dy ^l$ algebras.\end{enumerate}\end{remark}
\medskip

Note that, as particular cases of Remark \ref{rempropDyckalg} we get that for any $\Dy ^m$ algebra $D$ and any $0\leq k\leq m-1$,
the vector space $D$ equipped with the binary operations ${\displaystyle \succ ^k := \sum _{i=0}^k *_i}$ and 
${\displaystyle \prec ^k := \sum _{i = k+1}^m *_i}$,
is a dendriform algebra.
\medskip

The following result is immediate to verify.

\begin{lemma} \label{lemfinctorsFr} For integers $0\leq k< h <m$, let ${\underline r} = (r_0,\dots ,r_h)$ be a composition of $m+1$ and 
${\underline s}=(s_0,\dots ,s_k)$ be a 
composition of $h +1$. Let ${\underline s}\circ {\underline r}$ be the composition $(r_0+\dots +r_{s_0}, r_{s_0+1}+\dots +r_{s_0+s_1},\dots ,
r_{s_0+\dots +s_{k-1}+1}+\dots +r_h)$, we have that 
$${\mathbb F}_{\underline s}\circ {\mathbb F}_{\underline r} = {\mathbb F}_{{\underline s}\circ {\underline r}}.$$\end{lemma}
\bigskip

\begin{notation} Theorem \ref{theorem1} asserts that the graded vector space $\K [\Dy ^m]$ spanned by the set of all $m$-Dyck paths, equipped with the operations $*_i$ defined in Section $3$, is a $\Dy^m$ algebra, for all $m\geq 1$. From now on we denote this $\Dy^m$ algebra by $\D _m$. \end{notation}

As the relations of Definition \ref{defDyckmalgebra} keep the order of the variables, the algebraic operad (see \cite{LodVal}) of $\Dy ^m$ algebras is regular, which means that the operad is described completely by the free object on one generator.

We now turn to prove that $\D_m$ is in fact the free $\Dy^m$ algebra on one generator. Before doing it, let us describe a simple way to describe the free $\Dy ^m$ algebra.
\bigskip

\begin{remark}\label{rem:freeDyckwithbinarytrees} Let $A$ be a vector space, equipped with a family $*_0,\dots ,*_m$ of binary operations. Definition \ref{defDyckmalgebra} states that $(A, *_0,\dots ,*_m)$ is a $\Dy ^m$ algebra if, and only if, the operations $*_i$, $0\leq i\leq m$, satisfy the following relations:\begin{enumerate}
\item $(x*_i y)*_j z = x*_i (y*_j z)$, for $0\leq i< j\leq m$,
\item $(x*_i y)*_iz = \displaystyle \sum _{j = 0}^i x*_i (y*_j z)-\sum _{j = i+1}^m (x*_j y)*_i z$, for $0\leq i\leq m$, \end{enumerate}
for $x, y$ and $z$ in $A$.\end{remark}

\bigskip

For $n\geq 1$, let ${\mathcal Y}_n^m$ be the set of all planar binary rooted trees with $n+1$ leaves (and $n$ internal vertices), with the vertices colored by the elements of $\{ *_0,\dots ,*_m\}$. Given two colored trees, $t$ and $w$ and an integer $0\leq i\leq m$, we denote by $t\vee _{*_i} w$ the colored tree obtained by connecting the roots of $t$ and $w$ to a new root colored by $*_i$. 

For any internal vertex $v$ of a colored planar binary rooted tree $t\in {\mathcal Y}_n^m$, we denote by $t_v$ the colored subtree of $t$ whose root is $v$. 
\medskip

\begin{definition} \label{def:basisBm} For $n\geq 2$, define the set ${\mathcal B}_n^{m}$ as the subset of all the elements $t$ in  ${\mathcal Y}_{n-1}^m$ such that any subtree $t_v$ satisfies the condition:

$\bold {(C)}$ if $t_v = t_v^l \vee _{*_i} t_v^r$, then the color of the root of $t_v^l$ is $*_j$ for some $ j >i$.\end{definition}
\medskip

For instance, the tree  $t = $
\scalebox{0.8} 
{
\begin{pspicture}(0,-0.80668944)(2.5949707,0.78668946)
\psline[linewidth=0.04cm](0.042070314,0.74668944)(1.2420703,-0.45331055)
\psline[linewidth=0.04cm](1.2220703,-0.47331056)(2.4220703,0.72668946)
\psline[linewidth=0.04cm](1.4420704,0.7666895)(1.8420703,0.16668946)
\psline[linewidth=0.04cm](1.9220703,0.7066895)(1.6420703,0.46668947)
\psline[linewidth=0.04cm](1.2420703,-0.47331056)(1.2420703,-0.6733105)
\usefont{T1}{ptm}{m}{n}
\rput(0.6495557,-0.068310544){\small $*_2$}
\usefont{T1}{ptm}{m}{n}
\rput(0.40955567,0.21168946){\small $*_3$}
\usefont{T1}{ptm}{m}{n}
\rput(1.4495556,-0.58831054){\small $*_1$}
\usefont{T1}{ptm}{m}{n}
\rput(0.86955565,0.41168946){\small $*_1$}
\usefont{T1}{ptm}{m}{n}
\rput(2.1095557,0.07168946){\small $*_2$}
\usefont{T1}{ptm}{m}{n}
\rput(1.8895557,0.43168944){\small $*_0$}
\psline[linewidth=0.04cm](0.6820703,0.12668945)(1.2620703,0.72668946)
\psline[linewidth=0.04cm](0.8620703,0.7066895)(0.4420703,0.32668945)
\psline[linewidth=0.04cm](0.46207032,0.74668944)(0.7020703,0.48668945)
\end{pspicture} 
}
does not belong to ${\mathcal B}_6^3$, because in the subtree $t_v = $ \scalebox{0.6} 
{
\begin{pspicture}(0,-0.48668945)(1.3529004,0.46668947)
\psline[linewidth=0.04cm](0.0,0.44668946)(0.6,-0.15331055)
\psline[linewidth=0.04cm](0.62,-0.17331055)(1.22,0.42668945)
\psline[linewidth=0.04cm](0.62,0.40668947)(0.3,0.10668945)
\usefont{T1}{ptm}{m}{n}
\rput(0.56748533,0.051689453){\small $*_0$}
\psline[linewidth=0.04cm](0.62,-0.23331055)(0.64,-0.17331055)
\psline[linewidth=0.04cm](0.64,-0.23331055)(0.62,-0.43331054)
\usefont{T1}{ptm}{m}{n}
\rput(0.86748534,-0.26831055){\small $*_2$}
\end{pspicture} 
}
the root of $t_v^l$ is colored $*_0$, while the root of $t_v$ is colored with $*_2$. 
\medskip

For $n=1$, ${\mathcal B}_1^m$ is the set which has as unique element the tree with one leave and no vertex: $|$. 
Let ${\mathcal B}^m=\bigcup _{n\geq 1}{\mathcal B}_n^m$.

Note that for any $t = t^l\vee _{*_i}t^r\in {\mathcal B}^m$ the trees $t^l$ and $t^r$ belong to ${\mathcal B}^m$.
\medskip

For any set $X$, let ${\mathcal B}_n^m(X)$ denote 
the set of all trees in ${\mathcal B}_n^m$ with leaves colored by the elements of $X$. Let $\Dy^m(X)$ be the graded vector space whose basis is the set ${\displaystyle \bigcup _{n\geq 1}{\mathcal B}_n^m(X)}$. 

For any pair of trees $t\in {\mathcal B}_n^m(X)$ and $w\in {\mathcal B}_r^m(X)$, with $n,r\geq 1$, and any integer $0\leq i\leq m$, the product $t*_i w\in \Dy^m(X)$ is defined recursively on $n+r$ as follows,\begin{enumerate}
 \item for $ n = r = 1$, we have \scalebox{0.7} 
{
\begin{pspicture}(0,-0.52871096)(2.9949708,0.52871096)
\psline[linewidth=0.04cm](0.2220703,0.18466797)(0.2220703,-0.21533203)
\psline[linewidth=0.04cm](1.0220703,0.18466797)(1.0220703,-0.21533203)
\psline[linewidth=0.04cm](1.8220704,0.18466797)(2.2220702,-0.21533203)
\psline[linewidth=0.04cm](2.2220702,-0.21533203)(2.6220703,0.18466797)
\psline[linewidth=0.04cm](2.2220702,-0.21533203)(2.2220702,-0.41533202)
\usefont{T1}{ptm}{m}{n}
\rput(0.22955567,0.30966797){\small $x$}
\usefont{T1}{ptm}{m}{n}
\rput(0.98955566,0.30966797){\small $y$}
\usefont{T1}{ptm}{m}{n}
\rput(1.7895557,0.32966796){\small $x$}
\usefont{T1}{ptm}{m}{n}
\rput(2.6895556,0.34966797){\small $y$}
\usefont{T1}{ptm}{m}{n}
\rput(0.6335254,-0.010332031){$*_i$}
\usefont{T1}{ptm}{m}{n}
\rput(2.5095556,-0.31033203){\small $*_i$}
\usefont{T1}{ptm}{m}{n}
\rput(1.4635254,0.009667968){$:=$}
\end{pspicture} 
}

\item for $t = t^l\vee_{*_j} t^r$ or $n = 1$, with $i < j\leq m$, the product $*_i$ of $t$ and $w$ is $t*_i w := t\vee _{*_i} w \in {\mathcal B}_{n+r}^m(X)$,
\item for $t = t^l\vee_{*_j} t^r$, with $0\leq j\leq i$, we have that\begin{enumerate}
\item when $j < i$, the recursive hypothesis states that $t^r\vee _{*_i} w$ is defined, 
and we put $t *_i w := t^l\vee_{*_j} (t^r\vee_{*_i} w)$,
\item when $j = i$, by Remark \ref{rem:freeDyckwithbinarytrees}, we get
$$\qquad \qquad 
t*_iw := ( t^l *_i t^r)*_i w=
\sum _{k = 0}^i t^l *_k (t^r *_i w)- \sum _{k=i+1}^m ( t^l *_k t^r)*_i w$$
For the second sum, for any $i < k\leq m$, by a recursive argument we suppose that $t^l*_k t^r =\sum _{\alpha} {\tilde {t}}_{k\alpha}$. Moreover, Remark \ref{rem:freeDyckwithbinarytrees} implies that the root of any $\tilde {t}_{k\alpha}$ is colored by an $*_h$ with $i < h$. 

So, $( t^l *_k t^r)*_i w := \sum _{k\alpha} {\tilde {t}}_{k\alpha}\vee _{*_i} w$.

For the first sum, as $t $ belongs to ${\mathcal B}_n^m(X)$, we know that the root of $t^l$ is colored by a $*_h$ with $h > i$, and therefore $h > k$ for all $0\leq k\leq i$. On the other hand, the recursive hypothesis implies that $t^r *_i w = \sum _{\beta} {\tilde {w}}_{i\beta}$ is defined. 
So, $t^l *_k (t^r *_i w) := \sum _{\beta} t^l\vee _{*_k} {\tilde {w}}_{i\beta}$. 

Finally, the formula for $t*_iw = ( t^l *_i t^r)*_i w$ is
$$\qquad t*_iw :=  \sum _{k = 0}^i\bigl ( \sum _{\beta} t^l\vee _{*_k} {\tilde {w}}_{i\beta}\bigr )\ {-}\ \sum _{k=i+1}^m\bigl ( \sum _{k\alpha} {\tilde {t}}_{k\alpha}\vee _{*_i} w\bigr ),$$
where $t^l*_k t^r =\sum _{\alpha} {\tilde {t}}_{k\alpha}$, for $i+1\leq k\leq m$, and $t^r *_i w = \sum _{\beta} {\tilde {w}}_{i\beta}$.\end{enumerate}\end{enumerate}
\medskip

\begin{example} Let 
\medskip

\scalebox{1} 
{
\begin{pspicture}(0,-0.8287109)(6.193916,0.8287109)
\usefont{T1}{ptm}{m}{n}
\rput(0.2924707,0.009667968){$t=$}
\psline[linewidth=0.04cm](0.6810156,0.504668)(1.4810157,-0.49533203)
\psline[linewidth=0.04cm](1.4810157,-0.49533203)(2.2810156,0.504668)
\psline[linewidth=0.04cm](1.0810156,-0.015332031)(1.5010157,0.48466796)
\psline[linewidth=0.04cm](1.8810157,0.504668)(2.0810156,0.30466798)
\psline[linewidth=0.04cm](1.1210157,0.56466794)(0.8810156,0.30466798)
\psline[linewidth=0.04cm](1.4810157,-0.49533203)(1.4810157,-0.69533205)
\usefont{T1}{ptm}{m}{n}
\rput(1.748501,-0.610332){\small $*_1$}
\usefont{T1}{ptm}{m}{n}
\rput(1.00,-0.050332032){\small $*_2$}
\usefont{T1}{ptm}{m}{n}
\rput(0.78,0.20966797){\small $*_3$}
\usefont{T1}{ptm}{m}{n}
\rput(2.248501,0.14966796){\small $*_0$}
\usefont{T1}{ptm}{m}{n}
\rput(3.2950392,0.04966797){and}
\usefont{T1}{ptm}{m}{n}
\rput(4.182471,0.029667968){$w =$}
\psline[linewidth=0.04cm](4.4810157,0.504668)(5.0810156,-0.095332034)
\psline[linewidth=0.04cm](5.0810156,-0.095332034)(5.6810155,0.504668)
\psline[linewidth=0.04cm](5.0810156,-0.095332034)(5.0810156,-0.29533204)
\psline[linewidth=0.04cm](5.0610156,0.524668)(4.7610154,0.22466797)
\usefont{T1}{ptm}{m}{n}
\rput(4.59,0.20){\small $*_2$}
\usefont{T1}{ptm}{m}{n}
\rput(5.348501,-0.23033203){\small $*_1$}
\usefont{T1}{ptm}{m}{n}
\rput(0.648501,0.589668){\small $x_1$}
\usefont{T1}{ptm}{m}{n}
\rput(1.108501,0.589668){\small $x_2$}
\usefont{T1}{ptm}{m}{n}
\rput(1.588501,0.60966796){\small $x_3$}
\usefont{T1}{ptm}{m}{n}
\rput(1.948501,0.629668){\small $x_4$}
\usefont{T1}{ptm}{m}{n}
\rput(2.4485009,0.649668){\small $x_5$}
\usefont{T1}{ptm}{m}{n}
\rput(4.468501,0.569668){\small $y_1$}
\usefont{T1}{ptm}{m}{n}
\rput(5.088501,0.589668){\small $y_2$}
\usefont{T1}{ptm}{m}{n}
\rput(5.708501,0.589668){\small $y_3$}
\end{pspicture} 
}
\medskip

we get that
\medskip

\scalebox{1} 
{
\begin{pspicture}(0,-0.9015332)(9.544922,0.9015332)
\usefont{T1}{ptm}{m}{n}
\rput(0.7924707,-0.1231543){$t*_2 w =$}
\psline[linewidth=0.04cm](1.5410156,0.2318457)(2.1410155,-0.3681543)
\psline[linewidth=0.04cm](2.1410155,-0.3681543)(2.7410157,0.2318457)
\psline[linewidth=0.04cm](2.1410155,-0.3681543)(2.1410155,-0.5681543)
\psline[linewidth=0.04cm](2.1210155,0.2318457)(1.8210156,-0.0681543)
\usefont{T1}{ptm}{m}{n}
\rput(1.5245312,0.3318457){\footnotesize $x_1$}
\usefont{T1}{ptm}{m}{n}
\rput(2.1245313,0.3318457){\footnotesize $x_2$}
\usefont{T1}{ptm}{m}{n}
\rput(2.7645311,0.3318457){\footnotesize $x_3$}
\usefont{T1}{ptm}{m}{n}
\rput(1.75,-0.103154294){\small $*_3$}
\usefont{T1}{ptm}{m}{n}
\rput(2.388501,-0.4431543){\small $*_2$}
\usefont{T1}{ptm}{m}{n}
\rput(3.4324708,-0.063154295){$*_1$}
\psline[linewidth=0.04cm](3.8410156,0.2918457)(4.6410155,-0.5081543)
\psline[linewidth=0.04cm](4.6410155,-0.5081543)(5.4410157,0.2918457)
\psline[linewidth=0.04cm](4.6410155,-0.5081543)(4.6410155,-0.7081543)
\psline[linewidth=0.04cm](4.5810156,0.3118457)(5.0010157,-0.1881543)
\psline[linewidth=0.04cm](5.0210156,0.2918457)(4.801016,0.051845703)
\usefont{T1}{ptm}{m}{n}
\rput(4.888501,-0.6031543){\small $*_0$}
\usefont{T1}{ptm}{m}{n}
\rput(5.248501,-0.22315429){\small $*_1$}
\usefont{T1}{ptm}{m}{n}
\rput(5.048501,0.016845703){\small $*_2$}
\psline[linewidth=0.04cm](4.341016,0.2918457)(4.841016,-0.2881543)
\usefont{T1}{ptm}{m}{n}
\rput(5.068501,-0.40315428){\small $*_2$}
\usefont{T1}{ptm}{m}{n}
\rput(3.7845314,0.3918457){\footnotesize $x_4$}
\usefont{T1}{ptm}{m}{n}
\rput(4.224531,0.4118457){\footnotesize $x_5$}
\usefont{T1}{ptm}{m}{n}
\rput(4.724531,0.3918457){\footnotesize $y_1$}
\usefont{T1}{ptm}{m}{n}
\rput(5.1045313,0.3918457){\footnotesize $y_2$}
\usefont{T1}{ptm}{m}{n}
\rput(5.6045313,0.3918457){\footnotesize $y_3$}
\usefont{T1}{ptm}{m}{n}
\rput(6.3624706,-0.0831543){$=$}
\psline[linewidth=0.04cm](6.5810156,0.5918457)(7.7810154,-0.6081543)
\psline[linewidth=0.04cm](7.7810154,-0.6081543)(8.981015,0.5918457)
\psline[linewidth=0.04cm](7.7810154,-0.6081543)(7.7810154,-0.8081543)
\psline[linewidth=0.04cm](7.9610157,0.5718457)(8.421016,-0.008154297)
\usefont{T1}{ptm}{m}{n}
\rput(6.65501,0.3968457){\small $*_3$}
\usefont{T1}{ptm}{m}{n}
\rput(6.8501,0.1568457){\small $*_2$}
\usefont{T1}{ptm}{m}{n}
\rput(8.028501,-0.6831543){\small $*_1$}
\usefont{T1}{ptm}{m}{n}
\rput(8.428501,-0.2831543){\small $*_0$}
\usefont{T1}{ptm}{m}{n}
\rput(8.628501,-0.0831543){\small $*_2$}
\usefont{T1}{ptm}{m}{n}
\rput(8.708501,0.37684572){\small $*_2$}
\usefont{T1}{ptm}{m}{n}
\rput(8.828501,0.0968457){\small $*_1$}
\usefont{T1}{ptm}{m}{n}
\rput(6.5045314,0.7118457){\footnotesize $x_1$}
\usefont{T1}{ptm}{m}{n}
\rput(6.9245315,0.7118457){\footnotesize $x_2$}
\usefont{T1}{ptm}{m}{n}
\rput(7.244531,0.7318457){\footnotesize $x_3$}
\usefont{T1}{ptm}{m}{n}
\rput(7.5845313,0.7118457){\footnotesize $x_4$}
\usefont{T1}{ptm}{m}{n}
\rput(7.9445314,0.6918457){\footnotesize $x_5$}
\usefont{T1}{ptm}{m}{n}
\rput(8.284532,0.6918457){\footnotesize $y_1$}
\usefont{T1}{ptm}{m}{n}
\rput(8.664532,0.6918457){\footnotesize $y_2$}
\usefont{T1}{ptm}{m}{n}
\rput(9.104531,0.7118457){\footnotesize $y_3$}
\psline[linewidth=0.04cm](6.7210155,0.4318457)(6.881016,0.5918457)
\psline[linewidth=0.04cm](7.2610154,0.25184572)(7.2410154,0.3118457)
\psline[linewidth=0.04cm](6.9010158,0.2118457)(7.2410154,0.5918457)
\psline[linewidth=0.04cm](7.5810156,0.5918457)(8.201015,-0.1481543)
\psline[linewidth=0.04cm](8.261016,0.61184573)(8.641016,0.1918457)
\psline[linewidth=0.04cm](8.701015,0.6318457)(8.481015,0.4118457)
\end{pspicture} 
}

\end{example}

The result below follows immediately from Remark \ref{rem:freeDyckwithbinarytrees} and the construction above.

\begin{proposition}\label{prop:freeDyckmontrees} For any set $X$, the graded vector space generated by the graded set ${\displaystyle \bigcup _{n\geq 1}{\mathcal B}_n^m(X)}$ equipped with the binary products define above is the free $\Dy ^m$ algebra on $X$.\end{proposition}
\medskip

\begin{notation} We denote by $\Dy ^m(X)$ the free $\Dy ^m$ algebra generated by a set $X$.\end{notation}
\bigskip

In order to prove that the $\Dy^m$ algebra ${\mathcal D}_m$ is the free $\Dy ^m$ algebra on one element, we need the following Proposition.

\begin{proposition}\label{Dyckpathsobtainedfromsmallerones} Any element of $P\in \D_m$ is a linear combination of elements of the form $R_1*_i R_2$, where $0\leq i\leq m$ and the sizes of $R _1$ and $R_2$ are strictly smaller than the size of $P$.
\end{proposition}
\medskip

\begin{proof} Let us point out that for $m = 1$, the result has been proved in \cite{LodRon}. 
\medskip

For the general case, let $$P = \bigvee _u (P_0,\dots , P_m) = P_0\t_0 (((\rho _m\t _m P_1)\t _{m-1}\dots )\t _1 P_m)\in \Dy_m^n.$$

It is immediate to see that $P\rq := ((\rho _m\t _m P_1)\t _{m-1}\dots )\t _1 P_m$ is prime. So, if $P_0\neq \bullet$, then 
$$P = P_0\t _0 P\rq = P_0*_0 P\rq,$$
and we are done.

Now suppose that $P_0 = \bullet$. 
The maximal element $P_{max(n)} = (0,\dots ,0,nm)$ of the Tamari lattice $\Dy _n^m$ satisfies 
that $P_{max(n)} =\rho_m*_m P_{max(n-1)}$. 

We may assume that the result is also true for elements $Q$ of size $n$ such that $P < Q \leq P_{max(n)}$ in the $m$-Tamari lattice. 

For $P= \bigvee _u (\bullet, P_1,\dots,P_m)$, 
let $0\leq i\leq m$ be the largest integer such that $P_i\neq \bullet$. 

Let $P\rq := ((\rho _m\t _m P_1)\t _{m-1}\dots )\t _{m{-}i+2} P_{i-1}$ be the Dyck path obtained from $P$ by collapsing $P_i$ to a point.
We get that the sizes of both $P\rq $ and $P_i$ are smaller than $P$\rq s size, and that the last $m_i+1$ letters of the word $\omega _{n - n_i}(P\rq)$ are equal to $1$. 

So, we get that $P\rq *_{m-i+1}P_i = P+\sum_k Q_k$, with $P< Q_k$ for all $k$. 
As we have supposed that all Dyck paths $Q$ such that $P < Q < P_{max(n)}$ are linear combinations of elements of type $R_1*_i R_2$, where $0\leq i\leq m$ and the sizes of $R _1$ and $R_2$ are strictly smaller than the size of $Q$, the result also holds for $P$.
\end{proof}
\medskip

The following theorem states that the graded vector space $\D_m$ also describes the algebraic operad $\Dy ^m$.

\begin{theorem}\label{theorem: Dm is free on one generator}
The free $\Dy^m$ algebra on one generator is isomorphic to $(\D_m,*_0,\dots,*_m)$.
\end{theorem}
\medskip

\begin{proof}
Let $\Dy^m(a)$ be the free $\Dy^m$ algebra on one generator $a$. 
As $\D_m$ is a $\Dy ^m$ algebra, there exists a unique homomorphism $\phi : \Dy^m(a)\longrightarrow \D_m$ such that $\phi (a)$ is $\rho _m$, the unique $m$-Dyck path of size $1$.
Proposition \ref{Dyckpathsobtainedfromsmallerones} implies that $\phi$ is surjective. 

The subspace of homogeneous elements of degree $n$ of $\D_m$ is generated by the subset $\Dy_n^m$ of $m$-Dyck paths of size $n$. Let $\Dy^m(a)_n$ be the subspace of elements of degree $n$ of $\Dy^m(a)$.

As $\phi $ is surjective, to prove that $\phi$ is an isomorphism it suffices to show that the dimension of the vector space 
$\Dy ^m(a)_n$ is the number of elements of the set $\Dy _n^m$, that is 
$$ {\mbox {dim}_{\K} (\Dy ^m(a)_n)}= \vert \Dy_n^m\vert = d_{m,n}.$$
\medskip

From Proposition \ref{prop:freeDyckmontrees}, we know the underlying vector space of $\Dy^m(a)_n$ is 
generated by the set ${\mathcal B}_n^m(a)$ of planar binary rooted trees with $n$ leaves colored by $a$ 
and the $(n-1)$ vertices colored by the elements of $\{ *_0,\dots ,*_m\}$ satisfying condition 5.5. {\bf (1)}. 

So, the dimension of $\Dy ^m(a)_n$ over $\K$ is the number of elements of the set ${\mathcal B}_n^m$, 
which we denote by $b_{m,n}$, for $n\geq 1$. 
\medskip

 The generating series of the set $\{b_{m,n}\}_{n\geq 1}$ is
\begin{align} f_m(x)=\sum_{n\geq 1}b_{m,n}x^{n}.\end{align}

We need only to prove that $b_{m,n} = d_{m,n}$, the number of $m$-Dyck paths of size $n$, for $n\geq 1$.

From Remark \ref{remark: generating series for Dyck paths}, the generating series 
$d_m(x) = \sum _{n\geq 1} d_{m,n}x^n$ of the family of integers $\{ d_{m,n}\}_{n\geq 1}$ satisfies 
\begin{align}
x\cdot (1+d_m(x))^{m+1} = d_m(x).
\end{align}
Therefore, to end the proof, it suffices to show that the generating series of  $\{b_{m,n}\} _{n\geq 1}$ 
satisfies the same recursion formula. 

\noindent Note that $b_{m,1}=1=d_{m,1}$.
\medskip

For any colored tree $t \in {\mathcal B}_n^m$, there exists a unique integer $r$, a unique collection of colored trees $w^1,\dots ,w^r$ in ${\mathcal B}^m$ and a word $*_{i_1}\dots  *_{i_r}$ in the alphabet $\{*_0,\dots ,*_m\}$ such that $i_1 >\dots > i_r$ and
$$t = (((\vert \vee _{*_{i_1}} w^1)\vee _{*_{i_2}} w^2)\dots )\vee_{*_{i_r}} w^r,$$ 
which implies that 
\begin{align}
f_m(x) = x \cdot (1+f_m(x))^{m+1} 
\end{align}
\end{proof}
\medskip

\begin{corollary}\label{corfreeoverV} Let $V$ be a $\K$-vector space. The free $\Dy^m$ algebra on $V$ is the vector space
$$\Dy^m(V): = \displaystyle\bigoplus_{n\geq 1}\D_{m,n}\otimes V^{\otimes n},$$
equipped with the binary products given by:
$$ P\ot (v_1\ot \dots \ot v_{n_1})\ *_i\ Q\ot (w_1\ot \dots \ot w_{n_2}) := (P*_i Q)\ot (v_1\ot \dots \ot v_{n_1}\ot w_1\ot \dots \ot w_{n_2}),$$
for any integer $0\leq i\leq m$, any Dyck paths $P\in \Dy _{n_1}^m$ and $Q\in \Dy _{n_2}^m$, and elements $v_1,\dots ,v_{n_1},w_1,\dots , w_{n_2}\in V$.\end{corollary}

\bigskip

\bigskip

In Remark \ref{rempropDyckalg} , we showed that, for $0\leq h \leq m$ and any composition ${\underline r}=(r_0,\dots ,r_h)$ of $m$, there exists a functor 
${\mathbb F}_{\underline r}$ from the category of $\Dy ^m$-algebras into the category of $\Dy ^h$ algebras (which is equivalent to an operad homomorphism from $\Dy ^h$ to $\Dy ^m$). 
\medskip

We want to show that the image under ${\mathbb F}_{\underline r}$ of a free $\Dy ^m$ algebra is free as a $\Dy ^h$ algebra, too. From Corollary \ref{corfreeoverV}, we get that it suffices to prove that the image ${\mathbb F}_{\underline r}(\Dy ^m(a))$ of the free $\Dy ^m$ algebra over one element, is free as a $\Dy^h$ algebra.

In order to do that, we need to introduce new basis of the underlying vector space of $\Dy (a)$, by modifying the basis ${\mathcal B}^m$ described at Definition \ref{def:basisBm}.
\medskip

\begin{notation} \label{not:combs} Given a family of colored trees $t_1,\dots ,t_p$ and a family of integers $0\leq i_1,\dots , i_p\leq m$, we denote by \begin{enumerate}
\item $\Omega_{i_1,\dots ,i_p}^L(t_1,\dots ,t_p)$ the colored tree 
$$\Omega_{i_1,\dots ,i_p}^L(t_1,\dots ,t_p) := (((\vert \vee_{*_{i_p}} t_p)\vee_{*_{i_{p-1}}} t_{p-1})\dots )\vee_{*_{i_1}} t_1,$$
\item $\Omega_{i_1,\dots ,i_p}^R(t_1,\dots ,t_p)$ the colored tree 
$$\Omega_{i_1,\dots ,i_p}^R(t_1,\dots ,t_p) := t_1\vee _{i_1} (t_2\vee _{*_{i_2}}(\dots (t_{p-1}\vee _{*_{i_{p-1}}} (t_p\vee _{*_{i_p}} \vert )))).$$

That is
\medskip

\scalebox{0.7} 
{
\begin{pspicture}(0,-1.0820215)(17.721016,1.1020215)
\usefont{T1}{ptm}{m}{n}
\rput(3.3124707,-0.057021484){$\Omega_{i_1,\dots ,i_p}^L(t_1,\dots ,t_p)=$}
\psline[linewidth=0.04cm](5.0810156,0.7379785)(6.4810157,-0.8620215)
\psline[linewidth=0.04cm](6.5010157,-0.8620215)(7.9010158,0.7379785)
\psline[linewidth=0.04cm](6.4810157,-0.8620215)(6.4810157,-1.0620215)
\psline[linewidth=0.04cm](5.4810157,0.7579785)(5.2610154,0.47797853)
\psline[linewidth=0.04cm](5.801016,0.4979785)(5.5610156,0.21797852)
\psline[linewidth=0.04cm](6.5410156,-0.3420215)(6.2810154,-0.6220215)
\usefont{T1}{ptm}{m}{n}
\rput(7.868501,0.8629785){\small $t_1$}
\usefont{T1}{ptm}{m}{n}
\rput(6.038501,0.5829785){\small $t_{p-1}$}
\usefont{T1}{ptm}{m}{n}
\rput(5.488501,0.8829785){\small $t_p$}
\usefont{T1}{ptm}{m}{n}
\rput(6.528501,-0.15702148){\small $t_2$}
\usefont{T1}{ptm}{m}{n}
\rput(7.6124706,-0.15702148){$,$}
\usefont{T1}{ptm}{m}{n}
\rput(11.082471,-0.17702149){$\Omega _{i_1,\dots ,i_p}^R(t_1,\dots ,t_p) =$}
\psline[linewidth=0.04cm](13.141016,0.71797854)(14.521015,-0.90202147)
\psline[linewidth=0.04cm](17.681015,-0.8620215)(17.701015,-0.8820215)
\psline[linewidth=0.04cm](14.521015,-0.8620215)(15.981015,0.7579785)
\psline[linewidth=0.04cm](14.521015,-0.84202147)(14.521015,-1.0220215)
\psline[linewidth=0.04cm](15.501016,0.7579785)(15.741015,0.5179785)
\psline[linewidth=0.04cm](14.501016,-0.3420215)(14.781015,-0.6420215)
\usefont{T1}{ptm}{m}{n}
\rput(13.128501,0.8829785){\small $t_1$}
\usefont{T1}{ptm}{m}{n}
\rput(14.528501,-0.21702148){\small $t_2$}
\usefont{T1}{ptm}{m}{n}
\rput(15.238501,0.6029785){\small $t_{p-1}$}
\usefont{T1}{ptm}{m}{n}
\rput(15.528501,0.9229785){\small $t_p$}
\psdots[dotsize=0.08](6.0810156,0.11797851)
\psline[linewidth=0.04cm](15.121016,0.4979785)(15.461016,0.13797851)
\psdots[dotsize=0.08](5.9210157,0.11797851)
\psdots[dotsize=0.08](15.121016,0.077978514)
\psdots[dotsize=0.08](14.941015,0.077978514)
\psdots[dotsize=0.08](15.261016,0.077978514)
\psdots[dotsize=0.08](5.7410154,0.13797851)
\end{pspicture} 
}

\end{enumerate}\end{notation}

 Note first that for any tree $t\in {\mathcal Y}_{n-1}^m$ there exist unique non negative integers $p$ and $q$, such that:
$$t = \Omega_{i_1,\dots ,i_p}^L(t_1,\dots ,t_p) =\Omega_{j_1,\dots ,j_q}^R(w_1,\dots ,w_q) ,$$
for a unique families of colored trees $t_1,\dots, t_p$ and $w_1,\dots ,w_q$ and unique collections of integers $i_1,\dots ,i_p$ and $j_1,\dots ,j_q$ in $\{ 0,\dots ,m\}$ with $i_1=j_1$.
In particular, $t = t^l\vee _{*_{k_0}}t^r$, for 
$$t^l = \Omega_{i_2,\dots ,i_p}^L(t_2,\dots ,t_p) = w_1,\qquad {\rm and}\qquad t^r= t_1 = \Omega_{j_2,\dots ,j_q}^R(w_2,\dots ,w_q),$$
and $k_0 = i_1 = j_1$.
\medskip

\begin{definition} \label{def:baseBmk} Given $0\leq k\leq m$, define ${\mathcal B}_n^{m,k}$ to be the set of planar binary rooted trees with $n$ leaves, with the vertices colored by the elements of $\{ *_0,\dots ,*_m\}$ such that for any vertex $v$, the tree $t_v$ satisfies the following conditions:\begin{enumerate}[(i)]
\item if $t_v = \Omega_{i_1,\dots ,i_p}^L(t_1,\dots ,t_p) $, with the root colored by $*_{i_1}$ for $i_1\neq k$, then either $i_2 = k$ or $i_2 > i_1$,
\item if $t_v = \Omega_{k,\dots ,i_p}^L(t_1,\dots ,t_p) =  \Omega_{k,\dots ,j_q}^R(w_1,\dots ,w_q)$,  then $i_s\in \{ k+1,\dots ,m\}$ for $s\geq 2$, and $j_h\in \{0,\dots ,k\}$ for $h\geq 2$.\end{enumerate}\end{definition}

The basis ${\mathcal B}^m$ coincides with the set ${\mathcal B}^{m,m}$, under this notation. 
\medskip

\begin{proposition}\label{prop:basisk} For any $0\leq k\leq m$, the set ${\mathcal B}^{m,k} = {\displaystyle \bigcup _{n\geq 1} {\mathcal B}_n^{m,k}}$ is a basis of the underlying vector space of the free $\Dy ^m$ algebra $\Dy^m(a)$. \end{proposition}
\medskip

\begin{proo} We know that ${\mathcal B}^m$ is a linear basis of the $\K$-vector space $\Dy ^m(a)$. 
We want to prove that there exists a bijective map $\varphi : {\mathcal B}^m\longrightarrow {\mathcal B}^{m,k}$ satisfying that:\begin{enumerate}[(i)]
\item $\varphi  (t) = t$, for all $t\in {\mathcal B}^m\cap {\mathcal B}^{m,k}$,
\item if $t = t^l\vee_{*_i} t^r$, with $i\neq k$, then $\varphi(t) = \varphi (t^l)\vee _{*_i} \varphi (t^r)$, where $t^l, t^r\in {\mathcal B}^m$,
\item if $t = t^l\vee_{*_k} t^r$, then the root of $\varphi (t)$, is $*_s$, for some $s\geq k$.
\item $t$ and $\varphi (t)$ represent the same element in $\Dy^m(a)$. \end{enumerate}
\medskip

For a colored tree $t\in {\mathcal B}^m\cap {\mathcal B}^{m,k}$, we define $\varphi (t) := t$. Clearly ${\mathcal B}_n^m = {\mathcal B}_n^{m,k}$, for $n = 1, 2$.

If $t\notin {\mathcal B}^{m,k}$, we apply a recursive argument on $\vert t\vert > 2$.

$(1)$ For $t = t^l\vee _{*_s}t^r$, with $t^l$ and $t^r$ in ${\mathcal B}^m$ and $s\neq k$, define $\varphi (t) := \varphi (t^l)\vee _{*_s}\varphi (t^r)$.

Note that, as $t\in {\mathcal B}^m$, we know that the root of $t^l$ is colored by $*_h$, for some $h > s$. The recursive hypothesis states that the colored planar rooted trees $\varphi (t^l)$ and $\varphi (t^r)$ belong to ${\mathcal B}^{m,k}$ and the color of the root of $\varphi (t^l)$ is $*_p$ for some $s < h\leq p$. So, $\varphi (t)\in {\mathcal B}^{m,k}$.
\bigskip

$(2)$ If $t = t^l\vee _{*_k} t^r$, with $t^l$ and $t^r$ in ${\mathcal B}^m$, then there exist unique pair of positive integers $p, q$ such that  
\begin{enumerate}[(i)]
\item $t^l = \Omega _{i_2,\dots ,i_p}^L(t_2,\dots ,t_p)$, for a unique family of trees $t_2,\dots ,t_p$ in ${\mathcal B}^m$ and unique nonnegative integers $ i_2 <\dots < i_p$, 
\item $t^r = \Omega _{j_2,\dots ,j_q}^R(w_2,\dots ,w_q)$, for a unique family of trees $w_2,\dots ,w_q$ in ${\mathcal B}^m$ and unique nonnegative integers $ j_2,\dots , j_q$.\end{enumerate}
\medskip

\begin{enumerate} [(a)]\item  If $j_h \leq k$, for all $2\leq h\leq q$, then a recursive argument on $q$ shows that
 $$\varphi (t^r) = \Omega _{j_2,\dots ,j_q}^R(\varphi (w_2),\dots ,\varphi (w_q)),$$
 where $w_h\in {\mathcal B}^m$, for $2\leq h\leq q$.
 In this case, we define $\varphi (t):= \varphi (t^l)\vee_{*_k}\varphi (t^r)$. 
\item If there exist at least one $2\leq h\leq q$ such that $j_h > k$, let $s$ be the minimal integer such that $j_s >k$, $2\leq s\leq q$. For  $2\leq h\leq s-1$, we have that $j_h \leq k < j_s$. 

Applying that $x*_i (y*_j z) = (x*_i y)*_j z$ in $\Dy ^m(a)$, whenever $0\leq i < j\leq m$, we get that
 the tree  $t$ describes the same element than the tree:
$$( t^l\vee_{*_k}\ (w_2\vee_{*_{j_2}}(\dots (w_{s{-}1}\vee_{ *_{j_{s{-}1}}} w_s))))\vee_{ *_{j_s}}\Omega _{j_{s+1},\dots ,j_q}^R(w_{s+1},\dots ,w_q).$$

That is, we replace the tree 
\medskip

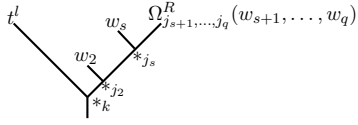
\scalebox{0.7} 
{
\begin{pspicture}(0,-1.0791992)(9.02291,1.0991992)
\usefont{T1}{ptm}{m}{n}
\rput(0.3224707,0.90580076){$t^l$}
\usefont{T1}{ptm}{m}{n}
\rput(2.2124708,-0.5141992){$*_{j_2}$}
\usefont{T1}{ptm}{m}{n}
\rput(2.8124707,0.06580078){$*_{j_s}$}
\usefont{T1}{ptm}{m}{n}
\rput(4.8624706,0.8858008){$\Omega_{j_{s+1},\dots ,j_q}^R(w_{s+1},\dots ,w_q)$}
\psline[linewidth=0.04cm](1.7210156,-0.65919924)(1.7210156,-1.0591992)
\usefont{T1}{ptm}{m}{n}
\rput(1.6924707,0.08580078){$w_2$}
\usefont{T1}{ptm}{m}{n}
\rput(2.2724707,0.72580075){$w_{s}$}
\usefont{T1}{ptm}{m}{n}
\rput(1.9924707,-0.7941992){$*_k$}
\psline[linewidth=0.04cm](1.7210156,-0.65919924)(3.1210155,0.7408008)
\psline[linewidth=0.04cm](1.7210156,-0.65919924)(0.32101563,0.7408008)
\psline[linewidth=0.04cm](2.3010156,0.60080075)(2.6210155,0.26080078)
\psline[linewidth=0.04cm](1.7210156,-0.059199218)(2.0210156,-0.35919923)
\end{pspicture} 
}
\medskip

by the tree
\medskip

\scalebox{0.7} 
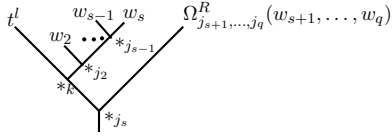
{
\begin{pspicture}(0,-1.2091992)(9.68291,1.2291992)
\psline[linewidth=0.04cm](0.34101564,0.8108008)(1.5410156,-0.38919923)
\usefont{T1}{ptm}{m}{n}
\rput(0.3224707,0.9958008){$t^l$}
\usefont{T1}{ptm}{m}{n}
\rput(1.8924707,-0.044199217){$*_{j_2}$}
\usefont{T1}{ptm}{m}{n}
\rput(2.2724707,-0.9441992){$*_{j_s}$}
\usefont{T1}{ptm}{m}{n}
\rput(5.5224705,1.0158008){$\Omega_{j_{s+1},\dots ,j_q}^R(w_{s+1},\dots ,w_q)$}
\psline[linewidth=0.04cm](1.5410156,-0.38919923)(1.9410156,-0.78919923)
\psline[linewidth=0.04cm](1.9410156,-0.78919923)(3.5410156,0.8108008)
\psline[linewidth=0.04cm](1.3410156,-0.18919922)(2.4210157,0.8908008)
\psline[linewidth=0.04cm](1.2810156,0.45080078)(1.6410156,0.07080078)
\psline[linewidth=0.04cm](1.8810157,0.9308008)(2.1810157,0.59080076)
\psdots[dotsize=0.08](1.6610156,0.5708008)
\psdots[dotsize=0.08](1.8410156,0.5708008)
\psdots[dotsize=0.08](2.0010157,0.59080076)
\psline[linewidth=0.04cm](1.9410156,-0.78919923)(1.9410156,-1.1891992)
\usefont{T1}{ptm}{m}{n}
\rput(1.2124707,0.59580076){$w_2$}
\usefont{T1}{ptm}{m}{n}
\rput(2.6124706,1.0358008){$w_{s}$}
\usefont{T1}{ptm}{m}{n}
\rput(1.3124707,-0.3441992){$*_k$}
\usefont{T1}{ptm}{m}{n}
\rput(1.8424706,1.0358008){$w_{s-1}$}
\usefont{T1}{ptm}{m}{n}
\rput(2.6424706,0.45580077){$*_{j_{s-1}}$}
\end{pspicture} 
}
\medskip

without changing the element in $\Dy^m(a)$.
\medskip

We cannot assume that the root of $\varphi (t^l\vee_{ *_k}(w_2\vee_{*_{j_2}}(\dots (w_{s{-}1}\vee_{ *_{j_{s{-}1}}} w_s))))$ is colored by $*_r$ with $r > j_s$. So, we have to work a bit more to define $\varphi (t)$.
\medskip

Suppose that $w_s = \Omega _{h_1,\dots ,h_u}^L(w_{s,1},\dots ,w_{s,u})$. As $w_s\in {\mathcal B}^m$, we get that $j_s < h_1<\dots < h_u$. 

The tree  $t^l\vee_{ *_k}(w_2\vee_{*_{j_2}}(\dots (w_{s{-}1}\vee_{ *_{j_{s{-}1}}} w_s)))$ represents the same element than the tree: 
$${\tilde {t^l}}:= (((t^l\vee_{ *_k}\Omega _{j_2,\dots ,j_{s{-}1}}^R(w_2,\dots ,w_{s{-}1}))\vee_{*_{h_u}} w_{s,u})\dots )\vee_{*_{h_1}} w_{s,1}=$$ 
\scalebox{0.7} 
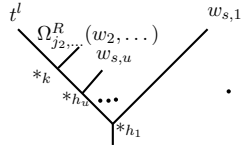
{
\begin{pspicture}(0,-1.3591992)(6.20291,1.3791993)
\psline[linewidth=0.04cm](1.4410156,0.8608008)(3.2410157,-0.9391992)
\psline[linewidth=0.04cm](3.2410157,-0.9391992)(5.0410156,0.8608008)
\usefont{T1}{ptm}{m}{n}
\rput(1.4024707,1.1858008){$t^l$}
\usefont{T1}{ptm}{m}{n}
\rput(1.8924707,-0.03419922){$*_k$}
\usefont{T1}{ptm}{m}{n}
\rput(3.5524707,-1.0541992){$*_{h_1}$}
\usefont{T1}{ptm}{m}{n}
\rput(2.5524707,-0.45419922){$*_{h_u}$}
\usefont{T1}{ptm}{m}{n}
\rput(5.372471,1.1458008){$w_{s,1}$}
\usefont{T1}{ptm}{m}{n}
\rput(3.2724707,0.32580078){$w_{s,u}$}
\psline[linewidth=0.04cm](2.2010157,0.12080078)(2.6010156,0.52080077)
\psline[linewidth=0.04cm](2.6610155,-0.35919923)(3.0410156,0.08080078)
\usefont{T1}{ptm}{m}{n}
\rput(2.97,0.74){$\Omega_{j_2,\dots }^R(w_2,\dots )$}
\psline[linewidth=0.04cm](3.2410157,-0.9391992)(3.2410157,-1.3391992)
\psdots[dotsize=0.08](3.1610155,-0.4991992)
\psdots[dotsize=0.08](3.0010157,-0.4991992)
\psdots[dotsize=0.08](3.3010156,-0.47919923)
\psdots[dotsize=0.08](5.4410157,-0.3191992)
\end{pspicture} 
}

Moreover, $t^l =\Omega _{i_2,\dots ,i_p}^L(t_2,\dots ,t_p)$ is such that $k < i_2<\dots < i_p$ and the root of $\varphi (t^l)$ is colored by $*_{i_2}$, which implies that 
$$\varphi ({\tilde {t^l}}) = (((\varphi (t^l)\vee _{*_k}\Omega _{j_2,\dots ,j_{s{-}1}}^R(\varphi(w_2),\dots ,\varphi (w_{s{-}1})))\vee _{*_{h_u}}\varphi (w_{s,u}))\dots )\vee _{*_{h_1}}\varphi (w_{s,1}),$$ and that the root of $\varphi ({\tilde{t^l}})$ is colored by $*_{h_1}$.
\medskip

Therefore the tree $t$ describes the same element than $$\tilde {t}:={\tilde {t^l}}\vee_{*_{j_s}}\Omega _{j_{s+1},\dots ,j_q}^R(w_{s+1},\dots ,w_q),$$
where the root of $\varphi ({\tilde {t^l}})$ is colored by $*_{h_1}$, with $h_1 > j_s > k$. We define
$$\varphi (t) := \varphi ({\tilde {t^l}})\vee _{*_{j_s}}\varphi (\Omega _{j_{s+1},\dots ,j_q}^R(w_{s+1},\dots ,w_q)).$$\end{enumerate}
\bigskip

To prove that $\varphi $ is bijective, we give an explicit description of $\varphi ^{-1}$.  Clearly, if $t\in {\mathcal B}^m\cap {\mathcal B}^{m,k}$, then $\varphi^{-1}(t) = t$.

For $t\in {\mathcal B}^{m,k}\setminus {\mathcal B}^m$, we use a recursive argument on the degree $\vert t\vert$. Suppose that for any $r < n$, the map $\varphi ^{-1}:{\mathcal B}_r^{m,k}\rightarrow {\mathcal B}_r^m$ is defined. From the conditions satisfied by $\varphi $, we get that its inverse satisfies that:\begin{enumerate}[(iv)]
\item if the root of $t$ is colored by $*_s$, for some $s \leq k$, then the root of $\varphi ^{-1}(t)$ is colored by $*_s$,
\item if the root of $t$ is colored by $*_s$, for some $s> k$, then the root of $\varphi ^{-1}(t)$ is colored by $*_s$ or by $*_k$.\end{enumerate}
\medskip

Let $t = t^l\vee _{*_{i_1}} t^r\in {\mathcal B}^{m,k}$. If $i_1\leq k$, then we know that the root of $\varphi^{-1}(t^l)$ is colored by $*_s$, for some $s  > i_1$. So, we get that:
$$\varphi ^{-1}(t) = \varphi ^{-1}(t^l)\vee _{*_{i_1}}\varphi ^{-1}(t^r).$$

Suppose that $i_1 > k$ and that $t= \Omega _{i_1,\dots ,i_p}^L(t^r,t_2,\dots ,t_p) \notin {\mathcal B}^m$. If $i_1< i_2 < \dots < i_p$, then it is immediate to see that
$$\varphi ^{-1}(t) = \varphi^{-1}(t^l)\vee_{*_{i_1}}\varphi^{-1}(t^r) = \varphi ^{-1}(\Omega _{i_1,\dots ,i_p}^L(\varphi ^{-1}(t^r), \varphi ^{-1}(t_2),\dots ,\varphi^{-1}(t_p))).$$
\medskip

Otherwise, there exists a unique integer $1\leq h \leq p$, such that $$k < i_1 <\dots < i_{h{-}1},\qquad {\rm and}\qquad i_h = k< i_{h+1}<\dots <i_p.$$ Moreover, as $t\in {\mathcal B}^{m,k}$, we have that $t_h = \Omega _{j_1,\dots ,j_q}^R(w_1,\dots ,w_q)\in {\mathcal B}^{m,k}$, with $j_1\leq k$.

From the definition of $\varphi$ and a recursive argument, we get that:\begin{enumerate}
\item the tree $t$ represents the same element than the tree
$${\tilde t} := \Omega_{i_{h+1},\dots ,i_p}^L(t_{h+1},\dots ,t_p)\vee _{*_k} (w_1\vee_{*_{j_1}}(\dots (w_q\vee _{*_{j_q}}\Omega _{i_1,\dots ,i_{h{-}1}}^L(t^r,t_2,\dots ,t_{h{-}1}))))=$$
\scalebox{0.7} 
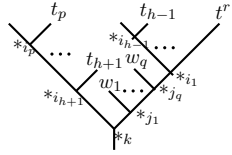
{
\begin{pspicture}(0,-1.4091992)(5.02291,1.4291992)
\psline[linewidth=0.04cm](0.5410156,1.0108008)(2.5410156,-0.9891992)
\psline[linewidth=0.04cm](2.5410156,-0.9891992)(4.5410156,1.0108008)
\psline[linewidth=0.04cm](2.5410156,-0.9891992)(2.5410156,-1.3891993)
\psline[linewidth=0.04cm](1.3410156,1.0108008)(0.9410156,0.6108008)
\psline[linewidth=0.04cm](2.2210157,0.13080078)(1.8210156,-0.26919922)
\usefont{T1}{ptm}{m}{n}
\rput(2.7124708,-1.1641992){$*_k$}
\psline[linewidth=0.04cm](2.9210157,0.15080078)(3.3210156,-0.24919921)
\psline[linewidth=0.04cm](3.6610155,0.030800782)(2.6610155,1.0308008)
\psline[linewidth=0.04cm](2.4410157,-0.26919922)(2.8610156,-0.6891992)
\psline[linewidth=0.04cm](3.1810157,1.0108008)(2.8810155,0.71080077)
\usefont{T1}{ptm}{m}{n}
\rput(1.4924707,1.2358007){$t_p$}
\usefont{T1}{ptm}{m}{n}
\rput(2.3624706,0.29580078){$t_{h+1}$}
\usefont{T1}{ptm}{m}{n}
\rput(2.4724708,-0.16419922){$w_1$}
\usefont{T1}{ptm}{m}{n}
\rput(2.9724708,0.29580078){$w_q$}
\usefont{T1}{ptm}{m}{n}
\rput(3.1924708,-0.7641992){$*_{j_1}$}
\usefont{T1}{ptm}{m}{n}
\rput(3.6324706,-0.30419922){$*_{j_q}$}
\usefont{T1}{ptm}{m}{n}
\rput(4.6124706,1.2358007){$t^r$}
\usefont{T1}{ptm}{m}{n}
\rput(3.3824706,1.2358007){$t_{h{-}1}$}
\usefont{T1}{ptm}{m}{n}
\rput(3.9324708,-0.044199217){$*_{i_1}$}
\usefont{T1}{ptm}{m}{n}
\rput(2.844707,0.6758008){$*_{i_{h{-}1}}$}
\psdots[dotsize=0.06](1.6810156,0.43080078)
\psdots[dotsize=0.06](1.3610156,0.43080078)
\psdots[dotsize=0.06](1.5210156,0.43080078)
\psdots[dotsize=0.06](2.7410155,-0.26919922)
\psdots[dotsize=0.06](2.9010155,-0.26919922)
\psdots[dotsize=0.06](3.0610155,-0.26919922)
\usefont{T1}{ptm}{m}{n}
\rput(0.85,0.55){$*_{i_{p}}$}
\usefont{T1}{ptm}{m}{n}
\rput(1.6,-0.4){$*_{i_{h+1}}$}
\psdots[dotsize=0.06](3.3410156,0.5508008)
\psdots[dotsize=0.06](3.5010156,0.5508008)
\psdots[dotsize=0.06](3.6610156,0.5508008)
\end{pspicture} 
}

\item 
$$\begin{array}{c}
\varphi ^{-1}(t) = \varphi ^{-1}(\Omega_{i_{h+1},\dots ,i_p}^L(t_{h+1},\dots ,t_p))\vee _{*_k}\qquad\qquad\qquad\hfill \\
\hfill\qquad\qquad\qquad(\varphi^{-1}(w_1)\vee_{*_{j_1}}(\dots (\varphi^{-1}(w_q)\vee _{*_{j_q}}\varphi^{-1}(\Omega _{i_1,\dots ,i_{h{-}1}}^L(t^r,t_2,\dots ,t_{h{-}1})))).\end{array}$$
\end{enumerate}

We have that $\varphi ^{-1}$ is well defined, a tedious but straightforwrad calculation shows that it is the inverse of $\varphi$.
\end{proo}
\medskip

\begin{lemma} \label{prop:DlDmfree} For any integer $0\leq k < m$, let
${\underline r}$ be the composition of $m+1$ in $m$ parts, given by $r_j = 1$ for $j\neq k$, and $r_k = 2$. 
For any set $X$, the image of $\Dy ^m(X)$ under the functor ${\mathbb F}_{\underline r}$ is generated as $\Dy ^{m-1}$ 
algebra by the graded set ${\mathcal A}^{m,k}(X)$ of all colored trees $t$ in 
${\displaystyle \bigcup _{n\geq 1}{\mathcal B}_n^{m,k}(X)}$, such that $n = 1$, or $n >1$ and the root of $t$ is colored by $*_k$. 
\end{lemma}
\medskip

\begin{proo} Again, from the description of $\Dy^m(X)$, we have that it suffices to prove the result for the set with one element $X=\{a\}$.

The $\Dy ^{m-1}$ algebra structure of $\Dy ^m(a)$ is given by the products ${\overline {*}}_j =\begin{cases} *_j,& {\rm for}\ 0\leq j< k,\\
*_{k}+*_{k+1},& {\rm for}\ j = k,\\
*_{j-1},& {\rm for}\ k < j < m.\end{cases}$

The underlying vector space of ${\mathbb F}_{\underline r}(\Dy ^m(a))$ is equal to $\Dy ^m(a)$. As the set ${\displaystyle \bigcup _{n\geq 1}{\mathcal B}_n^{m,k}}$ is a basis of $\Dy ^m(a)$ as a $\K$-vector space, it suffices to see that any element in ${\displaystyle \bigcup _{n\geq 1}{\mathcal B}_n^{m,k}}$ belongs to the $\Dy^{m-1}$ algebra generated by the set ${\mathcal A}^{m,k}$, under the operations ${\overline *}_0,\dots , {\overline *}_{m-1}$.

We proceed by induction on the degree $n$. For $n = 1,2$, the result is immediate.

For $t = t^l\vee _{*_i} t^r\in {\mathcal B}_n^{m,k}(X)$, the recursive recursive hypothesis states that the trees $t^l$ and $t^r$ are obtained by applying the products ${\overline {*}}_0,\dots ,{\overline {*}}_{m-1}$ 
to elements of the set ${\mathcal A}^{m,k}$ of degree smaller than $n$. 

We have to analize three different cases:\begin{enumerate}
\item for $i < k$, we have that $t = t^l {\overline {*}}_i t^r$, and as $t^l$ and $t^r$ are elements in the $\Dy ^{m-1}$ algebra generated by ${\mathcal A}^{m,k}$, so is $t$,
\item for $i = k$, as $t\in {\mathcal B}_n^{m,k}$, we get that $t\in {\mathcal A}^{m,k}$,
\item for $i = k+1$, we have that $t = t^l{\overline {*}}_{k} t^r\ {-}\ t^l*_{k} t^r$ and the root of $t^l$ is colored by $*_j$, with $j > k+1$ or $j =k$.

As $t^l$ and $t^r$ belong to $\Dy ^{m-1}({\mathcal A}^{m,k})$, the tree $t^l{\overline {*}}_{k} t^r$ is in $\Dy ^{m-1}({\mathcal A}^{m,k})$. 

On the other hand, either $t^l*_{k} t^r\in {\mathcal A}^{m,k}$, or $$t^l*_{k} t^r = (t^l\vee_{*_k}w^l)\vee _{*_h} w^r,$$ for some colored tree $w = w^l\vee {*_h} w^r$ and $h > k$.
 
 Applying a recursive argument to the degrees of the elements $t^l\vee_{*_k}w^l$ and $w^r$ the result follows.
\item For $i > k+1$, we have that $t = t^l\vee _{\overline {*}_{i{-}1}}t^r$, which
belongs to $\Dy ^{m-1}({\mathcal A}^{m,k})$ by recursive hypothesis.\end{enumerate}
\end{proo}
\medskip

Lemma \ref{prop:DlDmfree} states that ${\mathbb F}_{\underline r}(\Dy ^m(X))$ is a quotient of the free $\Dy^{m-1}$ algebra $\Dy^{m-1}({\mathcal A}^{m,k}(X))$. For $X$ finite, the subspace of homogeneous elements of degree $n$ in ${\mathbb F}_{\underline r}(\Dy ^m(X))$ is $d_{m,n}\vert X\vert ^n$. 

So, to prove that ${\mathbb F}_{\underline r}(\Dy ^m(X))$ is isomorphic to $\Dy^{m-1}({\mathcal A}^{m,k}(X))$, it suffices to show that the dimension of the subspace of homogeneous elements of degree $n$ in $\Dy^{m-1}({\mathcal A}^{m,k})$ is $d_{m,n}$, where ${\mathcal A}^{m,k}$ is the set of trees in ${\mathcal B}_n^m$ with the vertices colored by $*_0,\dots , *_m$ and the root colored by $*_k$.
\medskip

Recall that for any graded vector space $V=\bigoplus_{n\geq 1} V_n$ such that each $V_n$ is finite dimensional, the generating series of $V$ is $v(x) := {\displaystyle \sum_{n\geq 1} {\mbox{dim}_{\K}(V_n)} x^n}$.

\begin{lemma}\label{lemform}
Let $d_m(x)$ be the generating series of the free $\Dy^m$ algebra $\Dy^m(a)$. We have that:
$$d_m(x)=d_k(x\cdot (1+d_m(x)^{m-k}),$$
for all $0\leq k\leq m$.
\end{lemma}
\medskip

\begin{proo}
Clearly, it is enough to prove this for $k=m-1$. Let $g_m(x)$ be the inverse series of $d_m(x)$ ($g$ exists because $d(0)=0$). 

Since $x\cdot (1+d_m(x))^{m+1}=d_m(x)$, replacing $x$ by $g_m(x)$ we obtain that:
$$g_m(x)=\displaystyle\frac{x}{(1+x)^{m+1}},$$
which implies that  $(1+x)\cdot g_m(x) = g_{m-1}(x)$. So, replacing $x$ by $d_m(x)$ and applying $d_{m-1}(x)$ to both sides, we get the desired formula
$$d_{m-1} (x\cdot (1+d_m(x))) = d_m(x).$$
\end{proo}
\medskip

Applying Lemmas \ref{prop:DlDmfree} and \ref{lemform}, we get the following result.

\begin{proposition} \label{pro:freenessform} For a fixed $0\leq k\leq m-1$, let ${\underline r}$ be the composition of $m+1$ in $m$ parts, such that $r_i=1$ for $i\neq k$ and $r_k=2$. The $\Dy ^{m-1}$ algebra ${\mathbb F}_{\underline r}(\Dy ^m(X))$ is free.\end{proposition}
\medskip

\begin{proo} Applying Lemmas \ref{prop:DlDmfree} and \ref{lemform}, it suffices to prove that the number of elements in ${\mathcal A}_n^{m,k}$ is $d_{m,n-1}$, for $0\leq k\leq m$. 

The number of elements of ${\mathcal B}_{n-1}^{m,k}$ is $d_{m,n-1}$, to end the proof we define a bijective map $\theta _n$ from ${\mathcal B}_{n-1}^{m,k}$ to ${\mathcal A}_n^{m,k}$, for $n\geq 2$. 

For $n = 2$, $\theta _1(\vert )$ is the unique planar binary rooted tree with two leaves and the root colored by $*_k$.

Let $t =t^l\vee_{*_{h}}t^r$ be an element of ${\mathcal B}_{n-1}^{m,k}$. \begin{enumerate}
\item For $h > k$, let $t = \Omega_{h,i_2,\dots,i_p}^L(t^r,t_2,\dots ,t_p)$.\begin{enumerate}
\item If $i_p >\dots > i_2 > h > k$, then we define 
$\theta _n(t) := t\vee _{*_k} \vert .$
\item If there exists one integer $1\leq s\leq p$ such that $i_s = k$, then the $s$ is unique and $\theta _n(t)$ is defined by the formula:
$$\theta _n(t) :=\Omega_{h,i_2,\dots, i_{s{-}1}}^L(t^r,t_2,\dots ,t_{s{-}1})\vee_{*_k}\Omega_{k,i_{s+1},\dots ,i_p}^L(t_s,\dots, t_p).$$\end{enumerate}
\item For $h\leq k$, let $t = \Omega _{h,j_2,\dots ,j_q}^R(t^l,w_2,\dots ,w_q)$.\begin{enumerate}
\item If $j_i\leq k$ for any $2\leq i\leq q$, then we define:
$\theta _n(t) := \vert \vee _{*_k} t$.
\item Otherwise, let $2\leq s\leq q$ be the minimal integer such that $j_s >k$. In this case, as $t\in {\mathcal B}^{m,k}$, we know that $k\notin \{h,j_1,\dots ,j_{s{-}1}\}$.
We define $\theta _n(t)$ to be the element:
 $$ \theta _n(t):= \Omega_{j_s,\dots ,j_q}^R(w_s,\dots ,w_q)\vee_{*_k}\Omega_{h,j_2,\dots ,j_{s{-}1}}^R(t^l,w_2,\dots ,w_{s{-}1}.$$
\end{enumerate}
\end{enumerate}

It is not difficult to verify that $\theta _n$ is bijective for all $n\geq 2$. So, the result is proved.
\end{proo}

Applying Lemma \ref{lemfinctorsFr} , as a straightforward consequence of Proposition \ref{pro:freenessform}, we get the following result.

\begin{theorem} Let $0\leq i\leq m-1$ be an integer and let ${\underline r}$ be a composition of $m+1$ in $s+1$ parts. The image of a free $\Dy ^m$ algebra $\Dy ^m(X)$ under the functor ${\mathbb F}_{\underline r}$ is a free $\Dy ^{s}$ algebra.\end{theorem}
\medskip

Note that, in particular we get that, for any free $\Dy ^m$ algebra, the associative algebra $(\Dy ^m(X), *_0+\dots +*_m)$ is free.

\bigskip

\bigskip

\section{A diagonal on $m$-Dyck paths}
\label{coproductonmdyckpaths}
\medskip

As  $\Dy ^m$ is a regular operad,  given a $\Dy^m$ algebra $(A, \{ *_i\}_{0\leq i\leq m})$ and an associative algebra $(B, \circ )$, the tensor product $B\ot A$ has a natural structure of  $\Dy ^m$ algebra, where the products are given by the formula $*_i^{B\ot A} := \circ \ot *_i$, for $0\leq i\leq m$.  In particular, when $(B,\circ ) = (A, * :={\displaystyle \sum_{i=0}^m *_i})$, the tensor product $A\ot A$ is a $\Dy^m$ algebra. That is, the algebraic operad $\Dy ^m$ is a Hopf operad. 

However, there does not exist a good notion of unit for $\Dy ^m$ algebras, when $m\geq 1$. 

In this section, we introduce the notion of $\Dy ^m$ bialgebra, and give an explicit description of the coproduct on the free algebra $\D_m$, for $m\geq 1$. For $m=1$ it coincides, via the linear map induced by the applications $\Gamma _n:\Dy _n^1\longrightarrow {\mathcal Y}_n$, with the coproduct defined in \cite{LodRon} on the algebra $\K[{\mathcal Y}_*]$ of planar binary rooted trees.
\medskip

Given a vector space $V$, recall that $V^+$ is the vector space  $V^+:= \K \oplus V$ equipped with the usual augmentation map $\epsilon: V^+\longrightarrow \K$. 
Let ${\overline {V^+\ot V^+}}$ denote the vector space ${\overline {V^+\ot V^+}}:= V^+\ot V \oplus V\ot V^+$.
\medskip

Let $(A, \{*_i\}_{0\leq i\leq m})$ be a $\Dy^m$ algebra. The products $*_i$ are extended to linear maps
 $*_i: {\overline {A^+\ot A^+}} \longrightarrow A$, for $0\leq i\leq m$, by the formulas:
\begin{enumerate}
 \item $x*_0 1_{\K}=0$ and $1_{\K}*_0 x = x$,
 \item $x *_i 1_{\K} = 1_{\K}*_i x = 0$, for $0 < i < m$,
 \item $x*_m 1_{\K} = x$ and $1_{\K} *_m x = 0$,
\end{enumerate}
for $x\in A$.
\medskip

Note that the element $1_{\K}*_i1_{\K}$ is  not defined, for any $0\leq i \leq m$. 
\medskip

It is easily seen that the vector space ${\overline {A^+\ot A^+}}$, equipped with the operations $*_i$ given by: 
\begin{enumerate}
 \item $(x_1\ot x_2) *_i (y_1\ot y_2)=(x_1*y_1)\ot (x_2*_iy_2 )$, for $x_2\in A$ or $y_2\in A$;
 \item $(x_1\ot 1_{\K}) *_i (y_1\ot 1_{\K})=(x_1*_i y_1)\ot 1_{\K}$,
\end{enumerate}
for $x_1,x_2, y_1, y_2 \in A^+$, is a $\Dy^m$ algebra.
\bigskip

The previous construction motivates the following definition.

\begin{definition}\label{Dyckbialgebra}  A $\Dy ^m$ bialgebra over $\K$ is a $\Dy^m$ algebra $(A, \{*_i\}_{0\leq i\leq m})$ equipped with a linear map
$\Delta : A^+ \longrightarrow A^+\ot A^+$ satisfying that:\begin{enumerate}
\item the data $(A^+, * , \Delta ,\iota ,\epsilon )$ is a bialgebra in the usual sense, where\begin{enumerate}
\item the associative product $*$ is given by:
$$x* y:=\begin{cases} \sum _{i=0}^m x*_i y,& {\rm for}\ x,y\in A,\\
x\cdot y,& {\rm for}\ x\in\K\ {\rm or}\ y\in \K,\end{cases}$$
where $\cdot $ denotes indistinctly the product on $\K$ as well as the action of $\K$ on $A$, for $x,y\in A^+$.
\item $\iota :\K\hookrightarrow A^+$ is the canonical inclusion of $\K$ into $A^+$, and $\epsilon :A^+\longrightarrow \K$ is the canonical projection.\end{enumerate}
\item the restriction of $\Delta $ from $A$ to the subspace ${\overline {A^+\ot A^+}}$ is a homomorphism of $\Dy ^m$ algebras.\end{enumerate}\end{definition}
\bigskip

A standard argument shows that for any free $\Dy ^m$ algebra $\Dy^m(X)$, there exists a unique homomorphism 
$\Delta $ from $\Dy^m(X)^+$ into $\Dy^m(X)^+\ot \Dy^m(X)^+$ satisfying that:\begin{enumerate}
\item $\Delta (1_{\K}) = 1_{\K}\ot 1_{\K}$,
\item $\Delta (x) = x\ot 1_{\K} + 1_{\K}\ot x$, for $x\in X$,\end{enumerate}
giving $\Dy^m(X)$ a structure of $\Dy^m$ bialgebra.
\bigskip

Our aim is to give an explicit description, in terms of $m$-Dyck paths, of the coproduct $\Delta$ on the free $\Dy^m$ algebra $\D_m$.
\medskip

\begin{definition}\label{def:centralvertex}
Let $P$ be a $m$-Dyck path. A {\em central step} of $P$ is an up step of $P$ which is the initial step of $P$, or is preceded by another up step. \end{definition}
\medskip

\begin{example}
Consider the following $2$-Dyck path:

\medskip
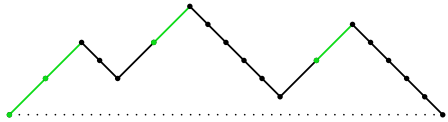
\scalebox{0.6} 
{
\begin{pspicture}(0,-1.28)(9.74,1.28)
\definecolor{color31}{rgb}{0.050980392156862744,0.8588235294117647,0.09411764705882353}
\psline[linewidth=0.04cm,linecolor=color31](0.06,-1.2)(0.86,-0.4)
\psline[linewidth=0.04cm,linecolor=color31](0.86,-0.4)(1.66,0.4)
\psline[linewidth=0.04cm](1.66,0.4)(2.06,0.0)
\psline[linewidth=0.04cm](2.06,0.0)(2.46,-0.4)
\psline[linewidth=0.04cm](2.46,-0.4)(3.26,0.4)
\psline[linewidth=0.04cm,linecolor=color31](3.26,0.4)(4.06,1.2)
\psline[linewidth=0.04cm](4.06,1.2)(6.06,-0.8)
\psline[linewidth=0.04cm](6.06,-0.8)(6.86,0.0)
\psline[linewidth=0.04cm,linecolor=color31](6.86,0.0)(7.66,0.8)
\psline[linewidth=0.04cm](7.66,0.8)(9.66,-1.2)
\psline[linewidth=0.04cm,linestyle=dotted,dotsep=0.16cm](0.06,-1.2)(9.66,-1.2)
\psdots[dotsize=0.12](0.06,-1.2)
\psdots[dotsize=0.12](0.86,-0.4)
\psdots[dotsize=0.12](1.66,0.4)
\psdots[dotsize=0.12](2.06,0.0)
\psdots[dotsize=0.12](2.46,-0.4)
\psdots[dotsize=0.12](3.26,0.4)
\psdots[dotsize=0.12](4.06,1.2)
\psdots[dotsize=0.12](4.46,0.8)
\psdots[dotsize=0.12](4.86,0.4)
\psdots[dotsize=0.12](5.26,0.0)
\psdots[dotsize=0.12](5.66,-0.4)
\psdots[dotsize=0.12](6.06,-0.8)
\psdots[dotsize=0.12](6.86,0.0)
\psdots[dotsize=0.12](7.66,0.8)
\psdots[dotsize=0.12](8.06,0.4)
\psdots[dotsize=0.12](8.46,0.0)
\psdots[dotsize=0.12](8.86,-0.4)
\psdots[dotsize=0.12](9.26,-0.8)
\psdots[dotsize=0.12](9.66,-1.2)
\psdots[dotsize=0.12,linecolor=color31](0.06,-1.2)
\psdots[dotsize=0.12,linecolor=color31](0.86,-0.4)
\psdots[dotsize=0.12,linecolor=color31](3.26,0.4)
\psdots[dotsize=0.12,linecolor=color31](6.86,0.0)
\end{pspicture} 
}

\medskip

The central steps are marked in green.
\end{example}
\medskip

\begin{notation} \label{not:Puv} Let $P$ be an $m$-Dyck path of size $n$. Given a pair of steps $(u,d)\in \up(P)\t \dw(P)$, such that the starting vertex of $u$ and the final vertex of $d$ belong to the same horizontal line, we denote by $P_{u,d}$ the (translated) $m$-Dyck path obtained from $P$ which starts with $u$ and ends with $d$.\end{notation}

\begin{definition}\label{def:admissiblecut} A {\it cut } of $P$ is an $m$-Dyck path $P_{u,d}$ such that $u$ is a central step of $P$ and $P_{u,d}\neq P$. An {\it admissible cutting} of $P$ is a non-empty family of cuts ${\mathcal P}= \{ P_{u_l,d_l}\}_{1\leq l\leq s}$ of $P$ such that $P_{u_l,d_l}$ and $P_{u_h,d_h}$ are disjoint whenever $l\neq h$.
\end{definition}
\medskip

\begin{remark}\label{rem:cuts} For any central step $u$ of $P$, the excursion $P_{u,w_u}$ of $u$ in $P$ is a cut of $P$. \end{remark}

\begin{notation}\label{not:ordercutting} Let ${\mathcal P}= \{ P^1,\dots ,P^s\}$ be an admissible cutting of an $m$-Dyck path $P$, such that $P^l=P_{u_l,d_l}$,  for $1\leq l\leq s$. 

Suppose that for any $1\leq l\leq s$, the starting vertex of $u_l$ has coordinates $(a_l,b_l)$ and the final vertex of $d_l$ is $(c_l,d_l)$, we shall always assume that ${\mathcal P}$ is ordered in such a way that $a_1 < a_2<\dots < a_s$, which implies that:
$$a_1 < c_1 < a_2< c_2 <\dots <a_s < c_s.$$\end{notation}
\medskip

\begin{example} \label{example:admiscut}
Consider the Dyck path of the preceding example. The admissible cuts are the paths above the dotted red lines.
\medskip

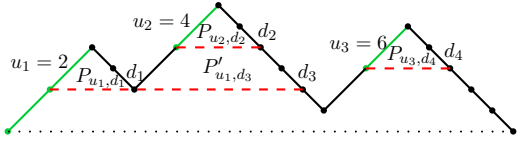
\scalebox{0.7} 
{
\begin{pspicture}(0,-1.28)(9.74,1.28)
\definecolor{color40}{rgb}{0.0,0.8,0.2}
\psdots[dotsize=0.12,linecolor=color40](0.06,-1.2)
\psdots[dotsize=0.12,linecolor=color40](0.86,-0.4)
\psline[linewidth=0.04cm,linecolor=color40](0.06,-1.2)(1.66,0.4)
\psdots[dotsize=0.12,linecolor=color40](3.26,0.4)
\psline[linewidth=0.04cm,linecolor=color40](3.26,0.4)(4.06,1.2)
\psdots[dotsize=0.12,linecolor=color40](6.86,0.0)
\psline[linewidth=0.04cm,linecolor=color40](6.86,0.0)(7.66,0.8)
\psdots[dotsize=0.12](1.66,0.4)
\psdots[dotsize=0.12](2.06,0.0)
\psdots[dotsize=0.12](2.46,-0.4)
\psdots[dotsize=0.12](4.06,1.2)
\psdots[dotsize=0.12](4.46,0.8)
\psdots[dotsize=0.12](4.86,0.4)
\psdots[dotsize=0.12](5.26,0.0)
\psdots[dotsize=0.12](5.66,-0.4)
\psdots[dotsize=0.12](6.06,-0.8)
\psdots[dotsize=0.12](7.66,0.8)
\psdots[dotsize=0.12](8.06,0.4)
\psdots[dotsize=0.12](8.46,0.0)
\psdots[dotsize=0.12](8.86,-0.4)
\psdots[dotsize=0.12](9.26,-0.8)
\psdots[dotsize=0.12](9.66,-1.2)
\psline[linewidth=0.04cm](1.66,0.4)(2.46,-0.4)
\psline[linewidth=0.04cm](2.46,-0.4)(3.26,0.4)
\psline[linewidth=0.04cm](4.06,1.2)(6.06,-0.8)
\psline[linewidth=0.04cm](6.06,-0.8)(6.86,0.0)
\psline[linewidth=0.04cm](7.66,0.8)(9.66,-1.2)
\psline[linewidth=0.04cm,linestyle=dotted,dotsep=0.16cm](0.06,-1.2)(9.66,-1.2)
\psline[linewidth=0.04cm,linecolor=red,linestyle=dashed,dash=0.16cm 0.16cm](0.86,-0.4)(5.66,-0.4)
\psline[linewidth=0.04cm,linecolor=red,linestyle=dashed,dash=0.16cm 0.16cm](3.26,0.4)(4.86,0.4)
\psline[linewidth=0.04cm,linecolor=red,linestyle=dashed,dash=0.16cm 0.16cm](6.86,0.0)(8.46,0.0)
\usefont{T1}{ptm}{m}{n}
\rput(2.491455,-0.075){$d_1$}
\usefont{T1}{ptm}{m}{n}
\rput(0.6614551,0.105){$u_1=2$}
\usefont{T1}{ptm}{m}{n}
\rput(2.981455,0.885){$u_2 = 4$}
\usefont{T1}{ptm}{m}{n}
\rput(4.991455,0.705){$d_2$}
\usefont{T1}{ptm}{m}{n}
\rput(5.7514553,-0.075){$d_3$}
\usefont{T1}{ptm}{m}{n}
\rput(6.721455,0.465){$u_3=6$}
\usefont{T1}{ptm}{m}{n}
\rput(8.531455,0.325){$d_4$}
\usefont{T1}{ptm}{m}{n}
\rput(1.8351074,-0.175){$P_{u_1,d_1}$}
\usefont{T1}{ptm}{m}{n}
\rput(4.241455,-0.015){$P\rq_{u_1,d_3}$}
\usefont{T1}{ptm}{m}{n}
\rput(4.121455,0.645){$P_{u_2,d_2}$}
\usefont{T1}{ptm}{m}{n}
\rput(7.741455,0.245){$P_{u_3,d_4}$}
\end{pspicture} 
}

\medskip

Observe that the cuts $P_{u_1,d_1}$ and $P_{u_1,d_3} = P_{u_1,d_1}\t_0 P_{3,d_3}$ (where $3$ denotes the third step of $P$) begin both with $u_1$, so there are two admissible cuts
corresponding to the lowest red dotted line. 

The admissible cuttings of $P$ are $\{ P_{u_1,d_1}\}$, $\{ P_{u_1, d_3}\}$, $\{ P_{u_2, d_2}\}$, $\{ P_{u_3, d_4}\}$, 
$\{ P_{u_1, d_1}, P_{u_2,d_2}\}$, $\{ P_{u_1, d_1}, P_{u_3,d_4}\}$, $\{ P_{u_2, d_2}, P_{u_3,d_4}\}$, $\{ P_{u_1, d_3}, P_{u_3,d_4}\}$ and 

\noindent $\{ P_{u_1, d_1}, P_{u_2,d_2}, P_{u_3,d_4}\}$.
\end{example}

Let ${\mbox {Ad}(P)}$ denote the set of admissible cuttings of $P$. 
\medskip

\begin{notation} \label{def:collapsing} Let $P$ be an element of $\Dy_n^m$. For any cut $P_{u,d}$ of $P$, denote by $P/ P_{u,d}$ the Dyck path obtained from replacing the path $P_{u,d}$ by a point in $P$, that is, by taking off all the steps of $P_{u,d}$ and gluing the initial vertex of $u$ with the final vertex of $d$.\end{notation}

\begin{remark} $(1)$ Suppose that $P=P_1\t_0\dots \t_0P_r$, with $P_i$ prime for $1\leq i\leq r$. The Dyck paths $P_1$, $P_1\t_0P_2$, \dots , $P_1\t_0\dots \t_0P_{r-1}$ are cuts of $P$.

$(2)$ Let ${\mathcal P} = \{P^1,\dots , P^s\}\in {\mbox {Ad}(P)}$ be an admissible cutting of $P$. For any $1\leq l\leq s$, the collection ${\mathcal P}\setminus \{ P^l\}$ is an admissible cutting of $P/P^l$.\end{remark}

For any admissible cutting ${\mathcal P} = \{P^1,\dots , P^s\}$ of a path $P\in \Dy _n^m$, the $m$-Dyck $P/\{P^1,\dots ,P^s\}$ is defined recursively by the formula:
$$P/\{P^1,\dots ,P^s\} := (P/\{P^1,\dots ,P^{s-1}\})/ \{P^s\}.$$

\begin{example} Let $P = (0, 2, 0, 5, 0,5)\in \Dy _6^2$ be the path of Example \ref{example:admiscut}, and consider the admissible cutting  $\{ P_{u_1, d_3}\}$ of $P$. The $2$-Dyck path $P/\{ P_{u_1, d_3}\}$ is the path:

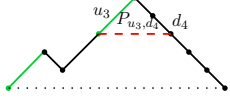
\scalebox{0.6} 
{
\begin{pspicture}(0,-1.08)(4.94,1.08)
\definecolor{color40}{rgb}{0.0,0.8,0.2}
\psdots[dotsize=0.12,linecolor=color40](0.06,-1.0)
\psdots[dotsize=0.12,linecolor=color40](0.86,-0.2)
\psline[linewidth=0.04cm,linecolor=color40](0.06,-1.0)(0.84,-0.2)
\psdots[dotsize=0.12,linecolor=color40](2.06,0.2)
\psline[linewidth=0.04cm,linecolor=color40](2.06,0.2)(2.86,1.0)
\psdots[dotsize=0.12](1.26,-0.6)
\psdots[dotsize=0.12](0.86,-0.2)
\psdots[dotsize=0.12](2.86,1.0)
\psdots[dotsize=0.12](3.26,0.6)
\psdots[dotsize=0.12](3.66,0.2)
\psdots[dotsize=0.12](4.06,-0.2)
\psdots[dotsize=0.12](4.46,-0.6)
\psdots[dotsize=0.12](4.86,-1.0)
\psline[linewidth=0.04cm](1.26,-0.6)(2.06,0.2)
\psline[linewidth=0.04cm](2.86,1.0)(4.86,-1.0)
\psline[linewidth=0.04cm,linestyle=dotted,dotsep=0.16cm](0.06,-1.0)(4.86,-1.0)
\psline[linewidth=0.04cm,linecolor=red,linestyle=dashed,dash=0.16cm 0.16cm](2.06,0.2)(3.66,0.2)
\usefont{T1}{ptm}{m}{n}
\rput(2.1314552,0.665){$u_3$}
\usefont{T1}{ptm}{m}{n}
\rput(3.8714552,0.485){$d_4$}
\usefont{T1}{ptm}{m}{n}
\rput(2.9414551,0.485){$P_{u_3,d_4}$}
\psline[linewidth=0.04cm,linestyle=dashed,dash=0.16cm 0.16cm](0.86,-0.2)(0.88,-0.24)
\psline[linewidth=0.04cm](0.86,-0.2)(1.26,-0.6)
\end{pspicture} 
}

\end{example}
\medskip

\begin{definition}\label{def:coprodpaths}
The (reduced) coproduct $\ov{\Delta}:\D_m\to \D_m\ot \D_m$ on $\D_m$ is defined by the following formula:
$$\ov{\Delta}(P)=\displaystyle \sum_{\mathcal P} P^1*\dots *P^s\ot P/\{P^1,\dots,P^s\},$$
for any $P\in \Dy_n^m$, where the sum ranges over all the admissible cuttings ${\mathcal P} = \{P^1,\dots, P^s\}\in Ad(P)$.
\end{definition}
\medskip

\begin{remark} \label{rem:L(P)mayorm} For any $m$-Dyck path $P\in \Dy _n^m$ and any admissible cutting ${\mathcal P} =\{P^1,\dots ,P^s\}$ of $P$ such that $P^i\in \Dy_{n_i}^m$, we have that $P/{\mathcal P}\in \Dy_N^m$ is a Dyck path, with $N = {n-n_1-\dots -n_s}$, so $L(P/{\mathcal P})\geq m$.\end{remark}
\medskip

The reduced coproduct extends to a coproduct $\Delta:\D_m\to \Dot2$ defining
$$\Delta(P)= P\ot 1_{\K} + \ov{\Delta}(P)+1_{\K}\ot P.$$

\begin{notation} \label{not:diag} Let $P$ be an $m$-Dyck path, \begin{enumerate} \item we use Sweddler\rq s notation for the coproduct, that is
$$\Delta (P) = \sum P_{(1)}\ot P_{(2)},$$
for any $P\in \D _m$, to denote the image of $P$ under the coproduct,
\item the image of $P$ under the reduced coproduct is denoted 
$${\ov{\Delta}}(P)=\sum {\ov {P}}_{(1)}\ot {\ov {P}}_{(2)},$$
\item for any integer $0\leq j\leq L(P)$, we denote by $\de_{L\geq j}(P)$ (respectively, $\de_{L=j}(P)$) 
the sum of the terms $P_{(1)}\ot P_{(2)}$ appearing in $\de(P)$ such 
that $L(P_{(2)})\geq j$ (respectively, $L(P_{(2)})= j$). 

\noindent We write $\sum P_{(1)}^{L\geq j}\ot P_{(2)}^{L\geq j}$ for $\de_{L\geq j}(P)$ (and similarly for $\de_{L=j}(P)$).

\item for the reduced coproduct, we denote ${\ov{\Delta}} _{L\geq j}(P) = \sum {\ov P}_{(1)}^{L\geq j}\ot {\ov P}_{(2)}^{L\geq j}$ (respectively, $ {\ov{\Delta}} _{L=j}(P) = 
\sum {\ov P }_{(1)}^{L= j}\ot {\ov P }_{(2)}^{L= j}$),
\item given an admissible cutting ${\mathcal P}=\{ P^1,\dots ,P^s\}$ of $P$, we use ${\mathcal P}_{(1)}$ to denote the sum of elements $P^1*\dots *P^s$ and ${\mathcal P}_{(2)}$ for the element $P/\{\mathcal P\}$.\end{enumerate}
\end{notation}
\medskip

From Remark \ref{rem:L(P)mayorm} we get that $\de_{L\geq m}(P) = \de(P) {-} P\ot 1_{\K},$
 for any Dyck path $P\in \Dy_n^m$. 
\bigskip

The main result of this section is the following Theorem.
\medskip

\begin{theorem}\label{theorem: diagonal respect the *_i products}
The coproduct $\Delta$ defined on $\D_m$ satisfies the relation:
$$\Delta(P*_i Q)=\Delta(P)*_i\Delta(Q) = \sum (P_{(1)} * Q_{(1)})\ot (P_{(2)}*_iQ_{(2)}),$$
for any integer $0\leq i\leq m$ and any pair of elements $P,Q\in \D_m$. In other words, the triple $(\D_m, \{*_i\}_{0\leq i\leq m}, \Delta)$ is a $\Dy^m$ bialgebra.
\end{theorem}
\medskip

The proof of Theorem \ref{theorem: diagonal respect the *_i products} requires to prove some additional results first. Let us begin by extending the products $\t_j$, defined in Section 2, to the ${\overline {\D _m^+\ot \D_m^+}}$ in a trivial way. 
\medskip 

\begin{definition} \label{def:timesjintensor} For any pair of $m$-Dyck paths $P$ and $Q$, and any integer $0\leq j\leq L(P)$, define:\begin{enumerate}
\item $P\t_j 1_{\K} := \begin{cases}0, & {\rm for}\ 0\leq j < L(P)\\
P, & {\rm for}\  j=L(P).\end{cases}$
\item $1_{\K}\t_jP := \begin{cases}0, & {\rm for}\ j > 0\\
P, & {\rm for}\  j=0.\end{cases}$
\item $(P*Q)\ot (1_{\K}\t_j1_{\K}): =(P\t_jQ)\ot 1_{\K}$.
\end{enumerate}
Extending by linearity, we get a well defined product $\t_j$ on $\overline{\D_m^+\ot\D_m^+}$, given by
$$(P\ot Q) \t_j (R\ot S) = (P*R)\ot Q\t_jS.$$
\end{definition}
\medskip

\begin{lemma}\label{lemma: the lambda set of P_2 can be obtained from the lambda set of P}
Let $P$ be an $m$-Dyck path and $P_{(2)}$ the result of collapsing a set of admissible cuts of $P$ to a point. 
 For $0<i<m$, we have $C_i(P)=C_i(P_{(2)})$ and $c_i(P)=c_i(P_{(2)})$. In particular, $C_{m-1}(P)<L(P_{(2)})$
\end{lemma}
\medskip

\begin{proof}
 Observe that the down steps of maximal level of $P_{(2)}$ are the last $L(P_{(2)})$ down steps of $P$ and 
 the colors of both differ only by a renaming of colors. Therefore, for $0<i<m$, we have that $C_i(P)=C_i(P_{(2)})$ and $c_i(P)=c_i(P_{(2)})$. Also, since $P_{(2)}$ is an $m$-Dyck path, it must have 
 a color repeated $m$ times, this implies that $C_{m-1}(P)=C_{m-1}(P_{(2)})<L(P_{(2)})$. 
\end{proof}
\medskip

\begin{proposition} \label{prop: formula for the diagonal *0} Let $P$ be an $m$-Dyck path and $Q$ a prime $m$-Dyck path. 
The coproduct $\Delta$ satisfies that 
$\Delta(P\t_0 Q)=\de(P)\t_0\de(Q).$ Moreover, we have that $\Delta(P*_0 Q)=\de(P)*_0\de(Q)$.\end{proposition}
\medskip

\begin{proo} Since $Q$ is prime, a cut of $P\times_0 Q$ is either a cut of $P$, or $P$ itself, or a cut of $Q$. So, an admissible cutting ${\mathcal R}$ of $P\t_0 Q$ satisfies one of the following conditions:\begin{enumerate}[(a)]
\item ${\mathcal R}\in {\mbox {Ad}(P)}$ and $(P\t_0Q )/{\mathcal R} = (P/{\mathcal R})\t_0 Q$, or ${\mathcal R} = \{P\}$ and 

\noindent $(P\t_0Q )/{\mathcal R} = Q$,
\item ${\mathcal R}\in {\mbox {Ad}(Q)}$ and $(P\t_0Q )/{\mathcal R} = P\t_0 (Q/{\mathcal R})$,
\item ${\mathcal R} = \{ P\}\cup {\mathcal Q}$, with ${\mathcal Q}\in {\mbox {Ad}(Q)}$. In this case 
$(P\t_0Q )/{\mathcal R} = Q/{\mathcal Q}$,
\item ${\mathcal R} = {\mathcal P}\cup {\mathcal Q}$, for a pair of admissible cuttings $ {\mathcal P}\in {\mbox {Ad}(P)}$ and $ {\mathcal Q}\in {\mbox {Ad}(Q)}$, and we get $(P\t_0Q )/{\mathcal R} = (P/{\mathcal P})\t_0  (Q/{\mathcal Q})$.\end{enumerate}
\medskip

Computing $\Delta (P\t_0Q)$, we get that:
\begin{align}
\Delta (P\t_0 & Q ) = P\t_0 Q\ \ot\  1_{\K} +  \sum {\ov P}_{(1)}\ \ot \ {\ov P}_{(2)}\t _0 Q\ +\nonumber \\
&\ P\ot Q\ + \sum  {\ov Q}_{(1)}\ \ot\ P\t _0 {\ov Q}_{(2)}\ +\ P*{\ov Q}_{(1)}\ \ot\ {\ov Q}_{(2)}\ +\nonumber \\
&\sum  {\ov {P}}_{(1)}*{\ov Q}_{(1)}\ \ot {\ov P}_{(2)}\t _0 {\ov Q}_{(2)}\ +\ 1_{\K}\ot \ P \t_0 Q.\nonumber 
\end{align}

Using that $\ov{\de}(P) \t_0 Q\ot 1_{\K} = 0$, we obtain $\Delta (P\t_0 Q ) = \de (P)\t _0 \de(Q).$

As $Q$ is prime, any $Q_{(2)}$ appearing in $\Delta(Q)$ is also prime.

 So, $P*_0Q=P\t_0Q$ and
$P_{(2)}*_0Q_{(2)}=P_{(2)}\t_0Q_{(2)}$, which implies that $\Delta (P*_0 Q ) = \de (P)*_0 \de(Q).$
\end{proo}
\medskip

\begin{lemma}\label{lemma: 2 formulas for the diagonal on m-Dyck paths}
Let $P,Q$ be two $m$-Dyck paths, with $P\in \Dy_{n_1}^m$ and $Q$ prime, and let $j$ be an integer $0< j \leq L(P)$. The coproduct $\Delta$ on the elements $P\t _jQ$ fulfills the following relation:
\begin{enumerate} 
\item if $0< j < L(P)$, then
\begin{multline}\nonumber
 \Delta(P\t_j Q)=  \ P\t_j Q\ \ot\  1_{\K} \  + \de_{L\geq j}(P)\t_j\de(Q) - \sum P_{(1)}^{L=j}* Q\ \ot \ P_{(2)}^{L=j} \\
 +\displaystyle{ \sum _{L({\mathcal P}_{(2)})\leq j} P^1*\dots *P^{s-1}*(P^s\t _{j{-}L({\mathcal P}_{(2)})}Q)\ 
\ot {\mathcal P}_{(2)}}
\end{multline}
where the sum is taken over all admissible cuttings ${\mathcal P} = \{P^1,\dots ,P^s\}$ of $P$ such that $L({\mathcal P}_{(2)})\leq j$.
\item if $j = L(P)$, then
\begin{multline}\nonumber
\Delta(P\t_{L(P)} Q)=  \ P\t_{L(P)} Q\ \ot\  1_{\K}  + \ \de_{L= L(P)}(P)\t_{L(P)}\de(Q)\\
  +\displaystyle{\sum_{L({\mathcal P}_{(2)}))< L(P)} P^1*\dots *P^{s-1}*(P^s\times_{L(P)-L({\mathcal P}_{(2)})}Q)\ot {\mathcal P}_{(2)}},
\end{multline}
where the sum is taken over all admissible cuttings ${\mathcal P} = \{P^1,\dots ,P^s\}$ of $P$ such that $L({\mathcal P}_{(2)}) < L(P)$.
\end{enumerate}
\end{lemma}
\medskip

\begin{proo} For $0< j\leq L(P)$, a cut $R$ of $P\t_j Q$ is of the form:\begin{enumerate}[(i)]
\item $R$ is a cut of $P$ such that $L(P/R)\geq j$. Note that it means that either the level of the last step of $R$ is smaller than $n_1$, or $R=P_{u,d_k^P}$ with $k\leq L(P)-j$
for some $d = d_k^P \in \dw_{n_1}(P)$.
\item $R$ is a cut of $Q$, for $j < L(P)$. For $j = L(P)$, $R$ is a cut of $Q$ or $Q$. 
\item $R = P_{u,d_k^P}\t_{k-L(P)+j}Q$, for some $d_k^P \in \dw_{n_1}(P)$ such that 

\noindent $k\geq L(P)-j$.
\end{enumerate}
So, any possible admissible cutting ${\mathcal R}$ of $P\t_jQ$ satisfies one of the following conditions:\begin{enumerate}[(a)]
\item ${\mathcal R}\in {\mbox {Ad}(P)}$ is such that  $L(P/{\mathcal R})\geq j$, and $(P\t_jQ )/{\mathcal R} = (P/{\mathcal R})\t_j Q.$
\item ${\mathcal R}\in {\mbox {Ad}(Q)}$, for $1\leq j\leq L(P)$, respectively ${\mathcal R} = \{Q\}$, for $j = L(P)$. 

\noindent In this case, $(P\t_jQ )/{\mathcal R} = P\t _j Q/{\mathcal R}$, respectively $(P\t_jQ )/\{ Q\} = P$.
\item ${\mathcal R}$ is the disjoint union of ${\mathcal P}\in {\mbox {Ad}(P)}$, such that $L(P/{\mathcal P})\geq j$, and ${\mathcal Q}\in {\mbox {Ad}(Q)}$, which does not contain $Q$. So, $(P\t_jQ )/{\mathcal R} = P/{\mathcal P} \t_j Q/{\mathcal Q}$.
\item ${\mathcal R} = \{ P^1,\dots , P^{s-1}, P^s\t _{j{-}L({\mathcal P}_{(2)})}Q\}$, where ${\mathcal P} = \{P^1,\dots ,P^s\}\in {\mbox{Ad}(P)}$ is such that 
$\begin{cases}L({\mathcal P}_{(2)})\leq j&\ {\rm for}\ j<L(P),\\
L({\mathcal P}_{(2)})< L(P)&\ {\rm for}\ j=L(P).\end{cases}$
\medskip

For the previous two cases, we get that $(P\t_jQ)/{\mathcal R}=P/{\mathcal P}$.\\ 
\item ${\mathcal R} = \{ P^1,\dots , P^{s-1}, P^s, Q\}$ for $j=L(P)$, where 

\noindent ${\mathcal P} = \{P^1,\dots ,P^s\}\in {\mbox{Ad}(P)}$ is such that $L({\mathcal P}_{(2)})= L(P)$. Again, we get that $(P\t_jQ)/{\mathcal R}=P/{\mathcal P}$.\\ 
\end{enumerate}
For any pair of Dyck paths $R, S$, we have that:

$\begin{cases}  1_{\K}\t_j R = 0,&\ {\rm for}\ 0 < j \leq L(R),\\
R\t _j Q = 0,&\ {\rm for}\ L(R) < j,\\ 
R\t_j 1_{\K} = 0,&\ {\rm for}\ 0 < j< L(R).\end{cases}$
\medskip

An easy calculation shows that
\begin{multline}
 \Delta (P\t_jQ) = P\t_j Q \ot 1_{\K} \ +\  \sum P_{(1)}^{L\geq j}\ot P_{(2)}^{L\geq j}\t_jQ\ +\\
\sum P_{(1)}^{L\geq j}*{\ov Q}_{(1)}\ot P_{(2)}^{L\geq j}\t_j {\ov Q}_{(2)}\  +\
\sum _{L({\mathcal P}_{(2)})\leq j} P^1*\dots *P^{s-1}*(P^s\t _{j{-}L({\mathcal P}_{(2)})}Q)
\ot {\mathcal P}_{(2)}= \\
P\t_j Q\ot 1_{\K} \ +\ \de_{L\geq j}(P)\t_j\de(Q) - \sum P_{(1)}^{L=j}* Q \ot P_{(2)}^{L=j}\ +\qquad \qquad\\
 \sum _{L({\mathcal P}_{(2)})\leq j} P^1*\dots *P^{s-1}*(P^s\t _{j{-}L({\mathcal P}_{(2)})}Q)
\ot {\mathcal P}_{(2)},\end{multline} for $1\leq j< L(P)$, and
\begin{multline} 
\Delta (P\t_{L(P)}Q) =  P\t_{L(P)} Q\ot 1_{\K}\  +\\
 \sum  P_{(1)}^{L=L(P)}\ot P_{(2)}^{L=L(P)}\t_{L(P)}Q\  +
\sum  P_{(1)}^{L=L(P)}*{\ov Q}_{(1)}\ot P_{(2)}^{L=L(P)}\t_{L(P)}{\ov Q}_{(2)}\ +\\ 
\sum  P_{(1)}^{L=L(P)}* Q \ot P_{(2)}^{L=L(P)}\ +\ 
\sum_{L({\mathcal P}_{(2)})< L(P)} P^1*\dots *P^{s-1}*(P^s\times_{L(P)-L({\mathcal P}_{(2)})}Q)\ot {\mathcal P}_{(2)}\ =\\
P\t_{L(P)} Q\ot 1_{\K} \ + \de_{L= L(P)}(P)\t_{L(P)}\de(Q)\ +\
\sum_{L({\mathcal P}_{(2)})< L(P)} P^1*\dots *P^{s-1}*(P^s\times_{L(P)-L({\mathcal P}_{(2)})}Q)\ot {\mathcal P}_{(2)},\end{multline}
for $j = L(P)$,
which ends the proof.
\end{proo}

\begin{proposition} \label{Lem:theorforprime}
Let $Q$ be a prime $m$-Dyck path and $P$ any $m$-Dyck path. For any $1\leq i\leq m$, the coproduct satisfies that:
$$\Delta(P*_i Q)= \sum P_{(1)}*Q_{(1)}\ot P_{(2)}*_i Q_{(2)}.$$
\end{proposition}
\medskip

\begin{proo} As $Q$ is prime, using the conventions of Notation \ref{notTamari}, we have that
${\displaystyle P*_iQ = \sum_{j = c_i(P)}^{C_i(P)} P\t _j Q}$.

For $1\leq i < m$, by Lemma \ref{lemma: the lambda set of P_2 can be obtained from the lambda set of P}, 
any $P_{(2)}$ coming from an admissible cutting of $P$
satisfies $L(P_{(2)}) \geq c_m(P) > j$, for $c_i(P)\leq j\leq C_i(P)$, which implies that 
$\sum P_{(1)}^{L=j}* Q\ot P_{(2)}^{L=j}=0$,
and
$$ {\displaystyle \sum _{L({\mathcal P}_{(2)})\leq j} P^1*\dots *P^{s-1}*(P^s\t _{j{-}L(P_{(2)})}Q)\ 
\ot {\mathcal P}_{(2)}}  =0.$$
So, $P_{(2)}$ satisfies $L(P_{(2)})\geq j$ 
(with $c_i(P)\leq j\leq C_i(P)$) if, and only if, $P_{(2)}\not=1_{\K}$.
Therefore, applying Lemma \ref{lemma: 2 formulas for the diagonal on m-Dyck paths}, we obtain
\begin{multline}\nonumber
 \de (P\t _j Q) = \sum _{L(P_{(2)})\geq j}P_{(1)}*Q_{(1)}\ \ot \ P_{(2)}\t_jQ_{(2)} +\ P\t_j Q\ \ot\ 1_{\K} =\\
 = \sum _{P_{(2)}\not =1_{\K}}P_{(1)}*Q_{(1)}\ \ot \ P_{(2)}\t_jQ_{(2)} +\ P\t_j Q\ \ot\ 1_{\K}.
\end{multline}
Applying the formula above and Lemma \ref{lemma: the lambda set of P_2 can be obtained from the lambda set of P}, we get
\begin{multline}\nonumber
 \Delta(P*_iQ) =\Delta \left(\sum_{j =c_i(P)}^{C_i(P)} P\times_jQ\right) =\\
\sum_{j =c_i(P)}^{C_i(P)}\bigl (\sum_{P_{(2)}\not =1_{\K}} P_{(1)} * Q_{(1)}\ \ot \  P_{(2)} \t_j Q_{(2)} \bigr )+ P*_iQ\ot 1_{\K}\ = \\
\sum_{P_{(2)}\not =1_{\K}} P_{(1)}*Q_{(1)}\ \ot \bigl (\sum_{j =c_i(P)}^{C_i(P)} P_{(2)}\times_j Q_{(2)}\bigr )\ +\ P*_iQ \ot 1_{\K} \ = \\
\sum_{P_{(2)}\not =1_{\K}} P_{(1)}*Q_{(1)}\ \ot \bigl (\sum_{j =c_i(P_{(2)})}^{C_i(P_{(2)})} P_{(2)}\times_j Q_{(2)}\bigr )\ +\ P*_iQ \ot 1_{\K}\ =\\
\sum_{P_{(2)} \neq 1_{\K}} P_{(1)}*Q_{(1)}\ot P_{(2)}*_i Q_{(2)}\ +\ P*_iQ \ot 1_{\K}\ =\\
\sum P_{(1)}*Q_{(1)}\ \ot\  P_{(2)}*_iQ_{(2)}.\end{multline}

To prove the formula for $i=m$, we use that
$$\de(P*_mQ)= \sum_{j = c_m}^{L(P)-1}\de(P\times_jQ)\ +\ \de(P\times_{L(P)}Q).$$

Applying Lemma \ref{lemma: 2 formulas for the diagonal on m-Dyck paths}, to both terms of the previous equality, we obtain
\begin{multline}\nonumber
 \sum_{j = c_m(P)}^{L(P)-1}\de(P\times_jQ)\ =\ 
 \sum_{j = c_m(P)}^{L(P)-1} \bigl( P\t_j Q\ \ot\  1_{\K} \  + \de_{L\geq j}(P)\t_j\de(Q) \ +\\ 
 \sum _{L({\mathcal P}_{(2)})\leq j} P^1*\dots *P^{s-1}*(P^s\t _{j{-}L({\mathcal P}_{(2)})}Q) \ot {\mathcal P}_{(2)} \ - 
 \sum P_{(1)}^{L=j}* Q\ \ot \ P_{(2)}^{L=j} \bigr),  \end{multline}
and 
\begin{multline}\nonumber
 \de(P\times_{L(P)}Q)\ =\  P\t_{L(P)} Q\ \ot\  1_{\K} \  + \de_{L=L(P)}(P)\t_{L(P)}\de(Q) \ +\\ 
\sum _{L({\mathcal P}_{(2)}) < L(P)} P^1*\dots *P^{s-1}*(P^s\t _{L(P){-}L({\mathcal P}_{(2)})}Q) \ot {\mathcal P}_{(2)}.
  \end{multline}
  
 We leave the proof of the following two equalities to the reader, from which the proof of the case $i=m$ is complete:
$$\Delta(P)*_m\Delta(Q) = \sum_{j = c_m(P)}^{L(P)} \bigl ( P\t_j Q\ \ot\  1_{\K} \  + \de_{L\geq j}(P)\t_j\de(Q) \bigr),$$
and
\begin{multline}\nonumber
\sum_{j = c_m(P)}^{L(P)-1}\bigl(  \sum P_{(1)}^{L=j}* Q\ \ot \ P_{(2)}^{L=j}\bigr)\ =\\
\sum_{j = c_m(P)}^{L(P)}  \sum _{L({\mathcal P}_{(2)})\leq min\{ j,L(P)-1\}} \! \! \! \! \! \!\!\!\!\!\!\!\!\! P^1*\dots *P^{s-1}*(P^s\t _{j{-}L({\mathcal P}_{(2)})}Q) \ot {\mathcal P}_{(2)}.
\end{multline}
\end{proo}
\medskip

 We may prove now Theorem \ref{theorem: diagonal respect the *_i products}.

\begin{proo} {\bf of Theorem \ref{theorem: diagonal respect the *_i products}}
We prove the result applying a recursive argument on the number of prime factors of $Q$. Suppose $Q=Q_1\t_0\dots\t_0Q_r$, 
where the $Q_i$s are prime Dyck-paths. 

For $r=1$, the result is proved in Proposition \ref{Lem:theorforprime}. 

Suppose that $r >1$, and let $R := Q_1\times_0\dots\times_0 Q_{r-1}$. Applying the relations satisfied by the products $*_i$'s, we get that:
$$P*_iQ = P*_i(R*_0Q_r) =\sum_{j=i}^m(P*_jR)*_iQ_r-\sum_{j=1}^i P*_i(R*_jQ_r).$$

Since $Q_r$ is prime and $R$ is the product of $r{-}1$ prime factors, the recursive hypothesis states that 
$$\Delta((P*_jR)*_iQ_r)=(\Delta(P)*_j\Delta(R))*_i \Delta(Q_r).$$

By Lemma \ref{lemma2}, if $j\geq 1$, then the element $R*_jQ_r$ has less than $r$ prime factors. So, we have that:
$$\Delta(P*_i(R*_jQ_r)) = \Delta(P)*_i(\Delta(R)*_j\Delta(Q)_r).$$

Since $\Dot2$ is a $\Dy ^m$ algebra, the substraction of these two terms gives exactly $\Delta(P)*_i(\Delta(R)*_0\Delta(Q_r))$ which is 
equal to $\Delta(P)*_i\Delta(Q)$, which ends the proof of the theorem. 
\end{proo}
\medskip

\begin{corollary}
The coproduct $\Delta$ (hence also $\overline{\de}$) is coassociative.
\end{corollary}
\medskip

\begin{proo} We need to show that the composition
$$(\Delta \ot {\mbox {Id}}{-}{\mbox {Id}}\ot \Delta)\circ \Delta: \D_m\longrightarrow \ov{\D_m^+\ot\D_m^+\ot\D_m^+},$$
is zero.

There is a $\Dy^m$-algebra structure on $\ov{\D_m^+\ot\D_m^+\ot\D_m^+}$ given by:
 $$(x_1\ot x_2\ot x_3)*_i(y_1\ot y_2\ot y_3)=(x_1*y_1)\ot (x_2*y_2)\ot (x_3*_iy_3),$$
and we make similar considerations as in the case of $\Dot2$ when $x_3=y_3=1$. 

As $\Delta$ is a $\Dy ^m$ homomorphism, it is easy to see that both $\Delta\ot 1, 1\ot \Delta$ are so. 

By Theorem \ref{theorem: Dm is free on one generator}, coassociativity of $\Delta$ follows from the fact that
$$(\Delta\ot {\mbox {Id}}-{\mbox {Id}}\ot \Delta)(\Delta \rho_m) = 0,$$
on the generator $\rho_m$ of $\D_m$.
\end{proo}

\end{document}